%% file: Carleson_submitMay2015.tex
\documentclass[oneside]{amsart}
\usepackage{amssymb,latexsym,amsmath,amsthm,enumerate,mathrsfs}
\usepackage{fancyhdr}
\usepackage{color}
\usepackage{stmaryrd}

\input{format}

	\newcommand{\xtra}[1]{}
	
	\newcommand{\cm}[1]{{\color{blue}{#1}}}

\numberwithin{equation}{section}

\newcommand{\uDel}{\underline{\Delta}}

\newcommand{\Kcal}{\mathcal{K}}
\newcommand{\Rcal}{\mathcal{R}}

\newcommand{\Rad}{\mathrm{Rad}}
\newcommand{\Ball}{\mathrm{Ball}}

\newcommand{\odd}{\mathrm{odd}}
\newcommand{\even}{\mathrm{even}}

\newcommand{\Del}{\Delta}

\newcommand{\Pscr}{\mathscr{P}}

\begin{document}

\title{A polynomial Carleson operator along the paraboloid}
\author{L. B. Pierce}
\address{Department of Mathematics, Duke University, 120 Science Drive, Durham NC 27708; Hausdorff Center for Mathematics, 62 Endenicher Allee, Bonn, Germany 53115}
\email{pierce@math.duke.edu; pierce@math.uni-bonn.de}

\author{P.-L. Yung}
\address{Department of Mathematics, The Chinese University of Hong Kong, Shatin, Hong Kong}
\email{plyung@math.cuhk.edu.hk}

\keywords{Carleson operator, Radon transform, oscillatory integrals, square function, van der Corput estimates}
\subjclass[2010]{  42B20, 43A50, 42B25, 44A12}

\date{}

\begin{abstract}
In this work we extend consideration of the polynomial Carleson operator to the setting of a Radon transform acting along the paraboloid in $\R^{n+1}$ for $n \geq 2$. Inspired by work of Stein and Wainger on the original polynomial Carleson operator, we develop a method to treat polynomial Carleson operators along the paraboloid via van der Corput estimates. A key new step in the approach of this paper is to approximate a related maximal oscillatory integral operator along the paraboloid by a smoother operator, which we accomplish via a Littlewood-Paley decomposition and the use of a square function. The most technical aspect then arises in the derivation of bounds for oscillatory integrals involving integration over lower-dimensional sets. The final theorem applies to polynomial Carleson operators with phase belonging to a certain restricted class of polynomials with no linear terms and whose homogeneous quadratic part is not a constant multiple of the defining function $|y|^2$ of the paraboloid in $\R^{n+1}$. \end{abstract}

\maketitle 


\section{Introduction}

A celebrated theorem of Carleson \cite{Car} proves an $L^2$ bound for the operator
\beq\label{Cop} f \maps \sup_{\lam \in \R} |T_\lam f(x)| ,\eeq
where for each $\lam \in \R$,
\[ T_\lam f (x) = \mathrm{p.v.} \int_{\mathbb{T}} f(x-y) e^{i \lam y } \frac{dy}{y};\]
here by convention $\mathbb{T} = [-\pi,\pi]$.
Carleson's work answered in the affirmative a long-standing question of Luzin on whether the Fourier series of an $L^2$ function must converge pointwise almost everywhere. Carleson's remarkable result inspired many generalizations: soon after, Hunt \cite{Hun68} proved $L^p$ bounds for $1<p<\infty$ for the Carleson operator (\ref{Cop}), and Sj\"{o}lin \cite{Sjo71} introduced the Carleson operator to a higher dimensional setting by defining
\beq\label{Sjolin_op}
 T_\lam f (x) =\int_{\mathbb{T}^n} f(x-y) e^{i \lam \cdot y } K(y) dy,
 \eeq
for $\mathbb{T}^n = [-\pi,\pi]^n$ and an appropriate class of Calder\'{o}n-Zygmund kernels $K$. Sj\"{o}lin's work proved that
\[ f \maps \sup_{\lam \in \R^n} |T_\lam f(x)| \]
is a bounded operator on $L^p(\mathbb{T}^n)$ for all $1<p<\infty$.
Further landmark approaches to the Carleson operator were then developed by Fefferman \cite{Fef73} and Lacey and Thiele \cite{LacThi}, and all together these works motivated the development of time-frequency analysis into a rich and active field of research.

This paper is inspired by a question of E. M. Stein, who asked if $L^p$ bounds continue to hold when the linear phase $\lam \cdot y$ in (\ref{Sjolin_op}) is replaced by a real-valued polynomial on $\R^n$ of the form
\beq\label{P_phase1}
\sum_{1 \leq |\al| \leq d} \lam_\al y^\al,
\eeq
of fixed degree $d$. If $\lam = (\lam_{\al})_{1 \leq |\al| \leq d}$ is the set of real coefficients in (\ref{P_phase1}), we denote this polynomial by $P_{\lam}(y)$ and define
 \beq\label{Tlam0}
 	 T_\lam f (x) =\int_{\R^n} f(x-y) e^{i P_\lam(y) } K(y) dy, \eeq
where $K$ is again  an appropriate Calder\'{o}n-Zygmund kernel.
Stein asked if one can show that
 \beq\label{goal} \| \sup_\lam |T_\lam f(x)| \|_{L^p(\R^n)} \leq A \|f\|_{L^p(\R^n)},
 \eeq
 for all $1<p<\infty$,
 where the supremum is now taken over all sets of coefficients $\lam = (\lam_\al)_{1 \leq |\al| \leq d} $ with each $\lam_\al$ ranging over $\R$. 
This operator, now known as the polynomial Carleson operator, remains mysterious in many cases.

A convenient way to formulate the supremum in (\ref{goal}) is to define a stopping-time function $\lam(x)$, which is taken to be any measurable function mapping $\R^n$ to the space of real coefficients for $P_\lam(y)$. Then the desired bound (\ref{goal}) would be implied by a bound of the form
\[ \| T_{\lam(x)} f(x) \|_{L^p(\R^n)}  \leq A \|f\|_{L^p(\R^n)} ,\]
in which the norm $A$ is independent of the choice of the stopping-time function $\lam(x)$. 

Along these lines, it was first shown by Ricci and Stein \cite{RSI} that for any fixed polynomial $P(x,y)$ the operator 
\[ Tf (x) = \int_{\R^n}f(x-y) e^{iP(x,y)}K(y)dy \]
is bounded on $L^p(\R^n)$ for $1<p<\infty$, with norm dependent only on the degree of $P(x,y)$ and not on the coefficients; this is the case of a polynomial stopping-time function, and hence does not imply the full inequality (\ref{goal}).  In another direction,  Stein \cite{Stein95} considered the specific case of $\R^1$ with purely quadratic phase polynomial $P_\lam(y) =\lam y^2$ and the corresponding operator
\beq\label{T_P_op}
 T_\lam f(x) = \mathrm{p.v.} \int_{\R} f(x-y) e^{i P_\lam(y)} \frac{dy}{y},
 \eeq
ultimately proving that
\[ f \mapsto \sup_{\lam \in \R} |T_\lam f(x)| \]
is bounded on $L^p(\R)$ for $1<p<\infty$. The methods used in this case hinged upon an asymptotic for the Fourier transform of the kernel $e^{i\lam y^2}/y$, which is not easily generalizable to higher powers in the phase, or to higher dimensions. 

In 2001, Stein and Wainger \cite{SWCarl} introduced a new approach to polynomial Carleson operators, based on van der Corput estimates for oscillatory integrals, as well as a Kolmogorov-Seliverstov stopping-time argument. Their result is as follows:
\begin{letterthm}\label{thm_SW}
Consider the operators $T_{\lam}$ as defined in (\ref{Tlam0}), where $P_{\lam}(y)$ is a real-valued polynomial with no linear terms, of the form
\beq\label{P_no_linear}
 P_\lam(y) = \sum_{2 \leq |\al| \leq d} \lam_\al y^\al.
 \eeq
Then 
\[  \| \sup_\lam |T_\lam f(x)| \|_{L^p(\R^n)} \leq A \|f\|_{L^p(\R^n)}, \]
for $1<p<\infty$, where the supremum is over all coefficients $\lam = (\lam_\alpha)_{2 \leq |\alpha| \leq d}$ of $P_\lam(y)$.
\end{letterthm}
The restriction that $P_\lam(y)$ omits linear terms is inherent to methods of van der Corput type, as will be seen explicitly later. 

Until very recently, bounds for the full polynomial Carleson operator remained unproved. Now, the case of the polynomial Carleson operator in dimension one has been resolved by Lie \cite{Lie11X} \cite{Lie09}, who proved by time-frequency analysis methods that for
\[ T_\lam f(x) = \mathrm{p.v.} \int_{\mathbb{T}} f(x-y) e^{i P_\lam(y)} \frac{dy}{y}\]
with $P_\lam(y)$ a phase polynomial including both linear and higher order terms,
the corresponding Carleson operator is a bounded operator on $L^p$ for $1<p<\infty$ with
\[ || \sup_{\lam} |T_\lam f| \|_{L^p(\mathbb{T})} \leq A \|f \|_{L^p (\mathbb{T})} .\]
The case of the polynomial Carleson operator involving both linear and higher order terms in dimensions $n \geq 2$ remains open.

\subsection{A new setting}
In the work of this paper, we take the polynomial Carleson operator in a new direction by introducing Radon-type behavior. 
Suppose $\{P_\lam(y)\}$ is a family of real-valued polynomials of $y \in \R^n$ indexed by a parameter $\lam$. We assume all polynomials $P_{\lam}$ in the family are of degree at most $d$ for a fixed positive integer $d$. For each $\lam$ in the parameter space, we define an operator $T_\lam$, initially acting on functions $f$ of Schwartz class, by
\beq\label{T_lam_par}
 T_\lam f (x,t)= \int_{\R^n} f(x-y, t-|y|^2)e^{i P_\lam(y)}K(y)dy.
 \eeq
Thus $T_\lam$ integrates $f$ along the paraboloid $(y,|y|^2) \subset \R^{n+1}$ against an oscillatory factor and a singular kernel. 
The kernel $K$ is a Calder\'{o}n-Zygmund kernel, that is, a tempered distribution agreeing with a $C^1$ function $K(x)$ for $x \neq 0$, such that it satisfies the differential inequalities
\beq\label{K_prop}
 | \partial_x^\al K(x)| \leq A |x|^{-n-|\al|} \quad \text{for $0 \leq |\al| \leq 1$, }
\eeq
and such that $\hat{K}$ is an $L^\infty$ function.

We then consider the polynomial Carleson operator of Radon type defined by 
\beq\label{op_Carl_par}
 f \mapsto \sup_\lam |T_\lam f(x,t)|,
 \eeq
in which the supremum is taken over all $\lam$ in a suitable parameter space. When the parameter space consists of all  possible coefficients $\lam = (\lam_{\alpha})_{1 \leq |\alpha| \leq d}$ of real-valued polynomials $P_\lam(y)$ of the form (\ref{P_phase1}), it is reasonable to expect that the following $\emph{a priori}$ estimate should hold for $f$ of Schwartz class:
\[ \| \sup_{\lam} |T_\lam f | \|_{L^p(\R^{n+1})} \leq A \| f \|_{L^p(\R^{n+1})},\]
for any $1<p< \infty$. 
Proof of a result of this type is expected to require inputs from multi-dimensional time-frequency analysis, and seems overly ambitious at this time.

Therefore we restrict ourselves to suprema over a narrower class of polynomial phases, and develop an approach for treating the operator (\ref{op_Carl_par}) with van der Corput estimates. Such methods are inspired by the work of Stein and Wainger in \cite{SWCarl}, and inherently require that the phase polynomial lack terms of certain lower orders. In the case of the operator (\ref{op_Carl_par}) on the paraboloid, a first guess might be that van der Corput methods would require $P_\lam(y)$ to lack both linear and quadratic terms; this may be seen intuitively by considering a model operator
\beq\label{Rop}
 R f(x,t) = \int_{\R^n} f(x-y,t-|y|^2) \eta(y)e^{i P_\lam(y)}dy,
 \eeq
where $\eta$ is a smooth bump function of compact support. If the parameter $\lam$ is regarded as fixed for the moment, we may compute 
the Fourier multiplier of $R$ to be
\beq\label{mop}
	m(\xi,\theta) = \int_{\R^n} \eta(y)e^{i P_\lam(y) -2\pi i \xi \cdot y - 2\pi i \theta |y|^2} dy.
	\eeq
One may obtain an upper bound for $m(\xi,\theta)$ by applying van der Corput bounds to the oscillatory integral (\ref{mop}), as long as one has lower bounds for the coefficients of the phase $P_\lam(y) - 2\pi \xi \cdot y - 2\pi \theta|y|^2$ as a polynomial in $y$; one would expect it to be difficult to obtain such lower bounds if certain linear or quadratic terms are present in $P_{\lam}$.  Of course, if $\lam = \lam(x,t)$ is a stopping-time function, we cannot compute a Fourier multiplier, so this is merely a heuristic.

Even considering the family $\{P_\lam(y)\}$ of all phase polynomials that vanish to order at least 2 at the origin is still overly ambitious; the presence of the Radon transform in the operator poses further limitations on the number of degrees of freedom our method can allow in the family of polynomial phases. On the other hand, we will show that with a sufficiently sharp scalpel in hand, one may allow certain quadratic terms, as long as the degree 2 homogeneous portion of the phase polynomial is not precisely a nonzero multiple of $|y|^2$. 

\subsection{Statement of results}
Our main result is the following:

\begin{thm}\label{thm_Carleson}
Fix a dimension $n \geq 2$ and a degree $d \geq 2$. Let 
\[\Pscr = \{ p_2(y), p_3(y),\ldots, p_d(y)\}\]
 be a set of real-valued polynomials on $\mathbb{R}^n$, where each $p_j(y)$ is homogeneous of degree $j$, and $p_2(y) \ne C|y|^2$ for any nonzero constant $C$. For each $\lam = (\lam_2, \dots, \lam_d) \in \R^{d-1}$, let
\begin{equation} \label{allowable_poly}
P_{\lam}(y) = \sum_{m=2}^d \lam_m p_m(y),
\end{equation}
and let $T_\lam$ be defined as in (\ref{T_lam_par}).
Then for every $1 <p<\infty$, we have an \emph{a priori} inequality for functions $f$ of Schwartz class:
\beq\label{Tbdf}
	 \|\sup_{\lam} | T_\lam f| \; \|_{L^p(\R^{n+1})} \leq A_p \|f\|_{L^p(\R^{n+1})},
	 \eeq
where the supremum is taken over all $\lam = (\lam_2, \dots, \lam_d) \in \R^{d-1}$.
\end{thm}
As a result of the \emph{a priori} inequality, it follows from conventional limiting arguments that the Carleson operator extends to a bounded operator on $L^p$ itself, for which the same bound (\ref{Tbdf}) holds.\xtra{The key point is that for a Cauchy sequence $\{ f_j\}$ in $L^p$, 
\[ |\sup_\lam |T_\lam f_j| - \sup_{\lam} |T_\lam f_k| \; |\leq \sup_\lam |T_\lam f_j - T_\lam f_k|,\]
hence we may apply (\ref{Tbdf}) to show that $\sup_\lam |T_\lam f_j|$ is a Cauchy sequence in $L^p$.
}

At its foundation, the idea of the proof is that for any value of the coefficient parameter $\lam$ one splits the integral defining $T_\lam$ in two parts: first, a part in which $y$ is small enough relative to the coefficients of $P_\lam(y)$ that $e^{iP_\lam(y)}$ does not contribute significant oscillation; this part is compared to a maximal truncated singular Radon transform, which may be bounded on $L^p$ (see the companion paper \cite{PY15b}). The second part of $T_\lam$ comprises the region in which $y$ is large enough relative to the coefficients of $P_\lam(y)$  that $e^{i P_\lam(y)}$ exhibits significant oscillation; this part is controlled via an auxiliary maximal oscillatory Radon transform, which is interesting in its own right. Precisely, our second main result is as follows:


\begin{thm}\label{Osc_Thm}
Fix a dimension $n \geq 2$ and a degree $d \geq 2$. Let 
\[\Pscr = \{ p_2(y), p_3(y),\ldots, p_d(y)\}\]
 be a set of real-valued polynomials on $\mathbb{R}^n$, where each $p_j(y)$ is homogeneous of degree $j$, and $p_2(y) \ne C|y|^2$ for any nonzero constant $C$. Let 
\beq \label{eq:Lambdadef} 
\Lambda = \Lambda (\Pscr) = \{ 2 \leq m \leq d \colon p_m(y) \text{ is not identically zero}\}.
\eeq
 For $\lam = (\lam_2,\dots,\lam_d) \in \R^{d-1}$, let $P_\lam(y)$ be as in (\ref{allowable_poly}), and let 
\beq \label{lam_Rn_def} \| \lam \| = \sum_{m \in \Lambda} |\lam_m|.\eeq
Then one can define an operator $\;^{(\eta)}I_a^\lam$ acting on $f$ of Schwartz class by
 \beq\label{I_op_dfn_0}
    \;^{(\eta)}I_a^\lam f(x,t) = \int_{\R^n}  f(x-y,t-|y|^2) e^{iP_{\lam}(y/a)}\frac{1}{a^n} \eta(\frac{y}{a}) dy,
   \eeq
for $\eta$ a $C^1$ bump function supported in the unit ball on $\R^n$, $\lam \in \R^{d-1}$, and $a > 0$.

Suppose furthermore that 
$\{\eta_{k}\}_{k \in \mathbb{Z}}$ is a family of bump functions supported in the unit ball $B_1(\R^n)$ with $C^1$ norm uniformly bounded by 1. Then there exists a fixed $\del>0$ such that for any Schwartz function $f$ on $\R^{n+1}$ and any $r \geq 1$, 
\beq\label{osc_decay_parabolic}
\| \sup_{k \in \Z} \sup_{\bstack{\lam \in \R^{d-1}}{r \leq \|\lam\| < 2 r}} |^{(\eta_{k})}I^\lam_{2^k} f(x,t)|\, \|_{L^2(\R^{n+1})} \leq A r^{-\del}  \|f\|_{L^2(\R^{n+1})}.
\eeq
\end{thm}


In order to put this result in context, it is illustrative to note that for fixed $k$ and $\lam$, it is simple to bound $^{(\eta_k)} I_{2^k}^\lam$ on $L^2$, since it has a bounded Fourier multiplier. Even after taking suprema over the parameters $k$ and $\lam$, the operator in question is bounded pointwise by the classical maximal Radon transform along the paraboloid, which acts on Schwartz functions $f$ by
\beq\label{MRad_dfn}
 \Mcal_{\Rad}f(x,t) = \sup_{a>0} \int_{\R^n} |f|(x-u,t-|u|^2) \frac{1}{a^n}\chi_{B_1}(\frac{u}{a}) du.
 \eeq
This is well-known to be bounded on $L^p$ for $1<p \leq \infty$ when acting on functions of Schwartz class (see for instance Chapter 11 of \cite{SteinHA}).
Thus the key feature of Theorem \ref{Osc_Thm} is the decay $r^{-\del}$ of the norm of the operator as $r \to \infty$. We remark that it will be clear from the proof that Theorem \ref{Osc_Thm} continues to hold for all $r \geq c$ for any fixed constant $c>0$.

As summarized thus far, the basic structure of reducing the treatment of the Carleson operator to the family of oscillatory integral operators $^{(\eta_k)} I_{2^k}^\lam$  is inspired by the work of Stein and Wainger in \cite{SWCarl}. Their original approach was to apply a $TT^*$ argument to $^{(\eta_k)}I^\lam_{2^k}$ and then apply van der Corput bounds to the kernel of $TT^*$ in order to show that $TT^*$ could be majorized by certain maximal functions for which $L^2$ bounds were known. In our setting of Radon transforms, this approach is not sufficient, and in order to prove Theorem \ref{Osc_Thm}, we must diverge significantly from the original work of Stein and Wainger and develop a new strategy that allows us to control the presence of Radon-type behavior.

To motivate our new approach, we recall that the classical treatment of the maximal Radon transform on the paraboloid introduces a smoother version of the operator, as well as an accompanying square function. In our setting,  we introduce a smoothed version of $^{(\eta)}I_a^\lam$ defined by
\beq\label{J0}
     ^{(\eta)}J_a^\lam f(x,t)  = \iint_{\R^{n+1}}  f(x-y,t-z) e^{iP_{\lam}(y/a)} \frac{1}{a^n} \eta(\frac{y}{a}) \frac{1}{a^2} \zeta(\frac{z}{a^2})dydz,
     \eeq
where $\zeta$ is a smooth function of compact support on $\R$, chosen such that $ \int \zeta (s)ds=1$.
We also utilize a square function defined (roughly speaking) by 
\beq\label{Sf0}
  S_r(f)(x,t) = \left[ \sum_{k=-\infty}^\infty \left( \sup_{r \leq \|\lam\| <2r} |(^{(\eta_{k})}I^\lam_{2^k} - {^{(\eta_{k})}J^\lam_{2^k}}) f(x,t)|\right)^2 \right]^{1/2}.
  \eeq
 (The precise definition of $S_r(f)$ is given in equation (\ref{Sfn_precise}).)
Our result in Theorem \ref{thm_Sf} is  that $S_r(f)$ is bounded on $L^2$, with a norm that decays as $r \to \infty$; this is the heart of proving Theorem \ref{Osc_Thm}.

\subsection{Outline of the paper}
Our approach succeeds by carefully intertwining the role of the square function, $TT^*$ methods, and van der Corput estimates. This work ultimately entails significant technical details, thus Section \ref{sec_anatomy} outlines the motivation for our approach, explaining why complications arise in the Radon case, and demonstrating why the tools we introduce allow us to circumvent these complications.
We then assemble in Sections  \ref{sec_prelim_VdC} and \ref{sec_LP} certain preliminary results on van der Corput estimates and a Littlewood-Paley decomposition.  In Section \ref{sec_reduction} we reduce the proof of Theorem \ref{thm_Carleson} to Theorem \ref{Osc_Thm}. In Section \ref{sec_outline_pf}, we define the relevant square function and reduce Theorem \ref{Osc_Thm} to an $L^2$ bound for the square function. In Section \ref{sec_diff_op} we prove Theorem \ref{Diff_Osc_Thm}, the main result that implies the boundedness of the square function, except for certain van der Corput estimates. These novel Van der Corput estimates are the most technical part of the paper, and for the sake of the reader we initially treat several examples in the case of dimension $n=2$ in Section \ref{sec_vdC_n2}, reserving the fully general higher dimensional case for Section \ref{sec_high_dim}. Finally in Section \ref{sec_gen_dim} we complete the proof of several key propositions in general dimension $n \geq 2$. An appendix in Section \ref{sec_appendix} provides a proof of an auxiliary theorem on maximal truncated singular Radon transforms (Theorem \ref{thm_trunc_Radon}).
 Readers interested in a first view of the key ideas may read Section \ref{sec_anatomy} for a motivation of our approach, and then focus on Sections \ref{sec_reduction}, \ref{sec_outline_pf}, \ref{sec_diff_op}, and \ref{sec_vdC_n2}.

\section{Anatomy of the proof}\label{sec_anatomy}

In the work of Stein and Wainger \cite{SWCarl}, the key oscillatory integral operator  takes the form 
\[ I_a^\lam f(x) = \int_{\R^n} f(x-y) e^{i P_\lam (y/a)} \frac{1}{a^n} \eta (\frac{y}{a}) dy,\]
with $\eta$ a $C^1$ bump function supported in the unit ball and $P_\lam(y)$ a polynomial phase of the form (\ref{P_no_linear}), i.e. with no linear terms.
Their analogue of Theorem \ref{Osc_Thm} is obtained by using stopping-times to linearize the maximal aspect of the operator $ f\mapsto \sup_{\lam,a} |I_a^\lam f(x)|$, so that this last operator is represented by
\[T: f \mapsto I_{a(x)}^{\lam(x)} f(x),\]
where $a(x)$ and $\lam(x)$ are arbitrary measurable stopping-time functions taking values of the form $2^k$ for $k \in \Z$ and $r \leq \| \lam \| <2r$, respectively. 
Stein and Wainger prove that
if we define a kernel $K^{\lam_1,\lam_2}_{a_1,a_2}$ by
$$
K^{\lam_1,\lam_2}_{a_1,a_2}(y) = \frac{1}{a_1^n}\frac{1}{a_2^n} \int_{\mathbb{R}^n} e^{iP_{\lam_1}(\frac{y+z}{a_1})-iP_{\lam_2}(\frac{z}{a_2})} \eta(\frac{y+z}{a_1}) \eta(\frac{z}{a_2}) dz,
$$
then $TT^*f(x)$ is given by
$$
TT^*f(x) = \int_{\mathbb{R}^n} f(x-y) K^{\lam(x),\lam(x-y)}_{a(x),a(x-y)}(y) dy.
$$ 
They then proceed to show, via a van der Corput estimate, that $K^{\lam_1,\lam_2}_{a_1,a_2}$ satisfies the following pointwise bound:
\beq\label{K_orig_bd}
 |K^{\lam_1,\lam_2}_{a_1,a_2}(x)| \leq C r^{-\del} \left( \frac{1}{a_1^n} \chi_{B_2}( \frac{x}{a_1}) + \frac{1}{a_2^n} \chi_{B_2}( \frac{x}{a_2} ) \right) 
 +
  \left( \frac{1}{a_1^n} \chi_{E_{\lam_1}}( \frac{x}{a_1}) + \frac{1}{a_2^n} \chi_{E_{\lam_2}}( \frac{x}{a_2} ) \right) ,
  \eeq
in which $B_2 = B_2(\R^n)$ is the ball of radius 2 and $E_{\lam_1},E_{\lam_2}$ are certain exceptional subsets of $B_2(\R^n)$ with small measure, namely $|E_{\lam_1}|, |E_{\lam_2}| \leq r^{-\del}$. 
 An $(L^2,L^2)$ norm for $TT^*$ with appropriate decay in $r$ can then be obtained by comparison to certain maximal functions, via a clever bilinear argument.

A straightforward generalization of this argument appears to fail in the setting of Theorem \ref{Osc_Thm}. Indeed, it is known that stopping-time arguments do not work well when Radon transforms are involved, and in order to illustrate this break-down and to motivate how we proceed instead, we briefly consider the following simple model of a maximal Radon transform (without oscillatory factor) along a paraboloid in $\R^3$:
\beq\label{Max_illus}
\mathcal{M}f(x,t) = \sup_{a > 0} \int_{y \in \mathbb{R}^2} |f|(x-y,t-|y|^2) \frac{1}{a^2} \eta(\frac{y}{a}) dy.
\eeq
Here $f$ is a function of $(x,t) \in \mathbb{R}^2 \times \mathbb{R}$, and $\eta$ is a fixed bump function supported in the unit ball in $\mathbb{R}^2$ such that $\eta=1$ in a neighborhood of the origin. This operator is of course known to be bounded on $L^2(\mathbb{R}^3)$ by standard results about maximal Radon transforms along a paraboloid. However, for the sake of illustration we will attempt an alternative proof using stopping-times and $TT^*$, and we will see what goes wrong. This will make it clear how this approach would fail for the operator considered in Theorem \ref{Osc_Thm}.

To begin our (doomed) alternative approach, we fix a non-negative function $f$ and pick any measurable stopping-time $a(x,t)$ taking positive real values, and consider the linear operator 
$$
Tf(x,t) = \int_{y \in \mathbb{R}^2} f(x-y,t-|y|^2) \frac{1}{a(x,t)^2} \eta(\frac{y}{a(x,t)}) dy;
$$
then $TT^*f(x,t)$ can be written as
\beq \label{simpleTT*}
TT^*f(x,t) = \int_{y \in \mathbb{R}^2} \int_{z \in \mathbb{R}^2} f(x-y+z,t-|y|^2+|z|^2) \frac{1}{a_1^2} \eta(\frac{y}{a_1}) \frac{1}{a_2^2} \eta(\frac{z}{a_2}) dy dz,
\eeq
where $a_1$, $a_2$ stand for the stopping-time functions
$$
a_1 := a(x,t), \quad a_2:=a(x-y+z,t-|y|^2+|z|^2).
$$
In the region in the $(y,z)$-plane where $a_1 \geq a_2$ we would sequentially make the change of variables $(y,z) \mapsto (u,\tau,\sigma)$, where 
$$
u = y-z, \qquad  \tau = \frac{u \cdot z}{|u|} = \frac{u_1z_1 + u_2z_2}{|u|}, \qquad \sig=  \frac{u_2z_1 - u_1z_2}{|u|}.
$$
\xtra{Here we are using the fact that $a_2 \leq a_1$ so that when we know $|z| \leq a_2$ and $|u+z| \leq a_1$ then we know $|u| \leq a_1 + a_2 \leq 2a_1$. } 
On the other hand, in the region in the $(y,z)$-plane where $a_1 \leq a_2$ we would make the change of variables $(y,z) \mapsto (u,\tau,\sigma)$, where
$$
u = z-y,  \qquad \tau = \frac{u \cdot y}{|u|} = \frac{u_1y_1 + u_2y_2}{|u|}, \qquad \sig=  \frac{u_2y_1 - u_1y_2}{|u|}.
$$
After this process, and a trivial integration with respect to $\sigma$, we would have
\begin{align*}
|TT^*f(x,t)| \lesssim
& \int_{u \in \mathbb{R}^2} \int_{\tau \in \mathbb{R}} f(x-u,t-|u|^2-2|u|\tau) \frac{1}{b_1^2} \chi_{B_2}(\frac{u}{b_1}) \frac{1}{b_2} \chi_{B_2}(\frac{\tau}{b_2}) du d\tau \\
 & \quad + \int_{u \in \mathbb{R}^2} \int_{\tau \in \mathbb{R}} f(x+u,t+|u|^2+2|u|\tau) \frac{1}{b_3^2} \chi_{B_2}(\frac{u}{b_3}) \frac{1}{b_4} \chi_{B_2}(\frac{\tau}{b_4}) du d\tau
\end{align*}
where 
$$
b_1:=a(x,t), \qquad b_2:=a(x-u,t-|u|^2-2|u|\tau),
$$
$$
b_3:=a(x+u,t+|u|^2+2|u|\tau), \qquad b_4:=a(x,t).
$$
The two terms on the right hand side above are similar; for the sake of illustration we will focus on the first term, which we temporarily denote by $I(x,t)$. Let $\mathcal{M}_{\Rad}$ denote the usual maximal Radon transform along the paraboloid (defined in (\ref{MRad_dfn})) and let $\mathcal{M}^{(3)}$ denote a one-dimensional maximal function in the $t$-variable alone. Then one could hope that the term $I(x,t)$ is bounded by
$$
|I(x,t)| \leq C \int_{u \in \mathbb{R}^2} \mathcal{M}^{(3)} f(x-u,t-|u|^2) \frac{1}{b_1^2} \chi_{B_2}(\frac{u}{b_1}) du \leq C \mathcal{M}_{\Rad} \mathcal{M}^{(3)} f (x,t).
$$
Unfortunately the first inequality above fails: the stopping-time function $b_2 = a(x-u,t-|u|^2-2|u|\tau)$ depends not only on $x$, $t$ and $u$, but also on $\tau$. This leads to a first crucial failure: indeed, if $b := b(\tau)$ is an unknown function of $\tau$, then 
\[\int_{\tau \in \mathbb{R}} F(t-\tau) \frac{1}{b} \chi_{B_1} (\frac{\tau}{b}) d\tau\]
 is not necessarily bounded by the Hardy-Littlewood maximal function of $F$. In such situations, it can be advantageous to dualize the statements, and integrate $I(x,t)$ against a test function $g(x,t)$. Applying this idea and then changing variables $t \mapsto t+2|u|\tau$ leads to 
 \begin{multline*}
\int_{\mathbb{R}^3} I(x,t) g(x,t) dx dt 
\\= \int_{(x,t) \in \mathbb{R}^3} \int_{u \in \mathbb{R}^2} \int_{\tau \in \mathbb{R}}  f(x-u,t-|u|^2) \frac{1}{c_1^2} \chi_{B_2}(\frac{u}{c_1}) \frac{1}{c_2} \chi_{B_2}(\frac{\tau}{c_2}) g(x,t+2|u|\tau) du d\tau dx dt,
\end{multline*}
where now the stopping-time parameters are
$$
c_1 := a(x,t+2|u|\tau), \qquad c_2 := a(x-u,t-|u|^2).
$$
It is true now that $c_2$ is independent of $\tau$, so one would hope that the above is bounded by 
$$
\int_{(x,t) \in \mathbb{R}^3} \int_{u \in \mathbb{R}^2}   f(x-u,t-|u|^2) \frac{1}{c_1^2} \chi_{B_2}(\frac{u}{c_1})  \mathcal{M}^{(3)} g(x,t) du dx dt;
$$
unfortunately this is still not true, since now $c_1$ depends on $\tau$, so that any integration with respect to $\tau$ would also be forced to consider $c_1$.

The failure of such $TT^*$ arguments for maximal Radon transforms (even without Carleson-type behavior) demonstrates the need for a different approach. In the case of the maximal Radon transform $\Mcal f(x,t)$ considered in (\ref{Max_illus}), a well-known successful strategy is to construct a smoother variant of the average over the paraboloid, and to compare the maximal Radon transform with this smoother maximal function via a square function. (One may see this for example in \cite{SWCurv} and \S 1.2 of Chapter 11 in \cite{SteinHA}.) More precisely, the maximal Radon transform $\mathcal{M}f(x,t)$  is  comparable (up to a constant) to $\sup_{k \in \mathbb{Z}} A_k f(x,t)$ where
$$
A_kf(x,t) := \int_{\mathbb{R}^2} f(x-y,t-|y|^2) 2^{-2k} \eta(2^{-k}y) dy.
$$
One then compares it to a smoother variant, namely $\sup_{k \in \mathbb{Z}} B_kf(x,t)$, where
$$
B_kf(x,t):= \int_{\mathbb{R}^2} \int_{\mathbb{R}} f(x-y,t-s) 2^{-2k} \eta(2^{-k}y) 2^{-2k} \zeta(2^{-2k} s) ds dy,
$$
in which $\zeta$ is a smooth bump function supported on $[-1,1] \subset \mathbb{R}$, such that $\int_{\mathbb{R}} \zeta(s) ds = 1$.

The operator $\sup_{k \in \mathbb{Z}} B_k f$ is known to be bounded on $L^2$ (e.g. by covering arguments), and the difference $\sup_{k \in \mathbb{Z}} |A_kf - B_kf|$ is then trivially dominated by the square function
$$
\left( \sum_{k \in \mathbb{Z}} |A_k f - B_k f|^2 \right)^{1/2}.
$$
The boundedness of the square function on $L^2$ can be proved using methods related to the Fourier transform. (The condition $\int_{\mathbb{R}} \zeta(s) ds = 1$ is used to guarantee the vanishing of the multiplier of $A_k - B_k$ near the origin.)

We now return to the setting of this paper, and in particular the oscillatory integral operator considered in Theorem \ref{Osc_Thm}. In order  to prove Theorem~\ref{Osc_Thm}, we will use the square function $S_r$ defined in (\ref{Sf0}) to control the supremum over $k$ in the $L^2$ bound (\ref{osc_decay_parabolic}). Since morally speaking $I^{\lam}_{2^k} - J^{\lam}_{2^k}$ only `sees' a function of `frequency' $2^k$, an appropriate Littlewood-Paley decomposition will allow us to reduce to a situation in which $k$ is fixed. We then need to control the supremum over $\lam$: this is handled using a $TT^*$ argument and stopping-times (in the phase only), which ironically are now the key to success, once the square-function has been applied. Roughly speaking, within the square-function the presence of the inconvenient stopping-times normalizing the supports of the bump functions have been removed, and so $TT^*$ and stopping-times may be used in order to apply van der Corput estimates.  But now we encounter another difficulty, and it is here that we need to restrict ourselves to the special class of polynomial phases considered in Theorem \ref{Osc_Thm}, namely those that are in the span of the fixed polynomials $p_2(y), \dots, p_d(y)$ provided in the theorem hypotheses. This restriction, which appears as well in the statement of Theorem \ref{thm_Carleson}, has its origin in the proof of Theorem \ref{Osc_Thm}.

Again for the sake of illustration, let us suppose we are more ambitious, and wish to take a supremum over all real polynomial phases of degree $\leq d$ with no linear or quadratic terms. 
We fix our attention for the moment on the case of dimension $n=2$ and let 
\begin{eqnarray*}
 \mathcal{Q}_d &=& \{ \text{real polynomials of degree $\leq d$ on $\mathbb{R}^2$}\},\\
  \mathcal{Q}_d^*& = & \{\text{$Q \in \mathcal{Q}_d$ that vanishes to order $\geq 2$ at the origin}\}.
  \end{eqnarray*}
   Let $\Gamma_d^*$ be the set of all coefficients of polynomials in $\mathcal{Q}_d^*$, so that any element $\lam \in \Gamma_d^*$ takes the form $\lam = (\lam_{\al})_{3 \leq |\al| \leq d}$. For each $\lam \in \Gamma_d^*$, we may define a polynomial in $\mathcal{Q}_d^*$ by
   \[Q_{\lam}(y) = \sum_{3 \leq |\al| \leq d} \lam_{\al} y^{\al},\]
   with coefficient norm $\|\lam\| = \sum_{\al} |\lam_{\al}|$. 
We may now accordingly define $T_{\lam}f$ and $^{(\eta)}I^{\lam}_a f$ by (\ref{T_lam_par}) and (\ref{I_op_dfn_0}) respectively, with the phase $Q_{\lam}(y)$ in place of $P_\lam(y)$. 
If we now want to prove that the polynomial Carleson operator $\sup_{\lam \in \Gamma_d^*} |T_{\lam} f|$ is bounded on $L^2(\mathbb{R}^3)$, we are led to consider the maximal oscillatory operator given by
 $$\sup_{k \in \mathbb{Z}} \sup_{ \bstack{\lam \in \Gamma_d^*}{r \leq \|\lam\| \leq 2r} } |^{(\eta)}I^{\lam}_{2^k} f|,$$
for which we must show that there exists a $\del>0$ such that its $L^2$ norm is $\leq C r^{-\del}$ for all $r \geq 1$. Since a square function argument will essentially allow one to handle the supremum over $k$,  we will temporarily fix $k = 0$ and consider only the supremum over $\lam$. 

Thus we restrict our attention to trying to prove
$$
\| \sup_{ \bstack{\lam \in \Gamma_d^*}{r \leq \|\lam\| \leq 2r} } |^{(\eta)}I^{\lam}_1 f|\, \|_{L^2(\mathbb{R}^3)} \leq C r^{-\delta} \|f\|_{L^2(\mathbb{R}^3)},
$$ 
for some $\delta > 0$. Motivated by Stein and Wainger \cite{SWCarl}, a natural approach is to combine stopping-times with a $TT^*$ argument: let $\lam (x,t) = (\lam_{\al}(x,t))_{3 \leq |\al| \leq d}$ be a measurable stopping-time function taking values in $\Gamma_d^*$, and consider 
$$Tf(x,t) = {}^{(\eta)}I^{\lam(x,t)}_1 f(x,t).$$
Then $TT^* f(x,t)$ can be written as
\beq \label{modelTT*}
TT^*f(x,t) = \int_{u \in \mathbb{R}^2} \int_{\tau \in \mathbb{R}} f(x-u,t-|u|^2-2|u|\tau) K_{\sharp}^{\lam(x,t),\lam(x-u,t-|u|^2-2|u|\tau)}(u,\tau) du d\tau,
\eeq
where for $\nu$, $\mu \in \Gamma_d^*$ we define the kernel
$$
K_{\sharp}^{\nu,\mu} (u,\tau):= \int_{\sigma \in \mathbb{R}} e^{iP_{\nu}(u+z) - iP_{\mu}(z)} \eta(u+z) \eta(z) d\sigma,
$$
in which  $z=(z_1,z_2)$ is a function of $(u,\tau,\sigma)$ defined implicitly by the relations
\beq\label{tausig}
 \tau = \frac{u \cdot z}{|u|} = \frac{u_1z_1 + u_2z_2}{|u|}, \qquad \sig=  \frac{u_2z_1 - u_1z_2}{|u|}.
 \eeq
Using van der Corput estimates, we can show the existence of a small $\delta>0$ such that the following holds: for any $r 
\geq 1$, and for each $\nu$, $\mu \in \Gamma_d^*$ with $r \leq \|\nu\|, \|\mu\| \leq 2r$, there exists a ``bad'' (but small) set of $u$, denoted $G^{\mu,\nu} \subset B_2(\R^2)$, and for each $u \in B_2(\R^2)$ there exists a ``bad'' (but small) set of $\tau$, denoted $F^{\nu,\mu}_u \subset B_2 (\R)$, such that
$$
|G^{\nu,\mu}| \leq C r^{-\delta}, \qquad |F^{\nu,\mu}_u| \leq C r^{-\delta},
$$
and
\beq\label{K_sharp_illus}
|K_{\sharp}^{\nu,\mu} (u,\tau)| \leq C [ r^{-\delta} \chi_{B_2}(u) \chi_{B_2}(\tau) + \chi_{G^{\nu,\mu}}(u) \chi_{B_2}(\tau) + \chi_{B_2}(u) \chi_{F^{\nu,\mu}_u}(\tau)].
\eeq
 At first glance this looks as good as the bound (\ref{K_orig_bd}) that is the keystone of Stein and Wainger's work. Yet in fact (\ref{K_sharp_illus}) actually fails to be effective since the small sets $G^{\nu,\mu}$ and $F^{\nu,\mu}_{u}$ depend simultaneously on $\nu$ and $\mu$, and this is just as deleterious as the simultaneous appearance of $a_1$ and $a_2$ in formula (\ref{simpleTT*}), again because of the presence of Radon-type behavior. More explicitly, if we apply the estimate (\ref{K_sharp_illus}) for $K_{\sharp}^{\nu,\mu}$ in (\ref{modelTT*}), it would show that $TT^*f(x,t)$ is bounded by a sum of three terms, of which we single out the term
$$
\int_{u \in \mathbb{R}^2} \int_{\tau \in \mathbb{R}} |f|(x-u,t-|u|^2-2|u|\tau) \chi_{G^{\lam(x,t),\lam(x-u,t-|u|^2-2|u|\tau)}}(u) \chi_{B_2}(\tau) du d\tau.
$$
We would like to say that this is bounded by a concatenation of maximal functions, including Radon-type adaptations of the ``small set maximal functions'' used by Stein and Wainger, however, this is not true, since when one integrates in $\tau$, one must remember that the small set $G^{\lam(x,t),\lam(x-u,t-|u|^2-2|u|\tau)}$ depends also on $\tau$. 

We circumvent all of these difficulties by intertwining the methods of square functions, $TT^*$, and stopping-times. First, we eliminate the presence of two independent stopping-times $a_1,a_2$ by passing to the square function (\ref{Sf0}); note that inside each term of $S_r(f)$, the scaling parameters are fixed and equal: $a_1=a_2=2^k$.
(Here it is important that we are able to restrict the supremum over the scaling factors $a$ to a supremum over a countable set of scaling factors $a=2^k$ for $k \in \Z$, rather than over all $a >0$, as in the work of Stein and Wainger.)
Second, we apply a $TT^*$ argument to each fixed summand within the square function, and apply van der Corput estimates to extract decay in $r$ from the kernel of $TT^*$.  (We note that our requirement that $n\geq 2$ arises here, since in the case $n=1$, the kernel of $TT^*$ does not take the form of an integral, and hence will not admit van der Corput estimates.) In particular, we are able to construct small bad sets $G^{\nu}$ and $F^{\nu}_u$ that depend on only $\nu$ but not on $\mu$ by an argument that explicitly uses the assumption that the polynomial phase $P_\lam(y)$ belongs to the restricted class specified in Theorems \ref{thm_Carleson} and \ref{Osc_Thm}.

Finally, as we have mentioned already, in order to estimate the square function $S_r$ (in particular, to carry out the sum over $k$ in (\ref{Sf0})), we need to know that $(I^{\lam}_a - J^{\lam}_a)f$ sees only the part of $f$ that has frequency $a$. To make this precise, we need to introduce a family of Littlewood-Paley projections $\Delta_j$ and to obtain some almost orthogonality for $(I^{\lam}_a - J^{\lam}_a) \Delta_j$ if $a$ is very different from $2^j$. This turns out to be rather tricky, especially when $2^j$ is much smaller than $a$; in that case one would like to write  $\Delta_j$ as a derivative, and integrate by parts, but then a derivative may land on the singular kernel of $I^{\lam}_a$. Since $I^{\lam}_a$ is a Radon transform that involves an integration over a lower dimensional submanifold, if the derivative is not tangential to the submanifold over which the integration takes place, then integration by parts does not work. It turns out that one must both define a correct Littlewood-Paley projection and use a $TT^*$ argument in order to carry the argument out rigorously. See Section~\ref{sec_I_2} for details, including a comment at the end of the same section on our choice of Littlewood-Paley projection.

\section{Preliminary lemmas: Van der Corput estimates}\label{sec_prelim_VdC}

We begin by recording several estimates of van der Corput type as  stated in \cite{SWCarl}. Let 
\[Q_\lam(x)  = \sum_{0 \leq |\al| \leq d} \lam_\al x^\al\]
 be a polynomial of degree $d$ in $\R^m$, with coefficients $\lam_\al \in \R$. Note that for the moment we allow constant, linear, quadratic, and higher order terms in $Q_\lam$.  Let 
 \beq\label{isotropic_norm}
 \|\lam\| = \sum_{1 \leq |\al| \leq d} |\lam_\al| 
 \eeq
  denote the isotropic norm of the coefficients of non-constant terms.
  Let 
   \beq\label{isotropic_norm_0}
 \llbracket \lam\rrbracket = \sum_{0 \leq |\al| \leq d} |\lam_\al|  = \| \lam \| + |Q_\lam(0)|
 \eeq
 denote the isotropic norm, including the constant term.
   We now quote Proposition 2.1 of \cite{SWCarl}:
\begin{lemma}\label{lemma_Prop2.1}
For any $C^1$ function $\psi$ defined on the unit ball $B_1(\R^m)$, with $\| \psi \|_{C^1} \leq 1$, and for any convex subset $\Om \subseteq B_1(\R^m)$, we have
\beq
\left| \int_\Om e^{iQ_\lam(x)} \psi(x)dx \right| \leq C  \| \lam \|^{-1/d} ,
\eeq
where the constant $C$ depends on the dimension $m$ and the degree $d$ of $Q_\lam$, but not otherwise on $Q_\lam$, $\psi$ or $\Omega$.
\end{lemma}
Note that this upper bound is of course independent of the constant term in $Q_\lam(x)$. The result of the lemma continues to hold with the unit ball $B_1(\R^m)$ replaced by a Euclidean ball of any other fixed radius, such as $B_2(\R^m)$.

\xtra{
We do not repeat the proof here, as it is carefully outlined in \cite{SWCarl}, but we do make one observation. At certain points we will apply this lemma when it is the coefficient of the linear terms in $P(x)$ that is known to be large. At first glance, it might seem this would require the additional assumption that the (first partial) derivative of the phase is monotone (as in the usual one variable van der Corput estimate; see Proposition 2 of Chapter 8.1.2 in \cite{SteinHA}). However, this condition is not required, because if it is the first derivative of $P_\lam$ that is large, say $|P'_\lam(x)| \geq c\|\lam\|$, one can simply estimate (say in the one-variable situation; the higher dimensional situation then follows by integrating over the extra variables): $$\int_{-1}^1 e^{iP_\lam(x)} \eta(x) dx = \int_{-1}^1 \frac{de^{iP_\lam(x)}}{dx}  \frac{1}{P'_\lam(x)} \eta(x) dx = \int_{-1}^1 e^{iP_\lam(x)} \frac{1}{iP'_\lam(x)} \left(\frac{P''_\lam(x)}{iP_\lam'(x)} \eta(x) - \eta'(x) \right) dx,$$ which is bounded in absolute value by $c\|\lambda\|^{-1}$ because $|P'_\lam(x)| \geq c \|\lambda\|$ and $|P''_\lam(x)| \leq c\|\lambda \|$.}

We will also require a lemma estimating the measure of the set where the real-valued polynomial $Q_\lam$ takes small values.

\begin{lemma}\label{lemma_Prop2.2}
Let $Q_\lam$ be a polynomial as above (possibly including a nonzero constant term). For every $\rho >0$,
\beq\label{poly_small_set_bound}
| \{x \in B_1( \R^m): |Q_\lam(x)| \leq \rho \} |\leq C \rho^{1/d}\|\lam\|^{-1/d},
\eeq
where the constant $C$ depends only on the dimension $m$ and the degree $d$ of $Q_\lam$.
\end{lemma}

This lemma is a slight variant of Proposition 2.2 of \cite{SWCarl}; in \cite{SWCarl} they assumed that $Q_\lam(y)$ has no constant term, but their proof carries over to the present case.

We will also require an improvement of this lemma when the polynomial $Q_\lam$ has a constant term that is large. (The estimate of the above lemma does not see the constant term of $Q_\lam$, by definition of $\|\lam\|$.) We prove:
\begin{lemma}\label{lemma_Prop2.2_0}
Let $Q_\lam$ be a polynomial as above (possibly including a nonzero constant term). For every $\rho >0$,
\beq\label{poly_small_set_bound_0}
| \{x \in B_1( \R^m): |Q_\lam(x)| \leq \rho \} |\leq C' \rho^{1/d}\llbracket \lam \rrbracket^{-1/d},
\eeq
where the constant $C'$ depends only on the dimension $m$ and the degree $d$ of $Q_\lam$.
\end{lemma}
We note that both Lemma \ref{lemma_Prop2.2} and \ref{lemma_Prop2.2_0} continue to hold if $B_1(\R^m)$ is replaced by any other Euclidean ball of fixed radius, such as $B_2(\R^m)$.

To prove Lemma \ref{lemma_Prop2.2_0},
first suppose that $\| \lam \| \geq \llbracket \lam \rrbracket /4$. Then  Lemma \ref{lemma_Prop2.2} shows
\[  |\{ x \in B_1 : |Q_\lam(x)| \leq \rho \}|  \leq C \rho^{1/d}\|\lam\|^{-1/d}  \leq C 4^{1/d} \rho^{1/d} \llbracket \lam \rrbracket^{-1/d}.\]
Next, we suppose that $\| \lam \| < \llbracket \lam \rrbracket /4$ so that necessarily 
\beq\label{Qbig}
 |Q_\lam(0) | = \llbracket \lam \rrbracket - \| \lam \| \geq \frac{3}{4}\llbracket \lam \rrbracket .
 \eeq
In this case we crucially use the fact that the set we consider lies inside a ball of fixed radius;  for $x \in B_1$, (\ref{Qbig}) implies that
\beq\label{Qball}
 |Q_\lam(x)|  \geq |Q_\lam(0)| - |\sum_{1 \leq |\al| \leq d} \lam_\al x^\al| 
	\geq  |Q_\lam(0)| - \sum_{1 \leq |\al| \leq d} |\lam_\al|  = |Q_\lam(0)| - \| \lam \| \geq \frac{1}{2} \llbracket \lam \rrbracket. 
	\eeq
In particular, if $\llbracket \lam \rrbracket > 2\rho$, this implies the set $\{ x \in B_1 : |Q_\lam(x) | \leq \rho \}$ is empty and hence (\ref{poly_small_set_bound_0}) trivially holds.
On the other hand if $\llbracket \lam \rrbracket  \leq 2\rho$ then we can prove (\ref{poly_small_set_bound_0}) directly. We need only note that in this case, 
\[ \left( \frac{\rho}{\llbracket \lam \rrbracket} \right)^{1/d} \geq \left( \frac{1}{2} \right)^{1/d}.\]
Thus  we may apply the trivial bound:
\[  | \{ x \in B_1 : |Q_\lam(x) | < \rho \} | \leq |B_1| \leq |B_1| 2^{1/d} \left( \frac{\rho}{\llbracket \lam \rrbracket} \right)^{1/d} .\]
This establishes (\ref{poly_small_set_bound_0}) with the constant $C' = \max \{ C 4^{1/d}, |B_1|2^{1/d} \}$. 

\xtra{Finally, we return to the proof of Lemma \ref{lemma_Prop2.2}.
We recall from Lemma 2.2 of \cite{SWCarl} that there exists an index $k$ with $1 \leq k \leq d$ and a unit vector $\xi \in \R^m$ such that for all $x \in B_1$,
\beq\label{xi_deriv}
 |(\xi \cdot \partial_x)^k Q_\lam(x)| \geq c \| \lam \|. 
 \eeq
We now pick an orthogonal coordinate system $x_1,\ldots, x_n$ so that $x_1$ lies in the direction of the unit vector $\xi$. Then $(\xi \cdot \partial_x)^k = \frac{\partial^k}{\partial x_1^k}$.
We call upon the following single-variable lemma  (see for example Christ \cite{Chr85}):
\begin{lemma}\label{lemma_Christ}
For any $k \in \N$ there exists $C< \infty$ such that for any interval $I \subset \R$, any $f \in C^{(k)}(I)$, and any $\rho>0$, 
\[ |\{ x \in I: |f(x)| \leq \rho \} | \leq C \rho^{1/k} ( \inf_{x \in I} |D^{(k)} f(x)| )^{-1/k},\]
where $D^{(k)}$ denotes the $k$-th derivative.
\end{lemma}
We apply this with $I$ being the unit interval with respect to $x_1$, and with the lower bound (\ref{xi_deriv}) (and use the fact that the remaining variables $x_2,\ldots, x_m$ also lie in unit intervals) to conclude that Lemma \ref{lemma_Prop2.2} holds.}

\xtra{We note that Lemma \ref{lemma_Prop2.2} and \ref{lemma_Prop2.2_0} continue to hold if $B_1$ is replaced by any other Euclidean ball of fixed radius, such as $B_2$. Indeed, to prove the analogue of Lemma \ref{lemma_Prop2.2_0} for $\{ x \in B_{R} \subset \R^m: |Q_\lam(x)| \leq \rho \}$ for a fixed $R \geq 1$, we would consider the cases $\| \lam \| \geq  (4R^d)^{-1}  \llbracket \lam\rrbracket$ and $\| \lam \| \leq(4R^d)^{-1}  \llbracket \lam \rrbracket$; in the latter case we would replace (\ref{Qbig}) by 
\[|Q_\lam(0)| =  \llbracket \lam \rrbracket  - \| \lam \| \geq (1 - (4R^d)^{-1}) \llbracket \lam \rrbracket \geq \frac{3}{4} \llbracket \lam \rrbracket. \]
We would then replace (\ref{Qball}) by 
\begin{eqnarray*}
 |Q_\lam(x)|  &\geq& |Q_\lam(0)| - |\sum_{1 \leq |\al| \leq d} \lam_\al x^\al| \\
	&\geq & |Q_\lam(0)| - R^d \sum_{1 \leq |\al| \leq d} |\lam_\al| \\
	&  =& |Q_\lam(0)| -R^d \| \lam \| \\
	& \geq & \frac{3}{4} \llbracket \lam \rrbracket - \frac{1}{4}\llbracket \lam \rrbracket   
	\geq \frac{1}{2} \llbracket \lam \rrbracket. 
	\end{eqnarray*}
	and the argument would proceed as before, proving the analogue of (\ref{poly_small_set_bound_0}) with the constant $C' = \max \{ C 4^{1/d}, |B_R|2^{1/d} \}$. 
}

\subsection{Novel Van der Corput estimates for kernels}
 We now state the novel van der Corput estimates we will apply to bound the kernels of various operators of the form $TT^*$; these kernel estimates are in fact the heart of the paper. 
\begin{prop}[$K_\flat^{\nu,\mu}$ van der Corput, $n \geq 2$]\label{prop_Kflat}
Fix any dimension $n \geq 1$, and any degree $d \geq 2$. 
Suppose $\nu = (\nu_{\alpha})_{2 \leq |\alpha| \leq d}$ denotes the coefficients of a real-valued polynomial
$$
Q_{\nu}(y) = \sum_{2 \leq |\alpha| \leq d} \nu_{\alpha} y^{\alpha}
$$
on $\R^n$ that has no linear terms, and similarly for $\mu = (\mu_{\alpha})_{2 \leq |\alpha| \leq d}$. Let 
$$
\|\nu\| = \sum_{2 \leq |\alpha| \leq d} |\nu_{\alpha}|, \qquad \|\mu\| = \sum_{2 \leq |\alpha| \leq d} |\mu_{\alpha}|.
$$
Given a $C^1$ function $\Psi(u,z)$ supported on $B_2(\R^n) \times B_1(\R^n)$, define
\[ K^{\nu,\mu}_{\flat}(u)
= \int_{\R^n} e^{iQ_{\nu}(u+z)-iQ_{\mu}(z)} \Psi(u,z) dz.\]
Suppose furthermore that 
\[ \| \Psi \|_{C^1} := \sup_{(u, z) \in B_2(\R^n) \times B_1(\R^n)} \left( |\Psi(u,z)| + |\nabla_z \Psi(u,z)|\right) \leq 1.
\]
Then there exists a small constant $\delta > 0$ (depending only on $d$), such that the following holds: 
if $\nu$, $\mu$ satisfy $$r \leq \|\nu\|,\|\mu\| \leq 2r$$ for some $r\geq1$,
then there exists a small measurable set $G^{\nu} \subset B_2(\R^n)$ depending on $\nu$ (but not on $\mu$ nor $\Psi$), with $$|G^{\nu}| \leq C r^{-\delta},$$ such that
$$|K^{\nu,\mu}_{\flat}(u)| \leq C  \left( r^{-\delta} \chi_{B_2}(u) + \chi_{G^{\nu}}(u) \right).$$ 
\end{prop}
This is effectively a result of Stein and Wainger (namely Lemma 4.1 in \cite{SWCarl}, in the case $h=1$, in their notation). It follows from a clever application of Lemmas \ref{lemma_Prop2.1} and \ref{lemma_Prop2.2}, and crucially uses the assumption that $Q_\lam (y)$ lacks linear terms.

In the Radon setting,
we require a more elaborate version; we state here a result in dimension $n=2$ that we prove in Section \ref{sec_vdC_n2}.  We reserve the more technical statement of the result for general dimensions $n \geq 2$ for Proposition \ref{prop_Ksharp_n} in Section \ref{sec_high_dim}.

\begin{prop}[$K_\sharp^{\nu,\mu}$ van der Corput, $n = 2$]\label{prop_Ksharp}
Fix the dimension $n = 2$ and a degree $d \geq 2$. Let 
\[\Pscr = \{ p_2(y), p_3(y),\ldots, p_d(y)\}\]
be a set of real-valued polynomials on $\mathbb{R}^2$, where each $p_j(y)$ is homogeneous of degree $j$, and $p_2(y) \ne C|y|^2$ for any nonzero constant $C$. Let 
\[\Lambda = \Lambda (\Pscr) = \{ 2 \leq m \leq d \colon p_m(y) \not\con 0\}.\]
 For $\nu = (\nu_2,\dots,\nu_d) \in \R^{d-1}$, let 
 \[ P_{\nu}(y) = \sum_{m=2}^d \nu_m p_m(y)\]
and
\[ \| \nu \| = \sum_{m \in \Lambda} |\nu_m|,\]
and define $P_\mu(y)$ and $\| \mu \|$ similarly for $\mu = (\mu_2, \dots, \mu_d) \in \R^{d-1}$.

 Given a $C^1$ function $\Psi(u,z)$ supported on $B_2(\R^2) \times B_1(\R^2)$, 
define
\[ K^{\nu,\mu}_{\sharp}(u, \tau)
= \int_\R e^{iP_{\nu}(u+z)-iP_{\mu}(z)} \Psi(u,z)d\sigma, \]
where the $z$ in the integral is defined implicitly in terms of $u, \tau,\sigma$ by 
\beq\label{z_dfns}
\tau = \frac{u_1z_1+u_2z_2}{|u|}, \qquad \sig = \frac{-u_1z_2 + u_2z_1}{|u|}.
\eeq
Suppose furthermore that 
\[\| \Psi \|_{C^1(\R)} := \sup_{(u,z) \in B_2(\R^2) \times B_1(\R^2)} \left( |\Psi(u,z)| + |\frac{\partial}{\partial \sig} \Psi(u,z)|\right) \leq 1.\]
Then there exists a small constant $\delta > 0$ (depending only on $d$) such that the following holds:
if $\mu$, $\nu$ satisfy
$$
r \leq \|\nu\|,\|\mu\| \leq 2r
$$
for some $r \geq 1$, then there exists a small set $G^\nu \subset B_2(\R^2)$, and for each $u \in B_2(\R^2)$ a small set $F_u^\nu \subset B_1(\R)$, such that
$$
|G^\nu| \leq C r^{-\del}, \qquad |F^\nu_u| \leq C r^{-\del} \quad \text{for all $u \in B_2(\R^2)$},
$$
and
\beq\label{K_sharp_bd}
|K^{\nu,\mu}_{\sharp}(u, \tau)| \leq C
	\left( r^{-\delta} \chi_{B_2}(u) \chi_{B_1}(\tau) + \chi_{G^{\nu}}(u) \chi_{B_1}(\tau) + \chi_{B_2}(u) \chi_{F_u^{\nu}}(\tau) \right).
\eeq
The choices of  the small sets $G^{\nu}$ and $F^{\nu}_u$ are independent of both $\mu$ and $\Psi$.
\end{prop}
The significance of this result relative to that of Proposition \ref{prop_Kflat} for $K_\flat^{\nu,\mu}$ when $n = 2$ is that the integral $K_\sharp^{\nu,\mu}$ is now only a one-dimensional integral, yet we are still able to capture decay in $r$, as long as the polynomial phase occurring in $K_\sharp^{\nu,\mu}$ is of a more restrictive form than in $K_\flat^{\nu,\mu}$.

\section{Preliminary lemmas: Littlewood-Paley decomposition}\label{sec_LP}

We now give an explicit construction of a Littlewood-Paley decomposition we will employ to bound the square function  (\ref{Sf0}) on $L^2$.
We start with a function $\chi \in C_c^\infty([-1,1])$ such that $\chi$ is identically 1 on a small neighborhood of the origin. We then set
\[ \Phi(\tau ) = \chi(\tau/4) - \chi(\tau),\]
so that $\Phi$ is a smooth function with support in $[-4,4]$ that vanishes in a small neighborhood of the origin. We observe that if we set $\Phi_j(\tau) = \Phi(2^{2j} \tau)$, then for any $\tau \neq 0$,
\[ \sum_{j \in \Z} \Phi_j(\tau) = 1,\]
by a telescoping argument.
We now define a Schwartz function $\Delta$ on $\R$ by 
\[\Delta (t)= \check{\Phi}(t).\]
 We also define corresponding scaled versions by $\Delta_j (t)=\check{\Phi}_j(t)$, so that 
 \[\Delta_j(t) = 2^{-2j}\Delta(2^{-2j}t).\]
We observe that 
\beq\label{Del_cancel}
 \int_\R \Delta (t) dt = \int_\R \Delta_j(t) dt = \Phi(0) = 0, \quad \text{for any $j$}.
 \eeq
Given any function $f \in \Scal(\R^{n+1})$ we define 
\[\Delta_j f (x,t) = \int f(x,t-s) \Del_j(s) ds,\]
so that $\Delta_j f$ is also of Schwartz class. 
Note that we use $\Del_j$ to denote both the kernel and the operator itself.
Also each $\Del_j$ extends to a bounded linear operator on $L^2(\mathbb{R}^{n+1})$. Finally, we have defined the Littlewood-Paley decomposition so that it is compatible with parabolic scalings $(x,t) \mapsto (\del x, \del^2 t)$.

The following is a standard result, which follows from an easy application of the Fourier transform, and Fubini's theorem:
\begin{prop}\label{prop_Del_N}
For all $f \in L^2(\R^{n+1})$, the partial sums $\sum_{j=-N}^N \Delta_j f$ converge to $f$ in $L^2(\R^{n+1})$ norm. 
Moreover, 
\beq\label{LNinf_norm}
 \left\| \left( \sum_{j=-\infty}^{\infty} |\Delta_j f|^2 \right)^{1/2} \right\|_{L^2(\R^{n+1})} \leq C \|f\|_{L^2(\R^{n+1})}.
\eeq
\end{prop}

\xtra{One can see the first result of the lemma using the Fourier transform: 
$$\left\|\sum_{j=-N}^N \Delta_j f - f\right\|_{L^2}^2 = \iint_{\R^{n+1}} |\widehat{f}(\xi,\tau)|^2 |\chi(2^{-2N-2} \tau) - \chi(2^{2N} \tau) - 1|^2 d\xi d\tau \to 0$$ as $N \to \infty$ by the dominated convergence theorem, since for $\tau \ne 0$, $\chi(2^{-2N-2} \tau) \to 1$ and $\chi(2^{2N} \tau) \to 0$ as $N \maps \infty$.\xtra{ the function $\chi(2^{-2N-2}\tau)$ will approach the (unit) value of $\chi$ at the origin, and $\chi(2^{2N} \tau)$ will escape the support of $\chi$. }

Next, the bound (\ref{LNinf_norm}) will follow if we can show 
\beq\label{phi_hat_sum}
\sum_{j=-\infty}^{\infty} |\Phi(2^{2j} \tau)|^2 \leq C \quad \text{uniformly in $\tau$};
\eeq
by scale invariance, one only needs to prove this last estimate when $1 \leq |\tau| \leq 4$, and we use
 \beq\label{phi_hat}
 |\Phi(\tau)| \leq C \min\{ |\tau|, |\tau|^{-1} \}.
 \eeq
Then for $1 \leq |\tau| \leq 4$, it follows from (\ref{phi_hat}) that  $|\Phi(2^{2j} \tau)| \leq C 2^{-2j}$ if $j > 0$, and $\leq C 2^{2j}$ if $j \leq 0$, which gives the desired bound (\ref{phi_hat_sum}). }

We will require an additional Littlewood-Paley decomposition in order to gain a reproducing property. Fix a smooth function $\tilde{\Phi}$ with compact support in $[-8,8]$ such that $\tilde{\Phi} \con 1$ on the support of $\Phi$ and $\tilde{\Phi} \con 0$ in a small neighborhood of the origin. Correspondingly, we define the scaled version $\tilde{\Phi}_j =  \tilde{\Phi}(2^{2j} \tau)$ and the kernels $\tilde{\Delta}_j$ and associated operators $\tilde{\Del}_j$ as above. In particular $\Phi\tilde{\Phi} = \Phi$ and thus as operators 
\beq\label{dup}
\Del_j \tilde{\Del}_j = \Del_j.
\eeq
\xtra{  Moreover, we remark that 
starting with any $f \in \Scal$ the resulting function $\Del_j \tilde{\Del}_j f$ remains of Schwartz class. }

We summarize the additional simple properties we require as follows:
\begin{prop}\label{tilde_del_prop}
For all $f \in L^2$, 
\beq\label{fDelP}
f =  \sum_{j=-\infty}^{\infty} \Delta_j \tilde{\Delta}_j f
\eeq
where the convergence of the sum on the right hand side is taken in the $L^2$ sense. 
In addition,
\beq\label{forwardLP}
 \left\| \left( \sum_{j=-\infty}^{\infty} |\tilde{\Delta}_j f|^2 \right)^{1/2} \right\|_{L^2(\R^{n+1})} \leq C \|f\|_{L^2(\R^{n+1})}
\eeq
and
\beq\label{LN_norm}
 \left\| \left( \sum_{j=-\infty}^{\infty} |\Delta_j \tilde{\Delta}_j f|^2 \right)^{1/2} \right\|_{L^2(\R^{n+1})} \leq C \|f\|_{L^2(\R^{n+1})}.
\eeq
 \end{prop}
The first and third properties are simple consequences of (\ref{dup}) and Proposition \ref{prop_Del_N}; the second property follows from a similar argument to that needed for (\ref{LNinf_norm}).

\subsection{Convolutions}\label{sec_uDel}
Later on we will require a third type of Littlewood-Paley operator, defined by $\Del_j * \Del^-_j$, where $\Del_j^-(x) := \Del_j(-x)$. Precisely, we define
\[ \uDel_j(t) = \int \Del_j(w+t) \Del_j(w) dw,\]
so that $\uDel_j$ is also a Schwartz function.
We note the following properties of $\uDel_j$:
\begin{lemma}\label{lemma_uDel}
 For all $j \in \Z$,
\beq\label{uDel_rescale}
 \uDel_j(t) = \frac{1}{a^2}\uDel_{j-k} (\frac{t}{a^2}) \qquad \text{if $a=2^k$}, 
 \eeq
and
\beq\label{Del_cancellation} \int_\R \uDel_j(t)dt    =   0. 
\eeq
In addition, the $\uDel_j$ are in $L^1(\R)$ with uniform norm, independent of $j$. Finally, for $\psi(t) = 1/(1+t^2)$ we set 
\beq\label{psi_normal}
\psi_j(t) = 2^{-2j}\psi(2^{-2j}t)
\eeq
 to be the $L^1$-normalization compatible with parabolic dilations. Then the derivative $(\uDel_j)'$ satisfies
\beq\label{uDel_prime}
 |(\uDel_j)'(t)| \leq C 2^{-2j} \psi_j(t).
\eeq
\end{lemma}
The first three properties of $\uDel_j$ are all simple consequences of the definition of $\Del_j$.
\xtra{In order to prove (\ref{uDel_prime}), 
we define the notation $(\Del')_j (t) = 2^{-2j} \Del' (2^{-2j}t)$ and note that 
\begin{eqnarray*}
(\uDel_j)'(t) &=& \frac{d}{dt} \int \Del_j(w+t) \Del_j(w)dw \\
&=& \int \frac{d}{dt} [2^{-2j} \Del (2^{-2j}(w+t)) ]\Del_j(w)dw \\
&=& 2^{-2j} \int 2^{-2j} \Del' (2^{-2j}(w+t)) \Del_j(w)dw \\
& = & 2^{-2j} \int (\Del')_{j} (w+t) \Del_j(w)dw \\
& = & 2^{-2j}\cdot 2^{-2j} \int_\R \Del' (w+2^{-2j}t) \Del(w)dw.
\end{eqnarray*}
We now consider only the integral, and apply the rapid decay of $\Del$ and $\Del'$ to deduce that for any $N \geq 1$,
\[  \int_\R \Del' (w+2^{-2j}t) \Del(w)dw \leq C \int_\R \frac{1}{1 + (w+2^{-2j}|t|)^{2N}} \frac{1}{1+w^{2N}} dw.\]
We consider the integral in three portions, first integrating over $|w| \geq 2 \cdot 2^{-2j}|t|$, from which the contribution is
\[ \lesssim \frac{1}{1 + (2^{-2j}|t|)^{2N}} \int_\R \frac{1}{1+w^{2N}} dw \lesssim \frac{1}{1 + (2^{-2j}|t|)^{2N}} .\]
We may bound the portion $|w| \leq (1/2) 2^{-2j}|t|$ similarly. Finally, in the range $(1/2)2^{-2j}|t| \leq |w| \leq 2 \cdot 2^{-2j}|t|$ we bound the integral by 
\[ \int_{|w| \approx 2^{-2j}|t|} \frac{1}{1+w^{2N}}dw \lesssim \frac{1}{1 + (2^{-2j}|t|)^{2N}} \cdot (2^{-2j}|t|),\]
which is sufficient, since $N \geq 1$ is arbitrary. This implies the estimate  (\ref{uDel_prime}) is valid.}
Also, note that by dilation invariance, it suffices to establish (\ref{uDel_prime}) when $j =0$. In that case, the desired estimate reduces to $|\uDel_0'(t)| \leq C / (1+t^2)$, which follows simply from the fact that $\uDel_0$ is Schwartz.

In addition, we will use the following consequence of the mean value theorem:
\begin{lemma}\label{lemma_uDel_diff}
Let $\psi(t) = 1/(1+t^2)$ and set $\psi_j(t) = 2^{-2j}\psi(2^{-2j}t)$.
For each $j  \geq 0$, for any $|\xi| \leq 2$, we have 
\[ |\uDel_j (t + \xi) - \uDel_j (t)|  \leq C 2^{-2j} |\xi| \psi_j(t).\]
In addition, we note that $\psi_j$ is a non-negative integrable function on $\R$, with $L^1$ norm uniformly bounded, independent of $j.$ 
\end{lemma}

Certainly, the mean value theorem applies to the Schwartz function $\uDel_j$, showing that for any fixed $t$ and any $|\xi| \leq 2$ we have 
\[ |\uDel_j (t + \xi) - \uDel_j (t)|  \leq |\xi| |(\uDel_j)'(t+\xi_0)|,\]
for some $|\xi_0| \leq |\xi| \leq 2$. We now apply the estimate (\ref{uDel_prime}) to conclude that 
\[ |\uDel_j (t + \xi) - \uDel_j (t)|  \leq C  2^{-2j} |\xi|   \psi_j(t+\xi_0).\]
\xtra{Now we need only verify  for $j \geq 0$  the comparison that}
It then remains to observe that
\[ \psi_j(t+\xi_0) \leq C \psi_j(t)
\]
for all $t \in \R$, $|\xi_0 | \leq 2$, as may be easily verified. \xtra{We need only confirm that
\beq\label{psi_comp}
 2^{-2j} \frac{1}{1 + 2^{-4j}(t+\xi_0)^2} \leq C 2^{-2j} \frac{1}{1 + 2^{-4j}t^2} .
 \eeq
Certainly, in the case $|\xi_0| \leq |t|/2$ we have $|t + \xi_0| \geq |t|/2$, so that (\ref{psi_comp}) is valid. On the other hand if $ |\xi_0| \geq |t|/2$ then necessarily $|t| \leq 4$ so that $|t + \xi_0| \leq 6$ and hence (\ref{psi_comp}) is implied by the  inequality 
\[  \frac{1}{ 1+ 2^{-4j} (t + \xi_0)^2} \leq C,\]
valid for all $j \geq 0$.}
 
 \subsection{Antiderivatives}
 
 Finally, we record for later use the following fact about antiderivatives of $\Del_j$ and $\uDel_j$:
\begin{lemma}\label{lemma_Del_anti}
There exist Schwartz functions $\widetilde{\Del}$ and $\widetilde{\uDel}$ on $\mathbb{R}$ such that upon setting
$\widetilde{\Del}_j(t)= 2^{-2j}\widetilde{\Del} (2^{-2j}t)$ and $\widetilde{\uDel}_j(t)= 2^{-2j}\widetilde{\uDel} (2^{-2j}t)$
we have for every $j \in \Z$
  \beq\label{DelDel}
  \Del_j(t) = 2^{2j} \left(\frac{d}{dt} \widetilde{\Del}_j \right)(t), \qquad  \uDel_j(t) = 2^{2j}  \left(\frac{d}{dt} \widetilde{\uDel}_j\right)(t).
  \eeq
In particular, $\widetilde{\Del}_j$ and $\widetilde{\uDel}_j$ are uniformly in $L^1(\R)$.
  \end{lemma}
Again by dilation invariance, one could reduce to the case $j = 0$. The lemma would then follow from the following claim: if $F$ is a Schwartz function with $\int_{\mathbb{R}} F(t) dt= 0$, then there exists another Schwartz function $\tilde{F}$ such that $\tilde{F}'(t) = F(t)$. To prove this claim, let  
  \[\widetilde{F}(t) = \int_{-\infty}^tF(\tau) d\tau.\]
  Then it is immediate that $F(t)  = \frac{d}{dt}\widetilde{F}(t)$, and we need only verify that $\widetilde{F}(t)$ is a Schwartz function. To see that $\widetilde{F}$ exhibits rapid decay as $t \maps -\infty,$ we use the fact that $F$ has rapid decay, so that 
  \[ |\widetilde{F}(t)| \leq C \int_{-\infty}^t \frac{1}{1+|\tau|^N} d\tau \leq C \frac{1}{1+|t|^{N-1}},\]
  as $t \maps - \infty$, for arbitrary large $N$. By the assumed fact that $F$ integrates to zero, we can also use the alternative representation 
  \[ \widetilde{F}(t) = -\int_t^\infty F(\tau)d\tau,\]
  which similarly allows us to conclude that $\widetilde{F}$ exhibits rapid decay as $t \maps +\infty$.
 Finally, we note that $\widetilde{F}(t)$ is infinitely differentiable and its derivatives exhibit rapid decay, since $(\widetilde{F})' = F$ and $F$ is Schwartz.
  Thus $\widetilde{F}$ is Schwartz.

\xtra{We will construct an appropriate $\widetilde{F}$ in general for any Schwartz function $F(t)$ such that $\int_\R F(t)dt=0$, with dilations $F_j(t) = 2^{-2j}F(2^{-2j}t)$. The result then applies to $\Del$ and $\uDel$, since they satisfy the key cancellation identities (\ref{Del_cancel}) and (\ref{Del_cancellation}), respectively.
  
We begin by defining
  \[\widetilde{F}(t) = \int_{-\infty}^tF(\tau) d\tau.\]
  Then it is immediate that $F(t)  = \frac{d}{dt}\widetilde{F}(t)$, and we need only verify that $\widetilde{F}(t)$ is a Schwartz function. To see that $\widetilde{F}$ exhibits rapid decay as $t \maps -\infty,$ we use the fact that $F$ has rapid decay, so that 
  \[ |\widetilde{F}(t)| \leq C \int_{-\infty}^t \frac{1}{1+|\tau|^N} d\tau \leq C \frac{1}{1+|t|^{N-1}},\]
  as $t \maps - \infty$, for arbitrary large $N$. By the assumed fact that $F$ integrates to zero, we can also use the alternative representation 
  \[ \widetilde{F}(t) = -\int_t^\infty F(\tau)d\tau,\]
  which similarly allows us to conclude that $\widetilde{F}$ exhibits rapid decay as $t \maps +\infty$.
 Finally, we note that $\widetilde{F}(t)$ is infinitely differentiable and its derivatives exhibit rapid decay, since $(\widetilde{F})' = F$ and $F$ is Schwartz.
  Thus $\widetilde{F}$ is Schwartz.

In order to accommodate scalings, we observe that
 \[\int_{-\infty}^t F_j(\tau)d\tau = \int_{-\infty}^{2^{-2j}t} F(\tau) d\tau = \widetilde{F}(2^{-2j}t) = 2^{2j} \cdot 2^{-2j}  \widetilde{F}(2^{-2j}t) .\]
  Upon setting $\widetilde{F}_j (t) =2^{-2j}  \widetilde{F}(2^{-2j}t)$, we may take derivatives on both sides of the above equation to see that
  $F_j(t)= 2^{2j} \frac{d}{dt} \widetilde{F}_j(t),$ as desired.}
  
  \xtra{We can't actually apply the following argument, because we do not yet know that $\widetilde{F}$ is in $L^1$ so we cannot take the F.t. But this is one style of argument that could work in some similar cases.
  By (\ref{DelDel}) we see that on the Fourier side,
  $\xi (\widetilde{F})\hat{\;}(\xi) = (F)\hat{\;}(\xi)$, so that  
  \[(\widetilde{F})\hat{\;}(\xi)= \frac{(F)\hat{\;}(\xi)}{\xi}.\]
  The function on the right hand side is a Schwartz function. Indeed, away from the origin it clearly is smooth with rapid decay; near the origin, it remains smooth since $(F)\hat{\;}(\xi)$ vanishes at the origin, by (\ref{Del_cancel})).}

\section{The proof of Theorem \ref{thm_Carleson}}\label{sec_reduction}

In this section, we show how to deduce our main result, Theorem~\ref{thm_Carleson}, from Theorem~\ref{Osc_Thm}. Suppose we are given $d-1$ real polynomials $p_2(y),\dots,p_d(y)$ on $\mathbb{R}^n$, where each $p_j$ is homogeneous of degree $j$, and $p_2(y) \ne C|y|^2$ for any nonzero constant $C$. We set $\Lambda = \{2 \leq m \leq d: p_m(y) \not\con 0 \}$. For $\lam = (\lam_2,\dots,\lam_d) \in \mathbb{R}^{d-1}$, we define $P_{\lam}(y)$ as in (\ref{allowable_poly}) by
\[P_{\lam}(y) = \sum_{m=2}^d \lam_m p_m(y),\]
and $T_{\lam} f$ as in  (\ref{T_lam_par}) by
\[T_{\lam}f(x,t) = \int_{\R^n} f(x-y,t-|y|^2) e^{iP_{\lam}(y)} K(y) dy.\]
We want to show that for all $1<p<\infty$ and all Schwartz functions $f$ on $\R^{n+1}$,
$$
\|\sup_{\lam \in \mathbb{R}^{d-1}} |T_{\lam} f|\, \|_{L^p(\mathbb{R}^{n+1})} 
\leq A_p \|f\|_{L^p(\mathbb{R}^{n+1})}.
$$
To do so, note that in the integral defining $T_{\lam}f(x,t)$, if $y$ is sufficiently small (with respect to $\lam$) then $P_{\lam}(y)$ is approximately zero and $e^{iP_{\lam}(y)}$ can be approximated by 1. To make this precise, recall we have already defined a homogeneous norm $\|\lam\|$ on $\R^{d-1}$ in (\ref{lam_Rn_def}). 
However, what is more relevant here is a non-isotropic norm on $\R^{d-1}$, which we define by
$$
N(\lam) = \sum_{m \in \Lambda} |\lam_m|^{1/m}
$$
for $\lam \in \R^{d-1}$. The key is that $|P_{\lam}(y)| \leq 1$ whenever $|y| \leq c N(\lam)$, where $c$ is some fixed constant dependent only on the dimension and the degree $d$. \xtra{This may be seen as follows: First suppose there is at least one $\lam_\al$ with $|\lam_\al| \geq 1$, so that $\| \lam \| \geq 1$. Then we write
\[ N(\lam) = \sum_{|\lam_\al| < 1} |\lam_\al|^{1/ |\al|} + \sum_{|\lam_\al| \geq 1} |\lam_\al|^{1/|\al|} \leq C_d + \sum_{|\lam_\al| \geq 1} |\lam_\al| \leq C_d + \| \lam \| \leq (C_d +1) \|\lam \|,\]
where $C_d$ is the number of multi-indices of degree up to $d$.
In the other case, when $|\lam_\al| <1$ for all $\al$, the assumption that $N(\lam) \geq 1$ still shows that there exists an $\al$ for which $|\lam_\al|^{1/|\al|} \geq C_d^{-1}$, and hence $\| \lam \| \geq |\lam_\al| \geq C_d^{-|\al|}$, or in other words $C_d^{|\al|} \|\lam \| \geq 1$, and so we may replace $\|\lam\|$ in the first argument by $C_d^{|\al|} \|\lam \|$.} We will make use of this observation very soon.

The Calder\'{o}n-Zygmund kernel $K$ admits a decomposition (see Chapter XIII of \cite{SteinHA}) as 
\[ K(x) = \sum_{j=-\infty}^\infty 2^{-nj}\phi_j(2^{-j}x) \]
where each $\phi_j$ has the following properties:
\newline \indent (i) $\phi_j$ is a $C^1$ function with support in $1/4 < |x| \leq 1$ 
\newline \indent (ii) $ | \partial_x^\al \phi_j(x)| \leq C$ for $0 \leq |\al| \leq 1$ for some constant ${C}$ that is uniform in $j$,
\newline \indent  (iii) $\int_{\R^n} \phi_j(x)dx =0$ for every $j$.\\
This allows us to decompose $K$ as in \cite{SWCarl}: precisely, given $\lam \in \R^{d-1}$, we split $K$ as 
\[ K= K_\lam^+ + K_\lam^- = \sum_{2^j< 1/N(\lam)} K_j +  \sum_{2^j \geq 1/N(\lam)} K_j,\]
where $K_j(x) = 2^{-nj} \phi_j(2^{-j}x)$.
  Accordingly we split $T_\lam$ into $T_\lam= T_\lam^+ + T_\lam^-$, where (respectively)
\[ T_\lam^\pm f(x,t) = \int_{\R^n} f(x-y,t-|y|^2)e^{iP_\lam(y)}K_\lam^\pm(y) dy.\]
To prove Theorem \ref{thm_Carleson}  it is sufficient to bound the $L^p$ norms of $\sup_\lam |T_\lam^+ f|$ and $\sup_\lam |T_\lam^- f|$ individually. In the support of $T_\lam^-$, where $2^jN(\lam) \leq1$, the phase $P_\lam(y)$ will not cause significant oscillation, and we will aim to remove the oscillatory factor and bound the remaining operator by maximal truncated singular Radon transforms. In the support of $T_\lam^+$, we would expect  the phase $P_\lam(y)$ to contribute significant oscillation; this portion of the operator leads to the operators $^{(\eta_k)}I_{2^k}^\lam$ and $^{(\eta_k)}J_{2^k}^\lam$ defined in (\ref{I_op_dfn_0}) and (\ref{J0}), and the square function (\ref{Sf0}).

\subsection{Bounding $T_\lam^-$}\label{sec_T_minus}
We first outline the treatment of $T_\lam^-$. Note that $K_\lam^-(y)$ is visibly supported where $|y| \leq N(\lam)^{-1}$, and in fact agrees precisely with $T_\lam$ if $|y| \leq (4N(\lam))^{-1}$. We may replace the oscillatory factor $e^{iP_\lam(y)}$ by $1$ with moderate error, since
\begin{multline}
 |e^{iP_\lam(y)}-1| \leq c \sum_{2 \leq m \leq d} |\lam_m| |p_m(y)| \leq c' \sum_{m \in \Lambda} N(\lam)^{m} |y|^{m}  \\
 =  c' \sum_{2 \leq m \leq d} (N(\lam)|y|)^m \leq c'' N(\lam) |y|.
 \end{multline}
Here we have used the fact that for every  $2 \leq m \leq d$, $|\lam_m| \leq N(\lam)^{m}$,\xtra{For indeed, for any exponent $\al_0$, 
\[
 N(\lam)^{|\al_0|} = \sum_{\al} |\lam_\al|^{|\al_0|/|\al|} = |\lam_{\al_0}| + \text{positive} \geq |\lam_{\al_0}|
\].}
followed by the assumption that $|y|N(\lam) \leq 1$. It follows from this and the usual bound (\ref{K_prop}) for $K$ that
\begin{align}
  T_\lam^-f(x,t) = & \int_{|y| \leq 1/N(\lam)} K_\lam^-(y)f(x-y,t-|y|^2) dy  
 \notag \\ & + O\left( N(\lam) \int_{|y| \leq 1/N(\lam)} |y|^{-n+1} |f|(x-y,t-|y|^2)dy \right). \label{T_lam_minus}
 \end{align}
Taking the supremum over $\lam$, the first term on the right hand side is dominated by a truncated maximal singular Radon transform;  the second is dominated by a maximal function along the paraboloid.

Precisely, we define the operators and state the results we require, beginning with the truncated maximal singular Radon transform.  Define a singular Radon transform along the paraboloid, initially acting on functions of Schwartz class, by
\beq\label{Hcal_dfn}
\cH f  (x,t)= \int f(x-y,t-|y|^2) K(y) dy,
\eeq
where $K(y)$ is a Calder\'{o}n-Zygmund kernel on $\mathbb{R}^n$. Then $\Hcal$ is known to be a bounded operator on $L^p$ (see for example \S4.5 of Chapter 11 of \cite{SteinHA}, or the earlier case for curves given in \cite{SWCurv}). For each $\ep>0$, let
$$\cH_{\ep} f (x,t) = \int_{|y| > \ep} f(x-y,t-|y|^2) K(y) dy$$
be a truncation of $\Hcal$.  The result we require is as follows: 
 
 \begin{thm}\label{thm_trunc_Radon}
 For every $1<p< \infty$, we have an \emph{a priori} inequality for $f$ of Schwartz class:
\beq\label{Hp}
\|  \sup_{\ep > 0} |\cH_{\ep} f| \; \|_{L^p(\R^{n+1})} \leq A_p \|f\|_{L^p(\R^{n+1})}.
\eeq
 \end{thm}
 As usual, by a limiting argument, it follows that the operator $\sup_{\ep >0} |\cH_{\ep} f|$ extends to a bounded operator acting on functions in $L^p$, satisfying the same bound (\ref{Hp}). 
 Theorem 5.1 can be deduced from a result of Duoandikoetxea and Rubio de Francia \cite{DuoRub86}; for completeness we adapt their proof to our setting, and provide the details in Section \ref{sec_appendix}.

It follows from Theorem \ref{thm_trunc_Radon} that the supremum over $\lam$ of the first term on the right hand side of (\ref{T_lam_minus}) is bounded on $L^p$, since 
 \begin{multline*}
  \| \sup_\lam \int_{|y| \leq 1/N(\lam)} K_\lam^-(y)f(x-y,t-|y|^2) dy \|_{L^p(\R^{n+1})} 
 \\
 	\leq \|\Hcal f\|_{L^p(\R^{n+1})} + \| \sup_\lam |\Hcal_\lam f| \|_{L^p(\R^{n+1})}.
	\end{multline*}

We now turn to the second term on the right hand side of (\ref{T_lam_minus}), which we compare to $\Mcal_\Rad f(x,t)$,  where $\Mcal_\Rad$ is the maximal Radon transform along the paraboloid defined in (\ref{MRad_dfn}), which is well-known to be bounded on $L^p(\R^{n+1})$ for $1<p<\infty$.
To make this comparison precise, let $ \chi_{A_{2^l}}(y)$ denote the characteristic function of the annulus $2^{l-1}< |y| \leq 2^{l}$ and $\chi_{B_{2^l}}(y)$  the characteristic function for the ball of radius $2^l$. Then
\begin{eqnarray*}
&  &	\hspace{-2cm}  \sup_{\lam} N(\lam) \int_{|y| \leq 1/N(\lam)} |y|^{-n+1} |f|(x-y,t-|y|^2)dy\\
	\hspace{2cm}	 &=& \sup_\lam N(\lam) \sum_{2^l \leq N(\lam)^{-1}} \int |f|(x-y,t-|y|^2) |y|^{-n+1} \chi_{A_{2^l}}(y) dy \\
		  & \leq & \sup_\lam N(\lam) \sum_{2^l \leq N(\lam)^{-1}} 2^{-(l-1)(n-1)} \int |f|(x-y,t-|y|^2) \chi_{B^{2^l}}(y) dy \\
		 & \leq &\sup_\lam N(\lam) N(\lam)^{-1} \Mcal_\Rad f(x,t)\\
		 & \leq & \Mcal_\Rad f(x,t)  .
		 \end{eqnarray*}
In total, we may conclude that 
\[ \| \sup_{\lam} |T_\lam^- f|\; \|_{L^p(\R^{n+1})} \leq \|f\|_{L^p(\R^{n+1})}, \]
for all $1<p<\infty$,
as desired. 

\subsection{Bounding $T_\lam^+$}\label{sec_Tlamplus}
We next consider $T_\lam^+$, which we recall is defined by
\beq\label{T_lam_plus}
T_\lam^+ f(x,t) = \sum_{2^j \geq 1/N(\lam)} \int_{\R^n}  f(x-y,t-|y|^2)e^{iP_\lam(y)} 2^{-nj}\phi_j(2^{-j}y) dy
\eeq
for $\lam \in \R^{d-1}$ and $f$ of Schwartz class. The key here is to write each term of the sum in the form of the operator $^{(\eta)}I^\lam_a$, for some suitable bump function $\eta$, scaling parameter $a$, and coefficient parameter $\lam$; we will then use Theorem~\ref{Osc_Thm} to obtain bounds for the $L^p$ norms that are summable in $j$. To do so, recall the auxiliary operator defined in (\ref{I_op_dfn_0})  by
 \beq\label{I_op_dfn}
\,^{(\eta)}I_ {a}^{\lam} f(x,t) = \int_{\R^n} f(x-y,t-|y|^2) e^{iP_{\lam}(y/a)}\frac{1}{a^n}\eta(\frac{y}{a}) dy
   \eeq
for $\lam \in \R^{d-1}$. We then introduce a non-isotropic scaling 
\[
2^j \circ \lam = (2^{jm}\lam_m)_m 
,\]
 where $\lam = (\lam_m)_{2 \leq m \leq d} $. Then $P_{2^j \circ \lam}(y/2^j) = P_{\lam}(y)$ since $P_{\lam}(y) = \sum_{m=2}^d \lam_m p_m(y)$ and the polynomials $p_m(y)$ are homogeneous of degree $m$.
We also note that the norm $N(\lam)$ we introduced earlier is homogeneous with respect to this dilation: 
 \beq\label{Nhom}
 N(2^j \circ \lam) = 2^j N(\lam).
 \eeq

Now upon choosing the bump function $\eta = \phi_j$ and the scaling $a = 2^j$, then the $j$-th summand in (\ref{T_lam_plus}) is precisely $^{(\eta)} I^{2^j \circ \lam}_a$; that is to say,
\beq \label{eq:Tlamplusbound}
T_\lam^+ f(x,t) = \sum_{2^j \geq 1/N(\lam)} \phantom{a}^{(\phi_j)}I^{2^j \circ \lam}_{2^j} f(x,t).
\eeq
In order to  bound the operators ${}^{(\phi_j)}I^{2^j \circ \lam}_{2^j}$ and hence $T_\lam^+$ on $L^p$, we will interpolate the key $L^2$ result of Theorem \ref{Osc_Thm} with the following trivial lemma:

\begin{lemma}\label{lemma_I_triv_Lp}
Let $\{\eta_k\}_{k \in \mathbb{Z}}$ be a family of $C^1$ bump functions with $\|\eta_k\|_{C^1} \leq 1$.
For any $1<p \leq \infty$ and any Schwartz function $f$ on $\R^{n+1}$,
\[ \| \sup_{\bstack{\lam \in \R^{d-1}}{k \in \mathbb{Z}}} | ^{(\eta_k)}I_{2^k}^\lam f | \|_{L^p(\R^{n+1})} \leq A_p \| f\|_{L^p(\R^{n+1})}.\]
\end{lemma}
This lemma is simply a trivial consequence of the fact that $ ^{(\eta_k)}I_{2^k}^\lam f(x,t)$ is majorized pointwise almost everywhere by the maximal Radon transform $\Mcal_{\Rad} f(x,t)$.
The precise interpolation result we now require is:
\begin{cor}
Under the hypotheses of Theorem \ref{Osc_Thm}, for any $1<p< \infty$ there exists $\del=\del(p)>0$ such that for all $r \geq 1$ and $f$ of Schwartz class,
\beq\label{Cor1}
	\| \sup_{\bstack{N(\lam)\geq r}{k \in \Z}}  | ^{(\eta_{k})}I_{2^k}^\lam f (x,t)|\, \|_{L^p(\R^{n+1})} \leq A r^{-\del}  \|f\|_{L^p(\R^{n+1})}. 
	\eeq
	\end{cor}
The $L^2$ case of this statement follows from Theorem \ref{Osc_Thm} immediately upon recalling that when $N(\lam) \geq 1$, then $N(\lam) \leq c\|\lam\|$ for some constant $c$ dependent only on the dimension and the degree $d$, so that  $\{ \lam : N(\lam) \geq r\} \subset \{\lam: c\|\lam \| \geq r \}$. Thus for any $r \geq 1$,
\[ \| \sup_{\bstack{N(\lam)\geq r}{k \in \Z}}  | ^{(\eta_{k})}I_{2^k}^\lam f (x,t)|\, \|_{L^2} 
\leq \|\sup_{\bstack{c\|\lam\| \geq r}{k \in \Z}} |^{(\eta_{k})}I^\lam_{2^k} f(x,t)|\, \|_{L^2} .\]
An application of Theorem \ref{Osc_Thm} then shows that
\begin{eqnarray}
 \|\sup_{\bstack{\|\lam\| \geq r/c}{k \in \Z}} |^{(\eta_{k})}I^\lam_{2^k} f(x,t)|\, \|_{L^2}	
	 	&\leq& \sum_{2^s \geq r/c} \|\sup_{\bstack{2^s \leq \|\lam\| < 2^{s+1}}{k \in \Z}} |^{(\eta_{k})}I^\lam_{2^k} f(x,t)|\, \|_{L^2} \nonumber \\
	&\leq &A \sum_{2^s \geq r/c} 2^{-s\del} \|f\|_{L^2} \nonumber \\
		&\leq &A' r^{-\del}\|f\|_{L^2},\label{IL2}
	\end{eqnarray}
	which proves (\ref{Cor1}) in the case $p=2$.
Once we have the crucial $L^2$ case of (\ref{Cor1}), we may interpolate it with the $L^p$ bound (without decay in $r$) of Lemma \ref{lemma_I_triv_Lp} to conclude that for any $1<p<\infty$ there exists a positive $\del(p)>0$ for which (\ref{Cor1}) holds.

	\xtra{Moreover, for any $\del_0< \del$,
\beq\label{Cor2}
	\| \sup_{\bstack{N(\lam)\geq 1}{k \in \Z}} N(\lam)^{\del_0}  | ^{(\eta)}I_{2^k}^\lam f (x)|\, \|_{L^2} \leq A \| \eta \|_{C^1} \|f\|_{L^2}. \eeq
	To prove (\ref{Cor2}), we note that 
\[ \sup_{\bstack{N(\lam)\geq 1}{k \in \Z}} N(\lam)^{\del_0} | ^{(\eta)}I_{2^k}^\lam f (x)| 
	\leq \sum_{l=0}^\infty 2^{l \del_0}  \sup_{\bstack{N(\lam)\geq 2^l}{k \in \Z}}  | ^{(\eta)}I_{2^k}^\lam f (x)| \]
and hence taking $L^2$ norms and applying (\ref{Cor1}),
\[ \| \sup_{\bstack{N(\lam)\geq 1}{k \in \Z}} N(\lam)^{\del_0} | ^{(\eta)}I_{2^k}^\lam f (x)|  \|_{L^2}
	\leq \sum_{l=0}^\infty 2^{l \del_0}  \| \sup_{\bstack{N(\lam)\geq 2^l}{k \in \Z}} | ^{(\eta)}I_{2^k}^\lam f (x)|  \|_{L^2}
	\leq  \sum_{l=0}^\infty 2^{l \del_0} 2^{-l \del} \|\eta\|_{C^1} \|f\|_{L^2}, \]
and the sum is finite since $\del_0 < \del$.}

We may now treat the operator $T_\lam^+$.\xtra{Previously, I tried to adapt Stein and Wainger's original statements. I wrote:
By definition,
	 \beq\label{Tsum'}
T_\lam^+ f(x,t)  =  \sum_{2^j \geq 1/N(\lam)} {\,}^{(\phi_j)}I_ {2^j}^{2^j \circ \lam} f(x,t) 
	=  \sum_{2^j \geq 1/N(\lam)} N(2^j \circ \lam)^{-\del_0} N(2^j \circ \lam)^{\del_0} {\,}^{(\phi_j)}I_ {2^j}^{2^j \circ \lam} f(x,t) .
 \eeq
Therefore,
\begin{eqnarray}
\sup_\lam |T_\lam^+ f(x,t) |
		&=& \sup_\lam |  \sum_{2^j \geq 1/N(\lam)} N(2^j \circ \lam)^{-\del_0} N(2^j \circ \lam)^{\del_0} {\,}^{(\phi_j)}I_ {2^j}^{2^j \circ \lam} f(x,t)| \nonumber \\
		& \leq & \sup_\lam   \sum_{2^j \geq 1/N(\lam)} N(2^j \circ \lam)^{-\del_0} | N(2^j \circ \lam)^{\del_0} {\,}^{(\phi_j)}I_ {2^j}^{2^j \circ \lam} f(x,t)|\nonumber \\
		& \leq & \sup_\lam   \sum_{2^j \geq 1/N(\lam)} N(2^j \circ \lam)^{-\del_0} \sup_{k \in \Z} | N(2^j \circ \lam)^{\del_0} {\,}^{(\phi_j)}I_ {2^k}^{2^j \circ \lam} f(x,t)| \nonumber \\
		& \leq & \sup_\lam   \sum_{2^j \geq 1/N(\lam)} N(2^j \circ \lam)^{-\del_0} \sup_{\bstack{N(\lam')\geq 1}{k \in \Z}} | N(\lam')^{\del_0} {\,}^{(\phi_j)}I_ {2^k}^{\lam'} f(x,t)| 	. \label{lastline}
 \end{eqnarray}
 It is legitimate to take the additional $\sup$ over $\lam'$ when passing to (\ref{lastline}), because all the terms in the sum are non-negative, and taking the sup over $\lam'$ doesn't change the number of $j$ in the summand, but simply enlarges each summand.
Taking $L^2$ norms, I wrote that:
\begin{eqnarray}
\| \sup_\lam |T_\lam^+ f | \|_{L^2}
	& \leq &  \| \sup_\lam   \sum_{2^j \geq 1/N(\lam)} N(2^j \circ \lam)^{-\del_0} \sup_{\bstack{N(\lam')\geq 1}{k \in \Z}} | N(\lam')^{\del_0} {\,}^{(\phi_j)}I_ {2^k}^{\lam'} f|  \|_{L^2} \\
	& \leq &  \sup_\lam   \sum_{2^j \geq 1/N(\lam)} N(2^j \circ \lam)^{-\del_0} \| \sup_{\bstack{N(\lam')\geq 1}{k \in \Z}} | N(\lam')^{\del_0} {\,}^{(\phi_j)}I_ {2^k}^{\lam'} f|  \|_{L^2} \\
	& \leq &  \sup_\lam   \sum_{2^j \geq 1/N(\lam)} N(2^j \circ \lam)^{-\del_0} A \| \phi_j \|_{C^1} \|f\|_{L^2} \\
	& \leq & \left(  \sup_\lam \sum_{2^j \geq 1/N(\lam)} N(2^j \circ \lam)^{-\del_0} \right) C \|f\|_{L^2} \\
	& = &\left(  \sup_\lam  N(\lam)^{-\del_0} \sum_{2^j \geq 1/N(\lam)} 2^{-j\del_0} \right) C \|f\|_{L^2} \\
	& \leq & C' \|f\|_{L^2}.
	\end{eqnarray}
\cm{But this is wrong:
Here, I thought we could pass the $L^2$ norm inside the supremum over $\lam$, but I do not think that we can. This is because for each $(x,t)$, $\lam(x,t)$ will determine which $j$ we sum over, and since we still have a dependence on $\phi_j$, we cannot pass the $L^2$ norm inside the sup over $\lam$ or inside the sum over $j$. }
	}
%
%
%
%
By (\ref{eq:Tlamplusbound}),
 \beq\label{Tsum}
 \sup_\lam |T_\lam^+ f(x,t)|  
 \leq \sup_\lam \sum_{2^j \geq 1/N(\lam)} | ^{(\phi_j)}I_ {2^j}^{2^j \circ \lam} f(x,t) | .
 \eeq
Now we make the key step that dissociates the coefficients of the phase from the scaling factor of the kernel by noting that we may  bound (\ref{Tsum}) by
\beq\label{Tsum2}
 \sup_\lam |T_\lam^+ f(x,t)|  
 \leq \sup_\lam \sum_{N(2^j \circ \lam) \geq 1} \sup_{k \in \Z} | ^{(\phi_k)}I_ {2^k}^{2^j \circ \lam} f(x,t) |
 \leq  \sum_{l=0}^\infty \sup_{\bstack{2^l \leq N(\lam') \leq 2^{l+1}}{k \in \Z}} | ^{(\phi_{k})}I^{\lam'}_{2^k} f(x,t)|.
 \eeq
 Here we have used the homogeneity relation (\ref{Nhom}).
Taking $L^p$ norms for any $1<p<\infty$, we may then apply Theorem \ref{Osc_Thm} in the form of its consequence (\ref{Cor1}) to conclude that
\begin{eqnarray*}
\| \sup_\lam |T_\lam^+ f(x,t)| \|_{L^p}
	& \leq &  \sum_{l=0}^\infty \| \sup_{\bstack{N(\lam') \geq 2^l}{k \in \Z}} | ^{(\phi_k)}I^{\lam'}_{2^k} f(x,t)| \|_{L^p} \\
	& \leq &  A  \sum_{l=0}^\infty 2^{-\del l} \|f\|_{L^p} \leq  A  \| f\|_{L^p},
	\end{eqnarray*}
thus proving Theorem \ref{thm_Carleson}.

\section{The proof of Theorem \ref{Osc_Thm}}\label{sec_outline_pf}
\subsection{A smoother operator}
To prove Theorem~\ref{Osc_Thm}, we require a smoother variant of $^{(\eta)}I_a^\lam$, which we denote by $^{(\eta)}J_a^\lam$. Since in treating the operator $^{(\eta)}J_a^\lam$ one can allow more general polynomial phases, we will first state a result for a more general operator denoted $^{(\eta)}\tilde{J}_a^\lam$, and then specialize to the case we need.

\begin{thm}\label{J_Osc_Thm_strong}
Fix a dimension $n \geq 1$, and a degree $d \geq 2$. Let
\[
Q_\lam(y) = \sum_{2 \leq |\alpha| \leq d} \lam_{\alpha} y^{\alpha}
\]
be a real-valued polynomial on $\R^n$ that has no linear terms, and write $\lam = (\lam_{\alpha})_{2 \leq |\alpha| \leq d}$ for its coefficient parameter. Also write $\Omega_d$ for the set of all such coefficients $\lam$. For $\lam \in \Omega_d$, define
$$
\|\lam\| = \sum_{2 \leq \al \leq d} |\lam_{\al}|.
$$
Fix a $C^1$ bump function $\eta$ supported in $B_1(\R^n)$, and another $C^1$ bump function $\zeta$ supported in $B_1(\R)$, with $\|\eta\|_{C^1}, \|\zeta\|_{C^1} \leq 1$. Then for $\lam \in \Omega_d$, and $a > 0$, one can define an operator $\tilde{J}_a^{\lam}$, acting on Schwartz functions on $\R^{n+1}$, by
$$
\tilde{J}_a^{\lam} f(x,t) = \int_{\R^{n+1}} f(x-y,t-z) e^{iQ_{\lam}(y/a)} \frac{1}{a^n} \eta (\frac{y}{a}) \frac{1}{a^2} \zeta( \frac{z}{a^2} ) dydz.
$$
Furthermore, there exists some fixed $\del>0$ such that for any any $r \geq 1$, 
\[ \|\sup_{\bstack{k \in \Z, \lam \in \Omega_d,}{r \leq \|\lam\| < 2 r}} |\tilde{J}^\lam_{2^k} f(x,t)|\, \|_{L^2(\R^{n+1})} \leq A r^{-\del}\|f\|_{L^2(\R^{n+1})}.\]
\end{thm}
\xtra{We note that the operator $\tilde{J}_a^{\lam}$ is well-defined acting on functions $f$ in $L^2$, although from now on  we will only consider $f$ of Schwartz class.}
In fact, this is effectively Theorem 1 of Stein and Wainger \cite{SWCarl}. The original stopping-time argument of Stein and Wainger, together with our Proposition~\ref{prop_Kflat}, is sufficient to treat $\tilde{J}_a^\lam$ since it does not exhibit Radon-type behavior. However, for precision we note that our statement of Theorem \ref{J_Osc_Thm} differs from Theorem 1 of \cite{SWCarl} in that the operator $\tilde{J}_a^\lam$ includes a product of bump functions $a^{-n}\eta(y/a) a^{-2}\zeta(z/a^2)$ scaled according to parabolic dilations, and the phase function $P_\lam(y)$ is independent of $z$. These are merely cosmetic differences, which Stein and Wainger's argument can handle with only minute changes. 

In addition, we remark that in Theorem \ref{J_Osc_Thm_strong} we could actually have replaced the bump function $\eta$ by a one-parameter family of bump functions $\eta_k$, as long as the $\eta_k$ are uniformly $C^1$ supported in the unit ball. This follows already from Stein-Wainger's original argument. To see this, one need only note that the key bound of Corollary 4.1  in \cite{SWCarl} depends only on the $C^1$ norm of the bump function, which is uniformly bounded by assumption. (For more details, see also Section \ref{sec_vdC_n2}, where the same phenomenon occurs with respect to the operator $^{(\eta_k)}I_{2^k}^\lam$.)

In our application, we only require a consequence of Theorem \ref{J_Osc_Thm_strong} in the special setting of the restricted class of polynomial phases considered in Theorem~\ref{Osc_Thm}. Thus we fix once and for all polynomials $p_2(y)$, $\dots$, $p_d(y)$ with $p_2(y) \neq C|y|^2$ and define for $\lam = (\lam_2, \dots, \lam_d) \in \R^{d-1}$ the polynomial $P_{\lam}(y)$ as in (\ref{allowable_poly}). Define also $\|\lam\|$ as in (\ref{lam_Rn_def}), and fix a $C^1$ function $\zeta$ supported on $[-1,1] \subset \R$, with $\|\zeta\|_{C^1} \leq 1$, and more importantly $$\int_{\R} \zeta(s) ds = 1.$$ From now on, for any $\lam \in \R^{d-1}$, $a > 0$, any $C^1$ bump function $\eta$ supported in the unit ball, and any Schwartz function $f$ on $\R^{n+1}$, we define
$$
^{(\eta)} J_a^\lam f(x,t) = \int_{\R^{n+1}} f(x-y,t-z) e^{iP_{\lam}(y/a)} \frac{1}{a^n} \eta(\frac{y}{a}) \frac{1}{a^2} \zeta( \frac{z}{a^2} ) dy dz.
$$ 

We will also fix a one-parameter family of bump functions $\{\eta_{k}\}_{k \in \Z}$ on $\R^n$ that are all supported in  $B_1(\R^n)$ and have $C^1$ norms uniformly bounded by 1. Then Theorem~\ref{J_Osc_Thm_strong} implies the existence of a small $\delta > 0$ such that for any Schwartz function $f$ on $\R^{n+1}$ and any $r \geq 1$,
\beq \label{J_Osc_Thm} 
\|\sup_{\bstack{k \in \Z, \lam \in \R^{d-1},}{r \leq \|\lam\| < 2 r}} |^{(\eta_{k})} J^\lam_{2^k} f(x,t)|\, \|_{L^2(\R^{n+1})} \leq A  r^{-\del}\|f\|_{L^2(\R^{n+1})}.
\eeq

In order to compare $^{(\eta)} J^{\lam}_a$ to $^{(\eta)} I_a^\lam$, it is helpful to display the singular support of the operator $^{(\eta)} I_a^\lam$ more explicitly by rewriting the operator in the form 
\[ ^{(\eta)}I_a^\lam f(x,t) =  \int_{\R^n} \int_\R f(x-y,t-z) e^{i P_\lam(y/a)} a^{-n} \eta_{a}(\frac{y}{a}) \del_{z=|y|^2} dydz,\]
where $\del$ is the Dirac delta function.
Temporarily, let $I_a^\lam$ and $J_a^\lam$ denote the kernels of the respective operators; then the key property of this construction is that for every fixed $y \in \R^n$, every $\lam \in \R^{d-1}$ and every $a > 0$,
\beq\label{IJ_cancel}
 \int_{\R} (I_a^\lam (y,z)- J_a^\lam(y,z)) dz=0.
 \eeq
We will exploit this property later.

\subsection{An approximation argument}
In order to define the relevant square function we will use to pass from the singular operator $I^\lam_a$ to the smoother operator $J^\lam_a$, it will be convenient to introduce the Littlewood-Paley decomposition constructed in Section \ref{sec_LP}.
For any $N \geq 1$, and any Schwartz functions $f$, we let
\[  L_N f = \sum_{j=-N}^N \Del_j\tilde{\Del}_j f.\]
The function $L_N f$ is then in the Schwartz class, by our observation in Section~\ref{sec_LP}. In order to avoid convergence issues, we will use the finite sum $L_N f$ in order to approximate $f$ in $L^2(\R^{n+1})$; all bounds will be independent of $N$, and morally speaking readers may regard $L_Nf$ as representing $f$.
 
In order to prove Theorem~\ref{Osc_Thm}, from this point on we fix a sequence $\{\eta_{k}\}_{k \in \mathbb{Z}}$ of $C^1$ bump functions with $\|\eta_{k}\|_{C^1} \leq 1$. We observe that trivially, for any Schwartz function $f$ and any fixed $N$,
\beq\label{III}
  \sup_{\bstack{\|\lam\| \approx r}{k \in \Z}} |^{(\eta_{k})} I_{2^k}^\lam f|
	\leq \sup_{\bstack{\|\lam\| \approx r}{k \in \Z}} |^{(\eta_{k})} I_{2^k}^\lam L_N f| + \sup_{\bstack{\|\lam\| \approx r}{k \in \Z}} |^{(\eta_{k})} I_{2^k}^\lam (f-L_Nf)|.\eeq
(From now on, we write $\|\lam\| \approx r$ as a shorthand for $r \leq \|\lam\| \leq 2r$ when $\lam \in \R^{d-1}$.)
We bound the second term on the right hand side of (\ref{III}) using the following trivial proposition:
\begin{prop}\label{prop_triv}
For any  Schwartz function $g$ on $\R^{n+1}$, 
\[ \| \sup_{\bstack{\|\lam\| \approx r}{k \in \Z}} |^{(\eta_{k})} I_{2^k}^\lam g|  \|_{L^2} \leq A \| g\| _{L^2}.\]
\end{prop}	
This follows immediately from the pointwise estimate
\[\sup_{\bstack{\|\lam\| \approx r}{k \in \Z}} |^{(\eta_{k})} I_{2^k}^\lam g| \leq A \Mcal_{\Rad} (g).
\]
Note that Proposition \ref{prop_triv} is extremely weak since it lacks the decay factor $r^{-\del}$, but it is acceptable when applied to the Schwartz class function $f-L_Nf$, which may be taken to have arbitrarily small $L^2$ norm, since $L_Nf$ converges to $f$ in $L^2$ norm, by (\ref{fDelP}). 

As a result, in order to prove Theorem \ref{Osc_Thm}, it now suffices to prove that there exists a $\del>0$ such that 
	\[ \|  \sup_{\bstack{\|\lam\| \approx r}{k \in \Z}} |^{(\eta_{k})} I_{2^k}^\lam L_N f|  \|_{L^2} \leq A r^{-\del} \|f\|_{L^2},\]
where all constants are independent of $N$.
	
\subsection{Introduction of the square function}\label{sec_sq_fn}
	Now we rigorously define the square function $S_r(f)$ acting on a Schwartz function $f$ on $\R^{n+1}$ by
\beq\label{Sfn_precise}
S_r(f) = 	\left( \sum_{k \in \Z} \left( \sup_{\|\lam\| \approx r} |(^{(\eta_{k})} I_{2^k}^\lam - {^{(\eta_{k})} J_{2^k}^\lam}) L_Nf| \right)^2 \right)^{1/2}.
\eeq
Then 
\beq\label{IJ}
\sup_{\bstack{\|\lam\| \approx r}{k \in \Z}} |^{(\eta_{k})} I_{2^k}^\lam L_Nf| \leq  \sup_{\bstack{\|\lam\| \approx r}{k \in \Z}} |^{(\eta_{k})} J_{2^k}^\lam L_N f|  + S_r (f).
\eeq

We bound the first term on the right hand side of (\ref{IJ}) directly by applying Theorem \ref{J_Osc_Thm} and the inequality $\| L_N(f) \|_{L^2} \leq C \|f \|_{L^2}$ provided by (\ref{LN_norm}), to conclude that 
\[  \| \sup_{\bstack{\|\lam\| \approx r}{k \in \Z}} |^{(\eta_{k})} J_{2^k}^\lam L_N f|  \; \|_{L^2} \leq C r^{-\del} \| f \|_{L^2}.\]
The main result for  the square function $S_r(f)$ is:
\begin{thm}\label{thm_Sf}
Let $\{ \eta_k \}_k$ be a one-parameter family of bump functions supported in the unit ball with $C^1$ norm uniformly bounded by $1$.
Then there exists some fixed $\del>0$ such that for any Schwartz function $f$ on $\R^{n+1}$ and any $r \geq 1$, 
\[ \|S_r (f) \|_{L^2(\R^{n+1})} \leq A r^{-\del} \|f \|_{L^2(\R^{n+1})},\]
uniformly in $N$.
\end{thm}

We will derive this as a consequence of the following key inequality:

\begin{thm}\label{Diff_Osc_Thm}
Let $\{\eta_k\}_k$ be as in Theorem \ref{thm_Sf}. Then there exist positive constants $\ep_0>0$ and $\del_0>0$ such that for any Schwartz function $F$ on $\R^{n+1}$, any $j,k \in \Z$, and any $r \geq 1$,
\[ \| \sup_{\|\lam\| \approx r} |(^{(\eta_{k})}I_{2^k}^\lam - ^{(\eta_{k})}J_{2^k}^\lam) \Delta_j F| \, \|_{L^2} \leq A r^{-\del_0} 2^{-\ep_0 |j-k|} \|F\|_{L^2}.\]
\end{thm} 

That this is sufficient to prove the boundedness on $L^2$ of the square-function $S_r(f)$ may be seen as follows. Temporarily define
$F_j = \tilde{\Delta}_j f$ and 
\[A_{j,k} (F)=  \sup_{\|\lam\| \approx r} |(^{(\eta_{k})}I_{2^k}^\lam - ^{(\eta_{k})}J_{2^k}^\lam) \Del_j F|.\]
 Fix an $\ep>0$ that satisfies $\ep < \ep_0$, where $\ep_0$ is the constant given in Theorem \ref{Diff_Osc_Thm}.\xtra{Here we will use the fact that we have a finite sum over $|j| \leq N$ so that we can pass this sum outside of the operators without any issues of convergence.}
 Then
\begin{eqnarray*}
S_r(f) & = & \left( \sum_{k \in \Z} \left( \sup_{\|\lam\| \approx r} | (^{(\eta_{k})}I_{2^k}^\lam - ^{(\eta_{k})}J_{2^k}^\lam) \sum_{|j|\leq N}\Del_j F_j| \right)^2 \right)^{1/2}\\
	& \leq &  \left( \sum_{k \in \Z} \left( \sum_{|j| \leq N} \sup_{\|\lam\| \approx r}  |(^{(\eta_{k})}I_{2^k}^\lam - ^{(\eta_{k})}J_{2^k}^\lam) \Del_j F_j| \right)^2 \right)^{1/2}\\
	&  = & \left( \sum_{k \in \Z} \left( \sum_{|j| \leq N} A_{j,k} (F_j) \right)^2 \right)^{1/2} \\
	 & = &  \left( \sum_{k \in \Z} \left( \sum_{|j| \leq N} 2^{-\ep |j-k|} 2^{\ep |j-k|} A_{j,k} (F_j) \right)^2 \right)^{1/2} \\
	 & \leq & \left( \sum_{k \in \Z} \left( \sum_{j \in \Z} 2^{-2\ep |j-k| }\right) \left(\sum_{|j| \leq N} 2^{2\ep |j-k|} (A_{j,k} (F_j))^2 \right) \right)^{1/2} \\
	 & \leq &  C \left( \sum_{k \in \Z} \sum_{|j| \leq N} 2^{2\ep |j-k|} ( \sup_{\|\lam\| \approx r} |(^{(\eta_{k})}I_{2^k}^\lam - ^{(\eta_{k})}J_{2^k}^\lam) \Del_jF_j|)^2  \right)^{1/2}.
	 \end{eqnarray*}
Thus, taking $L^2$ norms and applying Theorem \ref{Diff_Osc_Thm}, we have
\[ \|S_r(f)\|^2_{L^2} \leq C r^{-2\del_0}   \sum_{k \in \Z} \sum_{|j| \leq N} 2^{2\ep |j-k|} 2^{-2\ep_0 |j-k|} \|F_j\|_{L^2}^2  . \]
Now using the fact that $\ep< \ep_0$, we may sum first in $k$ (to obtain a constant coefficient) and then in $j$, using the property (\ref{forwardLP}) to obtain 
\[ \|S_r(f)\|^2_{L^2} \leq C r^{-2\del_0}  \|f\|^2_{L^2},\]
independent of $N$, thus proving Theorem \ref{thm_Sf} and hence Theorem \ref{Osc_Thm}. 
The remainder of the paper focuses on proving Theorem \ref{Diff_Osc_Thm}.

\section{Proof of Theorem \ref{Diff_Osc_Thm} for the difference operator $I_{2^k}^\lam - J_{2^k}^\lam$}\label{sec_diff_op}
\subsection{Division into Cases 1 and 2}
We recall that Theorem \ref{Diff_Osc_Thm} claims that 
there exists a positive constant $\ep_0>0$ and a constant $\del_0 > 0$ such that for any Schwartz function $F$ on $\R^{n+1}$, any $j,k \in \Z$, and any $r \geq 1$,
\beq\label{diff_claim}
 \| \sup_{\|\lam\| \approx r} |(^{(\eta_{k})}I_{2^k}^\lam - {}^{(\eta_{k})}J_{2^k}^\lam) \Delta_j F| \, \|_{L^2} \leq A  r^{-\del_0} 2^{-\ep_0 |j-k|} \|F\|_{L^2}.
\eeq
Note  that now that $k$ is fixed, the bump function $\eta_{k}$ is fixed, hence we will call the bump function $\eta$ and omit the superscript $\eta_k$ on the operators for simplicity.

We will divide the proof of (\ref{diff_claim}) into two cases: Case 1, when $j \geq k$, and Case 2, when $j <k$. The main propositions in these cases are the following:
\begin{prop}[Case 1: $j \geq k$]\label{prop_Case_1}
There exists a positive constant $\del_0>0$ such that for any Schwartz function $F$ on $\R^{n+1}$, $r \geq 1$ and $j,k \in \mathbb{Z}$ with $j \geq k$,
\[ \| \sup_{\|\lam\| \approx r} |(I_{2^k}^\lam - J_{2^k}^\lam) \Delta_j F| \, \|_{L^2}  \leq A r^{-{\del_0}} 2^{-(j-k)} \|F\|_{L^2}.\]
\end{prop}

\begin{prop}[Case 2: $j<k$]\label{prop_Case_2}
There exists a positive constant $\del_0>0$ and a positive constant $\ep_0 >0$ such that for any Schwartz function $F$ on $\R^{n+1}$, $r \geq 1$ and $j,k \in \mathbb{Z}$ with $j < k$,
\begin{eqnarray}
 \| \sup_{\|\lam\| \approx r} | I_{2^k}^\lam  \Delta_j F| \, \|_{L^2} & \leq & A r^{-\del_0} 2^{\ep_0 (j-k)}  \|F\|_{L^2} \label{case2I} \\
  \| \sup_{\|\lam\| \approx r} | J_{2^k}^\lam  \Delta_j F| \, \|_{L^2} & \leq & A r^{-\del_0} 2^{\ep_0(j-k)}  \|F\|_{L^2}. \label{case2J}
\end{eqnarray}
\end{prop}
Combining Propositions \ref{prop_Case_1} and \ref{prop_Case_2}  immediately proves Theorem \ref{Diff_Osc_Thm}.
Note that in each proposition there are two types of decay present: the first is decay in $r$, which indicates the size of the coefficient parameter $\lam$ controlling the phase, and the second is with respect to $|j-k|$.  In order to extract decay in $r$ we will apply a $TT^*$ argument with stopping-times. In order to extract decay with respect to $|j-k|$, we will use either the cancellation property (\ref{IJ_cancel}) of the kernels of $I_{2^k}^\lam$ and $J_{2^k}^\lam$, or the cancellation property (\ref{Del_cancel}) of the $\Del_j$ kernel in order to perform an integration by parts and pull out the desired factor. 

Schematically, we proceed as follows. To prove Proposition \ref{prop_Case_1} for Case 1 ($j \geq k$), we will use the cancellation property (\ref{IJ_cancel}) for $I^\lam_{2^k}-J^\lam_{2^k}$  in order to place a derivative on the $\Delta_j$ factor, thus enabling us to pull out a factor of $2^{-j}$, at the cost of a factor of $2^k$. We will simultaneously perform a $TT^*$ argument in order to extract the decay factor $r^{-\del}$.

To prove Proposition \ref{prop_Case_2} for Case 2 ($j<k$), we no longer need to exploit the cancellation property (\ref{IJ_cancel}), and therefore it is convenient to separate the terms via the trivial upper bound 
\[ \| \sup_{\|\lam\| \approx r} |(I_{2^k}^\lam - J_{2^k}^\lam) \Delta_j F| \, \|_{L^2} 
	\leq \| \sup_{\|\lam\| \approx r} | I_{2^k}^\lam  \Delta_j F| \, \|_{L^2}  + \| \sup_{\|\lam\| \approx r} | J_{2^k}^\lam  \Delta_j F| \, \|_{L^2} ,
	\]
and treat the two terms on the right hand side separately. We will bound each of these terms by interpolation between two bounds of quite different flavors. On the one hand we will prove the following: 
\begin{prop}[Case 2: $j<k$]\label{prop_Case_2a}
 There exists a positive constant $\del_0>0$ such that for any Schwartz function $F$ on $\R^{n+1}$, $r \geq 1$ and $j,k \in \mathbb{Z}$ with $j < k$,
\begin{eqnarray}
 \| \sup_{\|\lam\| \approx r} | I_{2^k}^\lam  \Delta_j F| \, \|_{L^2} & \leq & A r^{-{\del_0}} \|F\|_{L^2} \label{case2aI} \\
  \| \sup_{\|\lam\| \approx r} | J_{2^k}^\lam  \Delta_j F| \, \|_{L^2} & \leq & A r^{-{\del_0}} \|F\|_{L^2}. \label{case2aJ}
\end{eqnarray}
\end{prop}
This we prove via $TT^*$ arguments; we are allowed to focus purely on extracting decay from the oscillation of the exponential factor in the kernel, since we do not require any decay with respect to $j,k$. 
On the other hand, we will prove that:
\begin{prop}[Case 2: $j<k$]\label{prop_Case_2b}
For any Schwartz function $F$ on $\R^{n+1}$, $r \geq 1$ and $j,k \in \mathbb{Z}$ with $j < k$,
\begin{eqnarray}
 \| \sup_{\|\lam\| \approx r} | I_{2^k}^\lam  \Delta_j F| \, \|_{L^2} & \leq & A r^{1/2} 2^{ (j-k)} \|F\|_{L^2} \label{case2bI}\\
  \| \sup_{\|\lam\| \approx r} | J_{2^k}^\lam  \Delta_j F| \, \|_{L^2} & \leq & A 2^{(j-k)} \|F\|_{L^2}. \label{case2bJ}
\end{eqnarray}
\end{prop}

Note that in Proposition \ref{prop_Case_2b} we are even willing to lose by a factor of $r^{1/2}$, as long as we extract a decay factor $2^{(j-k)}$. Intuitively, we use the cancellation property (\ref{Del_cancel}) that $\int_{\mathbb{R}} \Del_j(\theta) d\theta =0$ in order to write $\Del_j$ as a derivative, and integrate by parts to place a derivative on the kernel of $I_{2^k}^\lam$ or $J_{2^k}^\lam$. In the case of $J_{2^k}^\lam$, it is then a straightforward matter to extract a factor of $2^{j-k}$. But in the case of $I_{2^k}^\lam$, the situation is more complicated since the kernel of $I_{2^k}^\lam$ is supported on the paraboloid and is not differentiable in all directions. However, we are able to proceed with a version of this argument by adapting a $TT^*$ argument, which allows us to differentiate in a way that contributes only an allowable singularity. We then lose by a factor of $r$ in the bound for $\| TT^* \|_{L^2}$, which comes from differentiating the oscillatory factor $e^{iP_\lam(y)}$. (We reiterate that in this case the $TT^*$ formulation aids in differentiating, but is not required in order to extract cancellation from the oscillatory phase.) 

Taking the geometric mean of the bounds in Propositions \ref{prop_Case_2a} and \ref{prop_Case_2b}, namely
\[  \| \sup_{\|\lam\| \approx r} | I_{2^k}^\lam  \Delta_j F| \, \|_{L^2}  \leq  (r^{-{\del_0}})^\theta (r^{1/2} 2^{ (j-k)})^{1-\theta} \|F\|_{L^2} ,\]
with a choice of $\theta$ sufficiently close to 1 (such that $\del_0 \theta > (1-\theta)/2$), finally yields the desired result of Proposition \ref{prop_Case_2}.

\subsection{Derivation of the generic kernel}\label{sec_kernel}
We now set the scene for proving Propositions \ref{prop_Case_1}, \ref{prop_Case_2a} and \ref{prop_Case_2b}. 
First, in Proposition \ref{prop_Case_1}, fix $j \geq k$ and temporarily set $a =2^k$ (for notational convenience). Let $\lam: (x,t) \mapsto  (\lam_2, \dots, \lam_d)$ be any measurable stopping-time function from $\R^{n+1}$ to $\R^{d-1}$. 
We define an operator $T$ acting on Schwartz functions $f$ by setting\xtra{We know that this mapping is well-defined on $L^2$ by the arguments extending the definition of $I_\lam$ to measurable functions, and we know that this mapping is bounded on $L^2$ because it can be majored by maximal functions.}
\beq \label{Tfdef}
Tf(x,t) = 
 (I^{\lam(x,t)}_a - J^{\lam(x,t)}_a) \Delta_j f(x,t).
\eeq
Now we apply the method of $TT^*$.
Explicitly,
$$Tf(x,t) 
= \iint_{\R^{n+2}} f(x-y, t-s-\theta) 
e^{iP_{\lam(x,t)}(\frac{y}{a})} \frac{1}{a^n} \eta(\frac{y}{a}) 
\left(\delta_{s=|y|^2} - \frac{1}{a^2} \zeta( \frac{s}{a^2} ) \right) 
\Delta_j(\theta) dy ds d\theta.$$
Formally changing variables $\theta \mapsto \theta-s$, we get
$$Tf(x,t) 
=\iint_{\R^{n+2}} f(x-y, t-\theta) 
e^{iP_{\lam(x,t)}(\frac{y}{a})} \frac{1}{a^n} \eta(\frac{y}{a}) 
\left(\delta_{s=|y|^2} - \frac{1}{a^2} \zeta( \frac{s}{a^2} ) \right) 
\Delta_j(\theta-s) dy ds d\theta.$$
Hence the operator $T^*$ acts on functions $f$ of Schwartz class by
\begin{multline*}
T^*f(x,t) 
=\iint_{\R^{n+2}} f(x+z, t+\omega) 
e^{-iP_{\lam(x+z, t+\omega)}(\frac{z}{a})} \frac{1}{a^n} \eta(\frac{z}{a}) 
\\ 
\cdot \; \left(\delta_{\xi=|z|^2} - \frac{1}{a^2} \zeta( \frac{\xi}{a^2} ) \right) 
\Delta_j(\omega-\xi) dz d\xi d\omega.
\end{multline*}
 It follows that 
\begin{align*}
TT^*f(x,t)= 
& \iint_{\R^{2n+4}} f(x-y+z,t-\theta+\omega) 
e^{iP_{\lam(x,t)}(\frac{y}{a})-iP_{\lam(x-y+z,t-\theta+\omega)}(\frac{z}{a})} \\
& \cdot \; \frac{1}{a^n} \eta(\frac{y}{a}) \frac{1}{a^n} \eta(\frac{z}{a}) \left(\delta_{s=|y|^2} - \frac{1}{a^2} \zeta( \frac{s}{a^2} ) \right)
  \left(\delta_{\xi=|z|^2} - \frac{1}{a^2} \zeta( \frac{\xi}{a^2} ) \right) \\
& \cdot \; 
\Delta_j(\theta-s)  \Delta_j(\omega-\xi) dy ds dz d\xi d\theta d\omega.
\end{align*}
Changing variables by setting $u=y-z$ and letting $\theta \mapsto \theta + \omega$, this becomes
\begin{align*}
TT^*f(x,t)= 
& \iint_{\R^{2n+4}} f(x-u,t-\theta) 
e^{iP_{\lam(x,t)}(\frac{u+z}{a})-iP_{\lam(x-u,t-\theta)}(\frac{z}{a})} \frac{1}{a^n} \eta(\frac{u+z}{a}) \frac{1}{a^n} \eta(\frac{z}{a}) \\
& \cdot \;  \left(\delta_{s=|u+z|^2} - \frac{1}{a^2} \zeta( \frac{s}{a^2} ) \right)  \left(\delta_{\xi=|z|^2} - \frac{1}{a^2} \zeta( \frac{\xi}{a^2} ) \right) \\
& \cdot \; \Delta_j(\theta-s+\omega)  \Delta_j(\omega-\xi) du ds dz d\xi d\theta d \omega.
\end{align*}
The integral in $\omega$ we recall defines the kernel $\uDel_j$ constructed in Section \ref{sec_uDel}; note that
$$\uDel_j(\theta-s+\xi) :=\int_\R \Delta_j(\theta-s+\omega)  \Delta_j(\omega-\xi) d\omega.$$
Then
\begin{multline}\label{TK}
TT^*f(x,t)= 
 \iint_{\R^{2n+3}} f(x-u,t-\theta) 
e^{iP_{\lam(x,t)}(\frac{u+z}{a})-iP_{\lam(x-u,t-\theta)}(\frac{z}{a})} \frac{1}{a^n} \eta(\frac{u+z}{a}) \frac{1}{a^n} \eta(\frac{z}{a}) \\
\cdot \, \left(\delta_{s=|u+z|^2} - \frac{1}{a^2} \zeta( \frac{s}{a^2} ) \right)  \left(\delta_{\xi=|z|^2} - \frac{1}{a^2} \zeta( \frac{\xi}{a^2} ) \right) 
\uDel_j(\theta-s+\xi) du ds dz d\xi  d\theta.
\end{multline}
Now recalling that $a=2^k$, we see that $\uDel_j(\omega) = \frac{1}{a^{2}} \uDel_{j-k} \left(\frac{\omega}{a^2} \right)$ by (\ref{uDel_rescale}). So we have
\beq\label{TK1}
TT^*f(x,t) = \iint_{\R^{n+1}} f(x-u,t-\theta) \frac{1}{a^{n+2}} \Kcal^{\lam(x,t),\lam(x-u,t-\theta)} \left(\frac{u}{a}, \frac{\theta}{a^2}\right) du d\theta 
\eeq
where for each $\nu$, $\mu \in \R^{d-1}$ we define the kernel
\begin{multline}\label{gen_kernel}
\Kcal^{\nu,\mu}(u,\theta) 
= \iint_{\R^{n+2}} e^{iP_{\nu}(u+z)-iP_{\mu}(z)} \eta(u+z) \eta(z) 
\left(\delta_{s=|u+z|^2} - \zeta(s) \right)  \left(\delta_{\xi=|z|^2} - \zeta(\xi) \right)  \\
\cdot \; \uDel_{j-k}(\theta-s+\xi) dz d\xi ds.
\end{multline}
Estimating this kernel is now the main focus of proving Proposition \ref{prop_Case_1}.

Next, in order to set the stage for proving Propositions \ref{prop_Case_2a} and \ref{prop_Case_2b}, we fix $a=2^k$ and fix a measurable stopping-time function $\lam(x,t): \R^{n+1} \maps \R^{d-1}$ and define two linear operators $T_1$ and $T_2$ respectively by 
\begin{eqnarray*}
T_1f(x,t) &=& I_a^{\lam(x,t)} \Del_j f(x,t)\\
T_2f(x,t) &=& J_a^{\lam(x,t)} \Del_j f(x,t).
\end{eqnarray*}
Then for each of $i=1,2$ we may compute as above that
\beq \label{eq:TiTi*defn}
T_iT_i^* f(x,t) =  \iint_{\R^{n+1}} f(x-u,t-\theta) \frac{1}{a^{n+2}} {\,}^{(i)}\Kcal^{\lam(x,t),\lam(x-u,t-\theta)} \left(\frac{u}{a}, \frac{\theta}{a^2} \right) du d\theta,
\eeq
where for $\nu$, $\mu \in \R^{d-1}$ the respective kernels are given by
\begin{multline}\label{gen_kernel_11}  
 ^{(1)}\Kcal^{\nu,\mu}(u,\theta) 
= \iint_{\R^{n+2}} e^{iP_{\nu}(u+z)-iP_{\mu}(z)} \eta(u+z) \eta(z) 
\delta_{s=|u+z|^2} \delta_{\xi=|z|^2} \\
 \cdot \; \uDel_{j-k}(\theta-s+\xi) dz d\xi ds 
 \end{multline}
 and 
 \begin{multline}\label{gen_kernel_22}
 ^{(2)}\Kcal^{\nu,\mu}(u,\theta) 
= \iint_{\R^{n+2}} e^{iP_{\nu}(u+z)-iP_{\mu}(z)} \eta(u+z) \eta(z) 
\zeta(s)   \zeta(\xi)
\uDel_{j-k}(\theta-s+\xi) dz d\xi ds .
\end{multline}
One then needs to estimate these kernels to complete the proofs of Propositions~\ref{prop_Case_2a} and~\ref{prop_Case_2b}; we return to this in Section \ref{sec_Case2}.

\subsection{Proof of Proposition \ref{prop_Case_1} (Case $j \geq k$)}
To prove Proposition \ref{prop_Case_1}, it is sufficient to prove that if $a=2^k$ and $\lam  = (\lam_2, \dots, \lam_d) \colon \R^{n+1} \to \R^{d-1}$ is a measurable stopping-time satisfying $\|\lam(x,t)\| \approx r$, then for $Tf$ defined as in (\ref{Tfdef}),
\beq\label{TTsuff}
 \| TT^* f \|_{L^2} \leq A^2 r^{-2\del_0}2^{-2\ep_0|j-k|} \|f\|_{L^2}.
 \eeq
We split the kernel (\ref{gen_kernel}) of $TT^*$ into $\mathrm{\bf I}+\mathrm{\bf II}$, where $\mathrm{\bf I}$ gives the contribution of $\delta_{s=|u+z|^2}$, and $\mathrm{\bf II}$ gives the contribution of $\zeta(s)$.

\subsubsection{The term $\mathrm{\bf I}$ $(n=2)$}\label{sec_term_I}
 Explicitly, we evaluate the delta function in $\mathrm{\bf I}$ and consider
\begin{align*}
\mathrm{\bf I}
&= \iint_{\R^{n+1}} e^{iP_{\nu}(u+z)-iP_{\mu}(z)} \eta(u+z) \eta(z) 
\left(\delta_{\xi=|z|^2} - \zeta(\xi) \right)  
\uDel_{j-k}(\theta - |u+z|^2 + \xi) dz d\xi.
\end{align*}
We split this integral into two terms, one coming from $\delta_{\xi=|z|^2}$, and the other coming from $\zeta(\xi)$. In the first term we evaluate the delta function to replace $\xi$ by $|z|^2$, and then reintroduce a trivial integration in $\xi$ via the property $1 = \int \zeta(\xi+|z|^2) d\xi$. In the second term we change variables: $\xi \mapsto \xi+|z|^2$. Grouping terms, we get
\begin{multline}\label{I_dfn_n2_2}
\mathrm{\bf I}
= \iint_{\R^{n+1}} e^{iP_{\nu}(u+z)-iP_{\mu}(z)} \eta(u+z) \eta(z) \zeta(\xi + |z|^2) \\
\cdot \;
\left( \uDel_{j-k}(\theta - |u+z|^2 + |z|^2) - \uDel_{j-k}(\theta - |u+z|^2 + |z|^2 + \xi) \right)  
dz d\xi .
\end{multline}

This is now the first point where we restrict our attention temporarily to the $(2+1)$-dimensional case, that is the $n=2$ case, in order to present the main ideas to the reader without unnecessary technical complications. We will return to analyze this term in general dimension $n \geq 2$ in Section \ref{sec_high_dim_termI}. Meanwhile, we will continue to write $\R^n$ with the understanding that $n=2$, in order to aid the reader conceptually.

We would  like to isolate an oscillatory integral within (\ref{I_dfn_n2_2}) that is independent of the Littlewood-Paley factors $\uDel_{j-k}$. 
To this end, we note that $|u+z|^2 -|z|^2 = |u|^2 + 2u\cdot z$. This motivates defining a new variable $\tau = \frac{u \cdot z}{|u|}$. In the case $n=2$, we would thus make the change of variables $z \mapsto (\tau, \sigma)$, where
\beq\label{tau_sig}
\tau = \frac{u \cdot z}{|u|}, \quad \sigma = \frac{u_2 z_1 - u_1 z_2}{|u|}.
\eeq
Note that this change of variables is simply a rotation, and hence has unit Jacobian. 
In particular, since the support of $\eta$ restricts to $|z| \leq 1$,  we also have $|\tau| \leq 1$, $|\sig| \leq 1$.
(In the case $ n \geq 2$ we will need to make a more general change of variables, and we postpone this discussion until Section \ref{sec_high_dim}.)
After this change of variables, we have
\begin{multline*}
|\mathrm{\bf I}| \leq \int_{\R^2} |K^{\nu,\mu}_{\sharp}(u, \tau; \xi)| \chi_{B_1}(\tau)  \\
\cdot \; \left| \uDel_{j-k}(\theta - |u|^2 - 2|u| \tau) - \uDel_{j-k}(\theta - |u|^2 - 2|u| \tau + \xi) \right| d\tau d\xi
\end{multline*}
where
$$
K^{\nu,\mu}_{\sharp}(u, \tau; \xi)
= \int_{\R^{n-1}} e^{iP_{\nu}(u+z)-iP_{\mu}(z)} \eta(u+z) \eta(z) \zeta(\xi + |z|^2) d\sigma ;
$$
here $z$ is defined implicitly by $u, \tau, \sig$\xtra{, and $\nu = \lam(x,t), \mu = \lam(x-u,t-\theta)$}.
From the supports of $\eta$ and $\zeta$, we see that $K^{\nu,\mu}_{\sharp}(u, \tau; \xi)$ has $\xi$ support in $[-2,2]$.  Thus it follows from the mean-value theorem, as recorded in Lemma \ref{lemma_uDel_diff}, that\xtra{technically there should be a factor of $|\xi|$ here from the MVT, but this is bounded by 2 so we leave it out.}
$$
|\mathrm{\bf I}| \leq C 2^{-2(j-k)} \int_{\R^2} |K^{\nu,\mu}_{\sharp}(u, \tau; \xi)| 
\chi_{B_2}(\xi) 2^{-2(j-k)}|\xi| \psi_{j-k}(\theta - |u|^2 - 2|u|\tau) \chi_{B_1}(\tau) d\tau d\xi ,
$$
where we recall from (\ref{psi_normal}) that $\psi_{j-k}(t)$ is a non-negative integrable function on $\R^1$ with $L^1$ norm uniformly bounded, independent of $j-k$.
We also note that  due to the support of $\eta$, $K_\sharp^{\nu,\mu}$ naturally has $u$ support inside $B_2(\R^n).$

We now apply Proposition \ref{prop_Ksharp} (which also specifies $n=2$) to bound $K_\sharp^{\nu,\mu}$ and conclude that there exists $\delta > 0$ and a small set $G^\nu \subset B_2(\R^n)$ with 
$|G^{\nu}| \leq C r^{-\delta} $, and for each $u \in B_2(\R^n)$ a small set $F^\nu_u \subset B_1(\R)$ with $|F^{\nu}_u| \leq C r^{-\delta},$
such that
\beq\label{K_sharp_bd_app}
|K^{\nu,\mu}_{\sharp}(u, \tau; \xi)| \leq C \left[ r^{-\delta} \chi_{B_2}(u) \chi_{B_1}(\tau) + \chi_{G^{\nu}}(u) \chi_{B_1}(\tau) + \chi_{B_2}(u) \chi_{F_u^{\nu}}(\tau) \right].
\eeq
Moreover, these estimates are uniform in $\xi$, as the small sets do not depend on $\xi$, and neither do the bounds.
Hence after applying this in $\mathrm{\bf I}$ and integrating trivially in $\xi$ we obtain
\begin{multline}\label{I_bound} 
|\mathrm{\bf I}| 
\leq C 2^{-2(j-k)} \int_{\R} \left( r^{-\delta} \chi_{B_2}(u) \chi_{B_1}(\tau) + \chi_{G^{\nu}}(u) \chi_{B_1}(\tau) + \chi_{B_2}(u) \chi_{F_u^{\nu}}(\tau) \right) 
\\ 
\cdot \; \psi_{j-k} (\theta - |u|^2 - 2|u| \tau) d\tau .
\end{multline}
We now plug this into (\ref{TK1}), recalling that $\nu = \lam(x,t)$, and see that the contribution of this part of the kernel to $TT^*f(x,t)$ is thus bounded by
\begin{multline}\label{f_chi_int}
 2^{-2(j-k)} \iint_{\R^{n+1}}|f|(x-u, t -\theta) 
 \int_{\mathbb{R}} \frac{1}{a^{n}} 
\left( r^{-\delta} \chi_{B_2}(\frac{u}{a}) \chi_{B_1}(\tau)  \right. \\
\left. + \; \chi_{G^{\lam(x,t)}}(\frac{u}{a}) \chi_{B_1}(\tau) 
+ \chi_{B_2}(\frac{u}{a}) \chi_{F_{( \frac{u}{a})}^{\lam(x,t)}}(\tau) \right) \\
 \quad \cdot \frac{1}{a^2} \psi_{j-k} (\frac{\theta - |u|^2 - 2a|u| \tau }{a^2})  d\tau du d\theta  .
\end{multline}
We now recognize this as a (non-maximal) averaging operator, with a fixed normalization $a=2^k$, to which we will apply the following simple lemma:
\begin{lemma}\label{lemma_average_op1}
Let $\chi(x,y)$ be an integrable function on $\R^m \times \R^m$ such that there are constants $C_0$ and $\lam$ for which
\beq\label{chi_L1_1}
 \|\sup_{x \in \R^m} |\chi(x,y)| \|_{L^1(\R^m(dy))} \leq C_0
 \eeq
and 
\beq\label{chi_L1_2}
\sup_{x \in \R^m} \| \chi(x,y) \|_{L^1(\R^m(dy))} \leq \lam.
 \eeq
Then for any $1\leq p \leq \infty$ there exists a constant $C = C(C_0,p)$ such that for any $f \in L^p(\R^m)$,  
\beq\label{f_conv_bound}
 \| \int_{\R^m} f(x-y) \chi(x,y) dy \|_{L^p(\R^m(dx))} \leq C \lam^{1/p'} \| f \|_{L^p(\R^m)},
 \eeq
where $p'$ is the conjugate exponent to $p$.
\end{lemma}
We delay the proof momentarily, and apply the lemma to (\ref{f_chi_int}) 
with the choice 
\begin{multline*}
 \chi ((x,t),(u,\theta)) 
 =  \int_{\mathbb{R}} \frac{1}{a^{n}} 
\left( r^{-\delta} \chi_{B_2}(\frac{u}{a}) \chi_{B_1}(\tau) 
+ \chi_{G^{\lam(x,t)}}(\frac{u}{a}) \chi_{B_1}(\tau) + \chi_{B_2}(\frac{u}{a}) \chi_{F_{( \frac{u}{a})}^{\lam(x,t)}}(\tau) \right)  \\
\cdot \;
  \frac{1}{a^2} \psi_{j-k} (\frac{\theta - |u|^2 - 2a|u| \tau }{a^2}) d\tau .
 \end{multline*}
To verify condition (\ref{chi_L1_1}), we note that since for all choices of $(x,t)$ we have the small sets $G^{\lam(x,t)} \subset B_2(\R^n)$ and $F_{(\frac{u}{a})}^{\lam(x,t)} \subset B_1(\R)$, then
\[ \sup_{(x,t)\in \R^{n+1}} |\chi((x,t),(u,\theta))| \leq  \int_{\mathbb{R}} \frac{1}{a^{n}} 
 \chi_{B_2}(\frac{u}{a}) \chi_{B_1}(\tau) 
  \frac{1}{a^2} \psi_{j-k} (\frac{\theta - |u|^2 - 2a|u| \tau }{a^2}) d\tau.\]
Taking the $L^1(du d\theta)$ norm of both sides by integrating first in $\theta$, we obtain
\[ \| \sup_{(x,t) \in \R^{n+1}} |\chi((x,t),(u,\theta))| \; \|_{L^1(du d\theta)} \leq C,\]
independent of $j,k$. 

To verify (\ref{chi_L1_2}), for each fixed $(x,t) \in \R^n$ we compute the $L^1(du d\theta)$ norm of $\chi((x,t),(u,\theta))$ by integrating successively in $\theta,\tau,$ and $u$, in that order. Then by the small measure of the sets $G^{\lam(x,t)}$ and $F_{(\frac{u}{a})}^{\lam(x,t)}$ (and in the first term, the factor $r^{-\del}$ out front), we have 
\[ \sup_{(x,t) \in \R^{n+1}} \| \chi((x,t),(u,\theta))\|_{L^1(du d\theta)} \leq Cr^{-\del},\]
independent of $j,k$. 
Thus by Lemma \ref{lemma_average_op1} we may conclude that 
the $L^2$ norm of the portion of $TT^*f$ considered in (\ref{f_chi_int}) is bounded above by 
\beq\label{ITT}
A 2^{-2(j-k)}r^{-\del/2},
 \eeq
which is sufficient for Proposition \ref{prop_Case_1}.

Finally, we briefly prove Lemma \ref{lemma_average_op1}. We first observe that under the hypothesis (\ref{chi_L1_1}) we have a trivial pointwise bound for all $x \in \R^m$:
\[ \left| \int_{\R^m} f(x-y) \chi(x,y) dy \right| \leq \|f \|_{L^\infty (\R^m)}  \sup_{x \in \R^m} \int_{\R^m} |\chi(x,y)| dy  \leq  \lam \| f \|_{L^\infty(\R^m)},\]
which proves the desired bound for $p= \infty$.
Now for $1 \leq p \leq \infty$, under the hypothesis (\ref{chi_L1_2}) we also have the bound 
\begin{eqnarray}
\|  \int_{\R^m} f(x-y) \chi(x,y) dy \|_{L^p(\R^m(dx))}  
	& \leq & \|  \int_{\R^m} |f|(x-y) \left( \sup_{z\in \R^m}  |\chi(z,y)| \right) dy \|_{L^p(\R^m(dx))} \notag \\
	& \leq &   \int_{\R^m} \|f\|_{L^p(\R^m)} \left( \sup_{z \in \R^m} |\chi(z,y)|\right) dy \label{Lptriangle} \\
	& \leq & C_0  \|f\|_{L^p(\R^m)} . \notag
	\end{eqnarray}
Interpolation between this bound and the $L^\infty$ bound gives the desired result of the lemma.
We remark that Lemma \ref{lemma_average_op1} replaces the so-called small set maximal functions used by Stein and Wainger in \cite{SWCarl}. In our setting, we still require the ability to track the $L^2$ norm relative to the size of the small exceptional sets, but we no longer have true maximal functions since the normalization factor $a=2^k$ is now fixed. (Here we again see the advantage of working inside the square function.)

\subsubsection{The term $\mathrm{\bf II}$ (general $n \geq 2$)}\label{sec_term_II}
Next we look at $\mathrm{\bf II}$, which is the contribution to the kernel (\ref{gen_kernel}) from $\zeta(s)$; we may now treat all dimensions $n \geq 2$ with no additional difficulty. Up to a sign,
\begin{align*}
\mathrm{\bf II} 
&= \iint_{\R^{n+2}} e^{iP_{\nu}(u+z)-iP_{\mu}(z)} \eta(u+z) \eta(z) 
\left(\delta_{\xi=|z|^2} - \zeta(\xi) \right) \zeta(s) 
\uDel_{j-k}(\theta - s + \xi) dz d\xi ds.
\end{align*}
We split this into two terms, one coming from $\delta_{\xi=|z|^2}$, another coming from $\zeta(\xi)$. In the first one we integrate the delta function to replace $\xi$ by  $|z|^2$, and change variables $s \mapsto s+|z|^2$. Then we write $1 = \int \zeta(\xi+|z|^2) d\xi$ to reintroduce a trivial integration in $\xi$. In the second term, we change variables $s \mapsto s+|z|^2$, $\xi \mapsto \xi+|z|^2$. We then get
\begin{multline*}
\mathrm{\bf II}
= \iint_{\R^{n+2}} e^{iP_{\nu}(u+z)-iP_{\mu}(z)} \eta(u+z) \eta(z) \zeta(\xi + |z|^2) \zeta(s + |z|^2)\\
\cdot \; \left( \uDel_{j-k}(\theta - s) - \uDel_{j-k}(\theta - s + \xi) \right)  
dz ds d\xi .
\end{multline*}
Hence
$$
|\mathrm{\bf II}| \leq \int_{\mathbb{R}^2} |K^{\lam(x,t),\lam(x-u,t-\theta)}_{\flat}(u; \xi, s)| \left| \uDel_{j-k}(\theta - s) - \uDel_{j-k}(\theta - s + \xi) \right|  ds d\xi,
$$
where
$$
K^{\nu,\mu}_{\flat}(u; \xi, s)
= \int_{\R^n} e^{iP_{\nu}(u+z)-iP_{\mu}(z)} \eta(u+z) \eta(z) \zeta(\xi + |z|^2) \zeta(s + |z|^2) dz.
$$
Note that $K^{\nu,\mu}_{\flat}(u; \xi, s)$ has compact $\xi$ and $s$ support in $B_2(\R)$, and moreover is an $n$-dimensional integral over the full variable $z$. 
We again apply the mean value theorem in the form of Lemma \ref{lemma_uDel_diff}, concluding that 
$$
|\mathrm{\bf II}| \leq 2^{-2(j-k)} \int_{\mathbb{R}^2}  |K^{\nu,\mu}_{\flat}(u; \xi, s)| \chi_{B_2}(\xi) \chi_{B_2}(s)  \psi_{j-k}(\theta - s) ds d\xi,
$$
where $\psi_{j-k}$ is a non-negative integrable function with $L^1$ norm uniformly bounded, independent of $j-k$. 
Now we apply Proposition \ref{prop_Kflat} (general $n \geq 2$) to bound $K_\flat^{\nu,\mu}$, so that
there exists $\delta > 0$ and a small set $G^{\nu} \subset B_2(\R^n)$ with $$|G^{\nu}| \leq C r^{-\delta}$$ such that
$$|K^{\nu,\mu}_{\flat}(u; \xi, s)| \leq C \left[r^{-\delta} \chi_{B_2}(u) + \chi_{G^{\nu}}(u)\right].$$ 
This estimate is uniform in $\xi$ and $s$, so plugging this back into $\mathrm{\bf II}$ and integrating trivially in $\xi$, we get
\[
|\mathrm{\bf II}|  \leq C 2^{-2(j-k)} \left( r^{-\delta} \chi_{B_2}(u) + \chi_{G^{\nu}}(u) \right)  \int_\R  \psi_{j-k}(\theta - s) \chi_{B_2}(s) ds  .
\]
Here for notational convenience, we will temporarily set
\[ \tilde{\psi}_{j-k} (\theta)=  \int_\R  \psi_{j-k}(\theta - s) \chi_{B_2}(s) ds  \]
which is itself an integrable function of $\theta$, with $L^1$ norm uniformly bounded independent  of $j-k$\xtra{, by Young's inequality}.
Recalling that $\nu = \lam(x,t)$, the contribution of the portion $\mathrm{\bf II}$ of the kernel to $TT^*f$ in (\ref{TK1}) is then bounded by
\beq\label{TTf_contribution2}
C2^{-2(j-k)} \iint_{\R^{n+1}} |f(x-u, t -\theta)| 
\frac{1}{a^{n}} 
 \left( r^{-\delta} \chi_{B_2}(\frac{u}{a}) + \chi_{G^{\lam(x,t)}}(\frac{u}{a}) \right) 
\frac{1}{a^2} \tilde{\psi}_{j-k} (\frac{\theta}{a^2})
d\theta du .
\eeq
We again call upon Lemma \ref{lemma_average_op1}, this time with the choice 
\[ 
\chi((x,t),(u,\theta)) = \frac{1}{a^{n}} 
 \left( r^{-\delta} \chi_{B_2}(\frac{u}{a}) + \chi_{G^{\lam(x,t)}}(\frac{u}{a}) \right) 
\frac{1}{a^2} \tilde{\psi}_{j-k} (\frac{\theta}{a^2}).\]
We may verify (\ref{chi_L1_1}) by noting that 
\[ \sup_{(x,t) \in \R^{n+1}} |\chi((x,t),(u,\theta))| \leq \frac{1}{a^n} \chi_{B_2}(\frac{u}{a}) 
\frac{1}{a^2} \tilde{\psi}_{j-k} (\frac{\theta}{a^2}),\]
and the right hand side is uniformly in $L^1(dud\theta)$ with a bounded norm. We may verify (\ref{chi_L1_2}) by noting that because of the factor $r^{-\del}$ out front of the first term of $\chi$, and the small measure of the set $G^{\lam(x,t)}$ in the second term of $\chi$, we have
\[ \sup_{(x,t) \in \R^{n+1}} \| \chi((x,t),(u,\theta)) \|_{L^1(dud\theta)} \leq Cr^{-\del}.\]
Hence by Lemma \ref{lemma_average_op1}, the $L^2$ norm of the portion of $TT^*f$ contributed by (\ref{TTf_contribution2}) is bounded above by $\leq C r^{-\del/2}2^{-2(j-k)}$.
Combining this with (\ref{ITT}), we have proved that when $j \geq k$, 
$$\|TT^*f\|_{L^2} \leq C r^{-\delta/2} 2^{-2(j-k)} \|f\|_{L^2},$$
which proves (\ref{TTsuff}) and hence Proposition \ref{prop_Case_1} with $\del_0 = \del/4$, in dimension $n = 2$.

To conclude the proof of Proposition~\ref{prop_Case_1}, it remains to prove (\ref{TTsuff}) in dimensions $n > 2$; this requires us to make a slightly different change of variables in estimating (\ref{I_dfn_n2_2}), and that is taken up in Section~\ref{sec_high_dim_termI}.

\subsection{Proof of Propositions \ref{prop_Case_2a} and \ref{prop_Case_2b} (Case $j < k$)}\label{sec_Case2}

In proving Propositions \ref{prop_Case_2a} and \ref{prop_Case_2b}, we will need to prove (\ref{case2aI}) and (\ref{case2bI}) for $I^{\lam}_a$, and (\ref{case2aJ}) and (\ref{case2bJ}) for $J^{\lam}_a$. Below in Section \ref{sec_I_2}, we will first prove (\ref{case2aI}) and (\ref{case2bI}) in dimension $n = 2$, since in this dimension the proof will be slightly cleaner. The case $n > 2$ is then deferred to Section \ref{sec_high_dim_2a2b}. Then in Section \ref{sec_J_n}, we prove (\ref{case2aJ}) and (\ref{case2bJ}), which we can do without any additional difficulty in general dimensions $n \geq 2$.

\subsubsection{Proof of (\ref{case2aI}) and (\ref{case2bI}) $(n = 2)$}\label{sec_I_2}
Suppose now we are in dimension $n = 2$. For $j < k$, we will prove (\ref{case2aI}) of Proposition~\ref{prop_Case_2a} by showing that if $\lam  = (\lam_2, \dots, \lam_d) \colon \R^{n+1} \to \R^{d-1}$ is a measurable stopping-time satisfying $\|\lam(x,t)\| \approx r$, then $T_1 f := I_{2^k}^{\lam(x,t)}f(x,t)$ satisfies 
\[ \|T_1T_1^* f(x,t)\|_{L^2} \leq C r^{-\del/2} \|f\|_{L^2},\]
for some small $\del>0$.

We recall from (\ref{gen_kernel_11}) that the kernel relevant to $T_1T_1^*$ is  
\[ ^{(1)}\Kcal^{\nu,\mu} (u,\theta) = \int_{\R^n} e^{iP_\nu(u+z) - iP_\mu(z)} \eta(u+z) \eta(z) \uDel_{j-k}(\theta - |u+z|^2 + |z|^2) dz,\]
where $\nu = \lam(x,t), \mu = \lam(x-u,t-\theta)$.
We make the change of variables $z \mapsto (\tau,\sig)$ as defined in (\ref{tau_sig}), so that 
\beq\label{11K}
 ^{(1)}\Kcal^{\nu,\mu} (u,\theta) = \int_{\R} K_\sharp^{\nu,\mu} (u,\tau)\uDel_{j-k}( \theta - |u|^2 - 2|u| \tau) d\tau,
 \eeq
where 
\[ K_\sharp ^{\nu,\mu}(u,\tau) = \int_{\R^{n-1}} e^{i P_\nu(u+z) - iP_\mu(z) } \eta (u+z) \eta(z) d\sig.\]
We apply the bound of Proposition \ref{prop_Ksharp} ($n=2$) to $K_\sharp^{\nu,\mu}$ to conclude that 
\begin{multline}\label{K_11_0}
| ^{(1)}\Kcal^{\nu,\mu} (u,\theta) | \leq C  \int_\R \left( r^{-\delta} \chi_{B_2}(u) \chi_{B_1}(\tau) + \chi_{G^{\nu}}(u) \chi_{B_1}(\tau) + \chi_{B_2}(u) \chi_{F_u^{\nu}}(\tau) \right) 
\\
\cdot \; \uDel_{j-k} (\theta - |u|^2 - 2|u| \tau) d\tau .
\end{multline}
Inserting this kernel bound into $T_1T_1^*$ via (\ref{eq:TiTi*defn}), we see that
\begin{eqnarray}
| T_1T_1^* f(x,t)| &\leq& \int_{\R^{n+1}} |f|(x-u,t-\theta) \frac{1}{a^{n+2}} | {\,}^{(1)}\Kcal^{\lam(x,t),\lam(x-u,t-\theta)} \left(\frac{u}{a}, \frac{\theta}{a^2}\right) |du d\theta \nonumber \\
	& \leq & C \int_{\R^{n+1}} |f|(x-u,t-\theta) \chi((x,t),(u,\theta)) du d \theta,  \label{T1T1_0}
	\end{eqnarray}
where we have defined $\chi((x,t),(u,\theta)) $ to be the function
\begin{multline*}
\chi((x,t),(u,\theta)) = \int_\R  \left( r^{-\delta} \frac{1}{a^n}\chi_{B_2}(\frac{u}{a}) \chi_{B_1}(\tau) + \frac{1}{a^n}\chi_{G^{\lam(x,t)}}(\frac{u}{a}) \chi_{B_1}(\tau) \right. \\
\left. +\; \frac{1}{a^n} \chi_{B_2}(\frac{u}{a}) \chi_{F_{(\frac{u}{a})}^{\lam(x,t)}}(\tau) \right)  
\left| \frac{1}{a^2}\uDel_{j-k} (\frac{\theta - |u|^2 - 2a|u| \tau}{a^2}) \right| d\tau. 
	\end{multline*}
Since $\uDel_{j-k}$ is uniformly in $L^1$ independent of $j-k$ (Lemma \ref{lemma_uDel}), we may use the same argument that we applied to (\ref{f_chi_int}) to show via Lemma \ref{lemma_average_op1} that $T_1T_1^*f$  has $L^2$ norm majorized by $C r^{-\del/2}$, which proves (\ref{case2aI}) with $\del_0 = \del/4$.

Next we prove (\ref{case2bI}) of Proposition~\ref{prop_Case_2b}. To do so, we still want to use the method of $TT^*$, not because we want to extract decay from the phase, but because we want to introduce the variables $u$ and $\tau$ and integrate by parts in $\tau$ to pick up a decay of $2^{2(j-k)}$. 
So again let $\lam  = (\lam_2, \dots, \lam_d) \colon \R^{n+1} \to \R^{d-1}$ be a measurable stopping-time satisfying $\|\lam(x,t)\| \approx r$, and define $T_1$ as before. We return to equation (\ref{11K}), from whence we note that since $\int \uDel_{j-k} (\tau) d\tau =0$, we may use (\ref{DelDel}) to write
\beq\label{delta_expression_0}
 \uDel_{j-k}(\theta - 2|u|\tau - |u|^2) = -\frac{2^{2(j-k)}}{2|u|} \frac{d}{d\tau}  \left[ \widetilde{\uDel}_{j-k}(\theta - 2|u|\tau-|u|^2) \right],
 \eeq
for the Schwartz function $\widetilde{\uDel}_{j-k}$ we constructed in Lemma \ref{lemma_Del_anti}.
Thus by integration by parts, 
\[  ^{(1)}\Kcal^{\nu,\mu} (u,\theta) = \frac{2^{2(j-k)}}{2|u|}  \int_{\mathbb{R}}  \partial_\tau K_\sharp^{\nu,\mu} (u,\tau)\widetilde{\uDel}_{j-k}(\theta -  2|u| \tau - |u|^2) d\tau.\]
We now note that since $K_\sharp^{\nu,\mu}(u,\tau)$ is supported where $|u|\leq 2$, $|\tau| \leq 1$ and is a smooth function of $\tau$, 
\[| \partial_\tau K_\sharp^{\nu,\mu}(u,\tau)| \leq C r\chi_{B_2}(u) \chi_{B_1}(\tau), \]
where the factor of $r$ comes from bringing down coefficients of size $\| \nu \|, \| \mu \| \approx r$ when differentiating the phase $P_\nu(u+z) - P_\mu(z)$ with respect to $\tau$. 
Hence
\[ |^{(1)}\Kcal^{\nu,\mu} (u,\theta) | \leq C \frac{ r 2^{2(j-k)}}{2|u|}\chi_{B_2}(u) \int_\R \chi_{B_1}(\tau) \widetilde{\uDel}_{j-k}(\theta - 2|u|\tau-|u|^2) d\tau.\]
Thus using (\ref{eq:TiTi*defn}) again, in total
\begin{eqnarray}
| T_1T_1^* f(x,t) |&\leq& \iint_{\R^{n+1}} |f| (x-u,t-\theta) \frac{1}{a^{n+2}}  |{\,}^{(1)}\Kcal^{\lam(x,t),\lam(x-u,t-\theta)} \left(\frac{u}{a}, \frac{\theta}{a^2}\right) |du d\theta \nonumber \\
 	& \leq  & C r 2^{2(j-k)} \iint_{\R^{n+1}} |f|(x-u,t-\theta) \chi((x,t),(u,\theta)) du d\theta,  \label{T1T1_f_2} 
 \end{eqnarray}
 where we have defined 
 \[ 	\chi((x,t),(u,\theta)) = \int_\R \left(\frac{2|u|}{a}\right)^{-1}\frac{1}{a^{n}}\chi_{B_2}(\frac{u}{a}) \chi_{B_1}(\tau)  
	 \frac{1}{a^2} \widetilde{\uDel}_{j-k}(\frac{\theta}{a^2} - 2a|u|\frac{ \tau}{a^2} - |\frac{u}{a}|^2 )d\tau .
	\]
\xtra{We will apply Lemma \ref{lemma_average_op1} to bound the averaging operator (\ref{T1T1_f_2}); we may verify both (\ref{chi_L1_1}) and (\ref{chi_L1_2}) with fixed upper bounds (not necessarily small) by integrating consecutively in $\theta, \tau, u$:
\begin{eqnarray*}
&&\sup_{(x,t) \in \R^{n+1}} \|\chi((x,t),(u,\theta))\|_{L^1(du d\theta)} \\
 & \leq &  \| \sup_{(x,t) \in \R^{n+1}} |\chi((x,t),(u,\theta))| \; \|_{L^1(du d\theta)}  \\
  & = &  \int_{\R^{n+2}} \left(\frac{2|u|}{a}\right)^{-1}\frac{1}{a^{n}}\chi_{B_2}(\frac{u}{a}) \chi_{B_1}(\tau)  
	 \frac{1}{a^2} \widetilde{\uDel}_{j-k}(\frac{\theta}{a^2} - 2a|u|\frac{ \tau}{a^2} - |\frac{u}{a}|^2 )d\tau du d\theta \\
	 & \leq & C \int_{\R^{n}} \left(\frac{2|u|}{a}\right)^{-1}\frac{1}{a^{n}}\chi_{B_2}(\frac{u}{a}) du \cdot \int_\R \chi_{B_1}(\tau)  
	d\tau  \\
	& \leq & C'.
  \end{eqnarray*}
Here we have used the facts that $\widetilde{\uDel}_{j-k}$ is uniformly in $L^1$ (Lemma \ref{lemma_Del_anti}), and  that $|u|^{-1}$ is locally integrable in $\R^n$ for $n \geq 2$.
 It follows from an application of Lemma \ref{lemma_average_op1} to (\ref{T1T1_f_2}) that 
  \beq\label{T1T1_final_2}
   \|T_1T_1^* f\|_{L^2(\R^{n+1})} \leq C  r 2^{2(j-k)},
   \eeq
as desired, in the case $n=2$.}

To bound the averaging operator (\ref{T1T1_f_2}), we only need to bound 
$$\left\|\sup_{(x,t) \in \R^{n+1}} |\chi((x,t),(u,\theta))|\right\|_{L^1(dud\theta)}$$
 and then proceed as in (\ref{Lptriangle}). Now
\begin{eqnarray*}
 & & \left\| \sup_{(x,t) \in \R^{n+1}} |\chi((x,t),(u,\theta))| \;\right\|_{L^1(du d\theta)}  \\
  & = &  \int_{\R^{n+2}} \left(\frac{2|u|}{a}\right)^{-1}\frac{1}{a^{n}}\chi_{B_2}(\frac{u}{a}) \chi_{B_1}(\tau)  
	 \frac{1}{a^2} \widetilde{\uDel}_{j-k}(\frac{\theta}{a^2} - 2a|u|\frac{ \tau}{a^2} - |\frac{u}{a}|^2 )d\tau du d\theta \\
	 & \leq & C \int_{\R^{n}} \left(\frac{2|u|}{a}\right)^{-1}\frac{1}{a^{n}}\chi_{B_2}(\frac{u}{a}) du \cdot \int_\R \chi_{B_1}(\tau)  
	d\tau  \\
	& \leq & C'.
  \end{eqnarray*}
Here we have used the facts that $\widetilde{\uDel}_{j-k}$ is uniformly in $L^1$ (Lemma \ref{lemma_Del_anti}), and  that $|u|^{-1}$ is locally integrable in $\R^n$ for $n \geq 2$.
 It follows 
 that 
  \beq\label{T1T1_final_2}
   \|T_1T_1^* f\|_{L^2(\R^{n+1})} \leq C  r 2^{2(j-k)},
   \eeq
as desired, in the case $n=2$.

We pause momentarily to remark on the necessity of applying the Littlewood-Paley decomposition carefully. First, it was necessary to compute $T_1T_1^*$ in order to exploit the cancellation property $\int_{\R} \Delta(\theta) d\theta = 0$. Indeed, suppose we computed the kernel of $T_1$ directly rather than taking $T_1T_1^*$; then for example when $0 = j < k$ we would observe that $T_1 f(x,t)$ is given by
$$
T_1 f(x,t) = \iint_{\R^{n+1}} f(x-u,t-\theta) e^{iP_{\lam(x,t)}(\frac{u}{2^k})} \frac{1}{2^{nk}} \eta(\frac{u}{2^k}) \Delta(\theta-|u|^2) du d\theta.
$$ 
It is not clear how one can gain a factor $2^{-k}$ from this, from the cancellation property of $\Delta$ directly.

Next, suppose we had chosen a different Littlewood-Paley projection, say an $(n+1)$-dimensional projection of the form
$$
P_j f(x,t) 
= \iint_{\R^{n+1}} f(x-y,t-s) \frac{1}{2^{(n+2)j}} P(\frac{y}{2^j},\frac{s}{2^{2j}}) dy ds,
$$  
with an appropriate function $P$.
With this projection in mind, if we then set $T_1 f(x,t) = I^{\lam(x,t)}_a P_j f(x,t)$ instead, one would have difficulty gaining decay in $|j-k|$ when one tries integrating by parts in the kernel of $T_1 T_1^*$. For instance, in the case $0 = j < k$ the kernel of $T_1T_1^*$ can be read off via
\begin{align} 
 T_1 T_1^* f(x,t) 
=& \int_{\R^{n+1}}  f(x-v,t-\theta)  \int_{\R^{n+1}} \int_{\R} e^{iP_{\lam(x,t)}(\frac{u+z}{a}) - i  P_{\lam(x-v,t-\theta)}(\frac{z}{a})} \eta(\frac{u+z}{a}) \eta(\frac{z}{a}) \notag \\
&\qquad \cdot (P*P)(v-u,\theta-|u|^2-2u \cdot z)  du dz dv d\theta. \label{T1T1*temp}
\end{align}
Assume $\int_{\R^{n+1}} P(y,s) dyds = 0$. Then one can write $P(y,s)$, or indeed $(P*P)(y,s)$, as
$$
\frac{\partial \tilde{P}^{(1)}}{\partial y_1} (y,s) + \cdots + \frac{\partial \tilde{P}^{(n)}}{\partial y_n} (y,s) + \frac{\partial \tilde{P}^{(n+1)}}{\partial s} (y,s)
$$
where $\tilde{P}^{(1)}, \dots, \tilde{P}^{(n+1)}$ are some Schwartz functions. If we apply the above expression for $(P*P)(y,s)$, then 
\[(P*P)(v-u,\theta-|u|^2-2u \cdot z) =
\left[\frac{\partial \tilde{P}^{(1)}}{\partial y_1} + \cdots + \frac{\partial \tilde{P}^{(n)}}{\partial y_n} + \frac{\partial \tilde{P}^{(n+1)}}{\partial s} \right](v-u,\theta-|u|^2-2u \cdot z).\]
It is not clear how one can plug this into (\ref{T1T1*temp}), and integrate by parts to gain a factor $a^{-1} = 2^{-k}$. So the choice of a good Littlewood-Paley projection here is very important; our choice $\Delta_j$ works for this particular purpose.

\subsubsection{Proof of (\ref{case2aJ}) and (\ref{case2bJ}) (general $n \geq 2$)} \label{sec_J_n}
We now turn to (\ref{case2aJ}) of Proposition~\ref{prop_Case_2a}, for which we will show that if $\lam  = (\lam_2, \dots, \lam_d) \colon \R^{n+1} \to \R^{d-1}$ is a measurable stopping-time satisfying $\|\lam(x,t)\| \approx r$, then for $T_2 f := J_{2^k}^{\lam(x,t)}f(x,t),$ 
\[ \|T_2T_2^* f(x,t)\|_{L^2} \leq C r^{-2\del_0} \|f\|_{L^2}\]
for some small $\del_0>0$. 
Here we can work in general dimension $n \geq 2$. We recall from (\ref{gen_kernel_22}) that the kernel relevant to $T_2T_2^*$ is 
\begin{eqnarray*}
	 ^{(2)}\Kcal^{\nu,\mu}(u,\theta) 
	&=& \int_{\R^{n+2}} e^{iP_{\nu}(u+z)-iP_{\mu}(z)} \eta(u+z) \eta(z) 
	\zeta(s)   \zeta(\xi) \uDel_{j-k}(\theta-s+\xi) dz d\xi ds \\
	& = & \int_{\R^2} K_\flat^{\nu,\mu} (u) \zeta(s) \zeta(\xi) \uDel_{j-k}(\theta-s+\xi)  d\xi ds,
	\end{eqnarray*}
where 
\[ K_\flat^{\nu,\mu}(u) = \int_{\R^n} e^{iP_{\nu}(u+z)-iP_{\mu}(z)} \eta(u+z) \eta(z)  dz.\]
We apply Proposition \ref{prop_Kflat} and conclude that there exists a $\del>0$ and a small set $G^\nu\subset B_2(\R^n)$ with $|G^\nu| \leq r^{-\del}$ so that
\[  |^{(2)}\Kcal^{\nu,\mu}(u,\theta)| \leq C (r^{-\del}\chi_{B_2}(u) + \chi_{G^\nu}(u)) \int_{\R^2} \left| \zeta(s) \zeta(\xi) \uDel_{j-k}(\theta-s+\xi) \right| d\xi ds . \]
This kernel bound shows that 
\begin{eqnarray*}
 |T_2T_2^* f(x,t) |&\leq & C \iint_{\R^{n+1}} |f|(x-u,t-\theta) \left|\frac{1}{a^{n+2}}  {\,}^{(2)}\Kcal^{\lam(x,t),\lam(x-u,t-\theta)} \left(\frac{u}{a}, \frac{\theta}{a^2}\right) \right|du d\theta \\
 	& \leq &  C \iint_{\R^{n+1}}  |f|(x-u,t-\theta) \chi((x,t),(u,\theta)) du d\theta
	\end{eqnarray*}
where 
\[ \chi((x,t),(u,\theta))  = \int_{\R^2} \frac{1}{a^{n}} \left(r^{-\del}\chi_{B_2}(\frac{u}{a}) + \chi_{G^{\lam(x,t)}}(\frac{u}{a})\right) \zeta(s) \zeta(\xi) 
 \frac{1}{a^2}| \uDel_{j-k}(\frac{\theta}{a^2}-s+\xi) | d\xi ds .\]
As in previous arguments, we may check that $\chi$ satisfies the hypotheses of Lemma \ref{lemma_average_op1}, that is
\[\| \sup_{(x,t) \in \R^{n+1}}  |\chi((x,t),(u,\theta)) | \; \|_{L^1(du d\theta)} \leq C \]
and 
\[  \sup_{(x,t) \in \R^{n+1}} \| \chi((x,t),(u,\theta)) \|_{L^1(du d\theta)} 
\leq Cr^{-\del},\]
since $\zeta$ is in $L^1$ and $\uDel_{j-k}$ is uniformly in $L^1$, independent of $j-k$ (Lemma \ref{lemma_uDel}). 
Applying Lemma \ref{lemma_average_op1}, we may conclude that 
\[\| T_2 T_2^* f\|_{L^2} \leq C r^{-\del/2}  \|f\|_{L^2},\]
 which proves (\ref{case2aJ}) with $\del_0 = \del/4.$

Next we prove (\ref{case2bJ}) of Proposition \ref{prop_Case_2b}. We want to bound $f \mapsto \sup_{\|\lam\|\simeq r} |J^{\lam}_{a} \Delta_j f|$ on $L^2$ where $a = 2^k$. Here we will use the fact that $\int_{\R} \Del_j(\theta) d\theta =0$ in order to place a derivative onto the kernel of $J_{2^k}^\lam$. We will not require a stopping time, nor a $TT^*$ argument, since we are not seeking decay with respect to $r$ and no singular support on the paraboloid is present. We first isolate the kernel. Note that $J^{\lam}_{a}$ is a convolution operator, whose convolution kernel is a (parabolic) dilation of $e^{iP_{\lam}(y)} \eta(y) \zeta(s)$ by $a$. Also, $\Delta_j$ is a convolution operator, whose convolution kernel is a (parabolic) dilation of $\delta_{y=0} \Delta(s)$ by $2^j$, which is the same as the (parabolic) dilation of $\delta_{y=0} \Delta_{j-k}(s)$ by $a=2^k$. Hence
$$
J^{\lam}_{a} \Delta_j f(x,t) = \iint_{\R^{n+1}} f(x-y,t-s) \frac{1}{a^{n+2}} L^{\lam} (\frac{y}{a}, \frac{s}{a^2}) dy ds,
$$
where $L^{\lam}(y,s)$ is the convolution of $e^{i P_{\lam}(y) }\eta(y)\zeta(s)$ with $\delta_{y=0} \Delta_{j-k}(s)$, namely
\[ L^{\lam} (y,s) = e^{i P_{\lam}(y) }\eta(y)  \int_{\R} \zeta(s-\theta) \Del_{j-k} (\theta) d\theta.\]
We now use the fact that $\int \Del_{j-k}(\theta) d\theta =0$ to apply Lemma \ref{lemma_Del_anti} and write
\[ L^{\lam} (y,s) = 2^{2(j-k)}e^{i P_{\lam}(y) }\eta(y) \int_{\R} \zeta(s-\theta) \left( \frac{d}{d\theta} \widetilde{\Del}_{j-k} \right)(\theta) d\theta,\]
with the antiderivative $\widetilde{\Del}_{j-k}$ provided by Lemma \ref{lemma_Del_anti}. Then after integration by parts,
\[ L^{\lam} (y,s) = - 2^{2(j-k)}e^{i P_{\lam}(y) }\eta(y) \int_{\R} \frac{d}{d\theta}\zeta(s-\theta)  \widetilde{\Del}_{j-k} (\theta) d\theta,\]
 so that after applying the trivial bound to $\eta$ we obtain
 \[ |L^{\lam} (y,s)| \leq C  2^{2(j-k)}L(y,s),\]
 say, where we define 
 \[ L(y,s)= \chi_{B_1}(y) |\zeta' * \widetilde{\Del}_{j-k}| (s).\]
 We note that this bound is independent of $\lam$,  and as a result,
 \[ \sup_{\lam} |J_{2^k}^{\lam} \Del_j f|   \leq C 2^{2(j-k)} |f| *  \left( \frac{1}{a^{n+2}}L (\frac{\cdot}{a}, \frac{\cdot}{a^2}) \right) .\]
 We now need only recall that $\|L\|_{L^1(\R^{n+1})} \leq C$ since $\widetilde{\Del}_{j-k}$ is uniformly in $L^1$ (Lemma \ref{lemma_Del_anti}).
 Thus by a further application of Young's inequality,
 \[ \| \sup_{\lam} |J_{2^k}^{\lam} \Del_j f| \; \|_{L^2(\R^{n+1})} \leq A 2^{2(j-k)} \|f\|_{L^2(\R^{n+1})},\]
 as required.

\subsection{Summary}
The work of this section has completed the proof of Propositions \ref{prop_Case_1}, \ref{prop_Case_2a} and \ref{prop_Case_2b}, except for the following:
\begin{enumerate}[(a)]
\item estimate of the term $\mathrm{\bf I}$ in (\ref{I_dfn_n2_2}) in dimensions $n > 2$;
\item proof of (\ref{case2aI}) and (\ref{case2bI}) in dimensions $n > 2$; and
\item two key oscillatory integral bounds for the kernels $K_\flat^{\nu,\mu}$ and $K_\sharp^{\nu,\mu}$.
\end{enumerate}
The last point (c) occupies the main body of the next two sections. The first two points (a) and (b) are then quickly resolved in Sections~\ref{sec_high_dim_termI} and \ref{sec_high_dim_2a2b} respectively.

\section{Van der Corput estimates for kernels: Part I}\label{sec_vdC_n2}
We  have now reached the heart of the matter: the van der Corput estimates for the kernels arising in the $TT^*$ arguments throughout Section \ref{sec_diff_op}. These estimates break into two cases: kernels of the form $K_\flat^{\nu,\mu}$, which arise from operators not involving Radon-type behavior, and kernels of the form $K_\sharp^{\nu,\mu}$, which arise from operators exhibiting Radon-type behavior along the paraboloid. 
	
\subsection{Proof of the Van der Corput estimate for $K^{\nu,\mu}_\flat$}\label{sec_vdC_sharp2}
 The van der Corput estimate Proposition \ref{prop_Kflat} for the kernel $K_\flat^{\nu,\mu}$  is implied by Lemma 4.1 of \cite{SWCarl}  (in the case $h=1$ in their notation).
For completeness, we briefly recall the proof here. We recall that $r \leq \| \nu \|, \| \mu \| \leq 2r$. We consider the terms in the phase $Q_\nu(u+z) - Q_\mu(z)$ of $K^{\nu,\mu}_\flat(u)$ that are linear in $z$; precisely, these terms contribute 
\[ \sum_{1 \leq j \leq n} Q_\nu^{(j)}(u) z_j\]
to the phase, where we set
\[Q_\nu^{(j)}(u)= \sum_{2\leq |\al| \leq d} \nu_\al \al_ju^{\al-e_j}.\]
(Here we note in particular that there is no contribution from $Q_\mu(z)$ since the original phase function is assumed to have no linear terms.) 
Thus by Lemma \ref{lemma_Prop2.1}, 
\beq\label{K_flat_bound1}
 |K_\flat^{\nu,\mu}(u)| \leq \left( \sum_{j=1}^n |Q_\nu^{(j)}(u)| \right)^{-1/d}.
 \eeq
It remains to show that the sum of linear coefficients in (\ref{K_flat_bound1}) is bounded below by a small power of $r$, except possibly for $u$ belonging to some small exceptional set. We fix $\rho>0$, to be specified later, and define 
\[ G^\nu = \{ u \in B_2(R^n) : \sum_{j=1}^n |Q_\nu^{(j)}(u)| < \rho \} .\]
Then for $u \in B_2 \setminus G^\nu,$ 
\[ | K_\flat^{\nu,\mu}(u) | \leq \rho^{-1/d}.\]
Moreover, $G^{\nu} \subset \cap_j \{ u \in B_2 : |Q_\nu^{(j)}(u)| < \rho \}$. Hence by Lemma \ref{lemma_Prop2.2}, 
\[| G^\nu| \leq c\rho^{1/d} \min_{1 \leq j \leq n}\left(\sum_{\bstack{2 \leq |\al| \leq d}{|\al-e_j| \geq 1}} |\nu_\al| \al_j \right)^{-1/d}  \simeq c \rho^{1/d}\left( \sum_{j=1}^n \sum_{\bstack{2 \leq |\al| \leq d}{|\al-e_j| \geq 1}} |\nu_\al| \al_j \right)^{-1/d} .\]
Here the key observation is that because there are no linear terms in the original phase polynomial, for every $\al$ such that $\nu_\al \neq 0$  there is at least one $j$ for which $|\al-e_j| \geq 1$, so it follows that
\[ \sum_{j=1}^n \sum_{\bstack{2 \leq |\al| \leq d}{|\al-e_j| \geq 1}} |\nu_\al| \al_j \geq  C \sum_{2 \leq |\al| \leq d} |\nu_\al| \geq C r.\]
Thus $|G^\nu| \leq C (r/\rho)^{-1/d},$ and so choosing $\rho = r^{\del_0}$ with $\del_0<1$ completes the proof of the proposition, with $\del = \min(\frac{1-\del_0}{d},\frac{\del_0}{d})$.

\subsection{Van der Corput estimates for $K_\sharp^{\nu,\mu}$: strategy and examples}
We now turn to the more challenging case of the kernel $K_\sharp^{\nu,\mu}$, which arises in $TT^*$ estimates involving Radon-type behavior; for this we recall that our goal (at least in the case $n=2$) is to prove Proposition \ref{prop_Ksharp} (and later its higher dimensional analogue Proposition \ref{prop_Ksharp_n}).
But in order to aid the reader, rather than proving Proposition \ref{prop_Ksharp} immediately, we first prove the desired bound for $K^{\nu,\mu}_{\sharp}$ in three example cases in dimension $n=2$; these examples already illustrate the main difficulties. We will then develop in Section \ref{sec_high_dim} the full $n$-dimensional argument that will prove both Proposition \ref{prop_Ksharp} and Proposition \ref{prop_Ksharp_n} at once.

 Recall that in dimension $n=2$, we have fixed a set $\Pscr$ of $d-1$ homogeneous polynomials on $\mathbb{R}^2$ with real coefficients, say
\beq\label{pm_dfn}
p_m(y) = \sum_{|\al|=m} c_{\al} y^{\al}, \quad \text{for $m = 2, \dots, d,$}
\eeq
where the coefficients $c_\al$ are fixed once and for all. 
We assume that $p_2(y) \neq C|y|^2$ for any nonzero constant $C$. 
We will also write $p_1(y) \equiv 0$ for convenience (i.e. we set $c_{\al} = 0$ whenever $|\al| = 1$).
Furthermore, for $\nu = (\nu_2, \dots, \nu_d)$, $\mu = (\mu_2, \dots, \mu_d)$, we have
\[
K^{\nu,\mu}_{\sharp} (u,\tau) = \int_{\mathbb{R}} e^{iP_{\nu}(u+z)-iP_{\mu}(z)} \Psi(u,z) d\sigma,
\]
where 
\beq\label{P_dfn}
P_{\nu}(y) = \sum_{m=2}^d \nu_m p_m(y), \quad P_{\mu}(y) = \sum_{m=2}^d \mu_m p_m(y),
\eeq
and $\Psi(u,z)$ is a $C^1$ function supported on $B_2(\R^2) \times B_1(\R)$ with $\|\Psi\|_{C^1} \leq 1$. Here $z$ is defined implicitly in terms of $u$, $\tau$, $\sigma$ by (\ref{tau_sig}), which we now write in the form
\beq z_1 = \frac{u_1\tau + u_2 \sig}{|u|} , \qquad
 z_2 = \frac{u_2\tau - u_1 \sig}{|u|}. \label{z_dfns_z1z2} 
 \eeq
Hence $K^{\nu,\mu}_{\sharp}$ is an oscillatory integral in $\sigma$, whose phase is a polynomial in $\sigma$, and with the trivial bound 
\[ |K^{\nu,\mu}_{\sharp} (u,\tau)| \leq C \chi_{B_2}(u) \chi_{B_1}(\tau).\]
 For fixed $u,\tau,$ in order to apply a van der Corput estimate to $K_\sharp^{\nu,\mu}(u,\tau)$ and deduce a bound at $(u,\tau)$ of the form 
 \[ |K_\sharp^{\nu,\mu}(u,\tau)| \leq A r^{-\del'}  \quad \text{for some $\del'>0$},\] 
we must show that within the phase $P_\nu(u+z) - P_\mu(z)$ (regarded as a polynomial in $\sigma$), the coefficient of at least one monomial $\sig^l$ with $1 \leq l \leq d$ is bounded below by $r^\del$ for some $\del>0$. The coefficients of this polynomial (with respect to $\sigma$) are functions of $u,\tau$. Thus our strategy is to show that for ``most'' $u,\tau$, at least one such coefficient with respect to the variable $\sigma$ is large; the remaining exceptional $u,\tau$ are shown to belong to a sufficiently small exceptional set.

 It is convenient to define the following notation to indicate the norm of the coefficients of $P_\nu(u+z) - P_\mu(z)$ as a polynomial in $\sig$:
\[ \| P_\nu(u+z) - P_\mu(z)\|_\sig := \sum_{1 \leq l \leq d} | C[\sig^l](u,\tau)|,\] 
where $C[\sig^l](u,\tau)$ is defined to be the coefficient of $\sig^l$ when one expands $P_\nu(u+z) - P_\mu(z)$ as a polynomial in $\sigma$. 
We also note that given a polynomial of $\tau$, say $R(\tau)$, we will let $\llbracket R \rrbracket_\tau$ denote the norm of the coefficients (including constant term) of $R$ as a polynomial in $\tau$; note that this norm was considered before in (\ref{isotropic_norm_0}). Analogously, for a polynomial $W(u)$, we will let  $\llbracket W \rrbracket_u$ denote the norm of the coefficients (including constant term) of $W$ as a polynomial in $u$.

Fix a small positive real number $0< \ep_1<1$. Then for a given $(u,\tau)$, if there is an index $l$ for which 
\beq\label{Cl}
|C[\sig^l](u,\tau)| \geq r^{\ep_1},
\eeq
 then 
 \beq\label{Pcoeff}
 \|P_\nu(u+z) - P_\mu(z) \|_\sig\geq r^{\ep_1},
 \eeq
  so that by Lemma \ref{lemma_Prop2.1} we would already know that
\beq\label{big_l}
| K_\sharp^{\nu,\mu}(u,\tau)| 	\leq r^{-\ep_1/d} ,
 \eeq
which would be sufficient for our purposes. Our strategy is to reduce to considering the set of $(u,\tau)$ for which 
\beq\label{C_small}
|C[\sig^l](u,\tau) |\leq r^{\ep_1} \qquad \text{for all $l=3,\ldots, d$},
\eeq
and then deduce certain information about the coefficient of $\sigma$ or $\sigma^2$. 
 At this point we recall that the natural ``enemy'' is a phase component that is precisely a multiple of $|y|^2$, that is, the defining function of the paraboloid. Thus we will say that any homogeneous polynomial of the form $C|y|^{2k}$ for some $k \geq 1$ is \emph{parabolic}; polynomials that are not a multiple of a power of $|y|^2$ will be termed \emph{non-parabolic}.

\subsection{Case A example: no parabolic term is present}
In our first example we consider the case where
$\Pscr = \{ p_3(y), p_4(y)\}$ with $p_3(y) = y_1^3$ and $p_4(y) = y_1^4$.
Then the phase polynomial is given by
\[ P_\lam(y) = \lam_4y_1^4 + \lam_3y_1^3.\]
We will use downward induction to reduce our consideration to the coefficient of $\sigma$ in $P_\nu(u+z) - P_\mu(z)$. 

We note that by the relations (\ref{z_dfns_z1z2}) giving $z_1,z_2$ in terms of $\sig,\tau$, we may easily compute
\begin{align*}
u_1+ z_1 = \frac{u_1}{|u|} \left(|u|+\tau\right) + \frac{u_2}{|u|}  \sigma, \\
u_2+ z_2 = \frac{u_2}{|u|} \left(|u|+ \tau\right) - \frac{u_1}{|u|}  \sigma.
\end{align*}
It will also be helpful for later use to note that $|z|^2 = \tau^2 + \sig^2$, while $|u+z|^2 = (|u| + \tau)^2 + \sig^2$.

We recall that by assumption $r \leq \| \nu \|, \|\mu \| \leq 2r$. We compute explicitly that
\[ P_\nu(u+z) - P_\mu(z) = \sum_{l=0}^4 C[\sig^l](u,\tau) \sig^l \]
where
\begin{eqnarray*}
C[\sig^4](u,\tau) & = & (\nu_4 - \mu_4)\frac{u_2^4}{|u|^4} \\
C[\sig^3](u,\tau) & = & \nu_4 4(|u|+\tau)\frac{u_1u_2^3}{|u|^4}  + \nu_3 \frac{u_2^3}{|u|^3} -  \mu_4 4\tau \frac{u_1u_2^3}{|u|^4} - \mu_3  \frac{u_2^3}{|u|^3}  \\
C[\sig^2](u,\tau) & = & \nu_4 6(|u|+\tau)^2\frac{u_1^2u_2^2}{|u|^4}  + \nu_3 3(|u|+\tau) \frac{u_1u_2^2}{|u|^3} -  \mu_4 6\tau^2 \frac{u_1^2u_2^2}{|u|^4} - \mu_3 3 \tau  \frac{u_1u_2^2}{|u|^3} \\
C[\sig](u,\tau) & = & \nu_4 4(|u|+\tau)^3\frac{u_1^3u_2}{|u|^4} + \nu_3 3(|u|+\tau)^2 \frac{u_1^2u_2}{|u|^3} -  \mu_4 4\tau^3\frac{u_1^3u_2}{|u|^4} - \mu_3 3 \tau^2  \frac{u_1^2 u_2}{|u|^3}.
\end{eqnarray*}
(We may clearly disregard the terms that are constant with respect to $\sig$.)

Our strategy is to use downward induction to eliminate the presence of $ \mu = (\mu_3, \mu_4)$ by  approximating $\nu_4 - \mu_4$ and $\nu_3 - \mu_3$ (up to small error) by functions independent of $\mu_3,\mu_4$. We will then be able to represent the coefficient of $\sig$ independently of $\mu_3,\mu_4$ (up to small error) and show that except for a small set of $u,\tau$ the coefficient of $\sig$ is sufficiently large (that is, at least a small power of $r$).

Precisely, we proceed as follows.
Fix $0 < \ep_1 < 1$. We may assume that $|C[\sig^4](u,\tau)| \leq r^{\ep_1}$, since otherwise (\ref{Pcoeff}) and hence (\ref{big_l}) would already be known. Under this assumption,
\beq\label{P14}
|(\nu_4 - \mu_4)\frac{u_2^4}{|u|^4} | \leq r^{\ep_1}.
\eeq
We may next assume that $|C[\sig^3](u,\tau)| \leq r^{\ep_1}$, otherwise (\ref{big_l}) would already be known. 
We re-write $C[\sig^3](u,\tau)$ as
\[ C[\sig^3](u,\tau) = ( \nu_3  - \mu_3)\frac{u_2^3}{|u|^3}+ (\nu_4 -\mu_4) 4\tau \frac{u_1u_2^3}{|u|^4} + \nu_4 4|u|\frac{u_1u_2^3}{|u|^4}   , \]
and under the assumption that this is $O( r^{\ep_1})$ we see that 
\beq\label{P13}
( \nu_3  - \mu_3)\frac{u_2^3}{|u|^3}  =  - (\nu_4 -\mu_4) 4\tau \frac{u_1u_2^3}{|u|^4}  -  \nu_4 4|u|\frac{u_1u_2^3}{|u|^4}  + O(r^{\ep_1}).
\eeq
We then multiply through by $u_2/|u|$ and apply (\ref{P14}) to deduce that 
\beq\label{P13'}
( \nu_3  - \mu_3)\frac{u_2^4}{|u|^4}  =   -  \nu_4 4|u|\frac{u_1u_2^4}{|u|^5}  + O(r^{\ep_1}).
\eeq
In this example it is unnecessary to consider the coefficient of $\sig^2$; instead we turn immediately to the coefficient of $\sig$, which we re-write as 
\begin{multline}\label{P11}
C[\sig](u,\tau) = (\nu_4  -\mu_4) 4\tau^3\frac{u_1^3u_2}{|u|^4}  +( \nu_3 - \mu_3) 3\tau^2 \frac{u_1^2u_2}{|u|^3}\\
+ \nu_4 4(|u|^3 + 3|u|\tau^2 + 3|u|^2 \tau) \frac{u_1^3u_2}{|u|^4}
+ \nu_3 3(|u|^2 + 2|u| \tau) \frac{u_1^2u_2}{|u|^3}.
\end{multline}
Note that we have isolated all presence of $\mu = (\mu_3,\mu_4)$ into the first two terms.
We multiply both sides of the identity by $u_2^3/|u|^3$ and apply (\ref{P14}) and (\ref{P13'}) to the first two terms on the right hand side to see that 
\[ C[\sig](u,\tau) \frac{u_2^3}{|u|^3}= R^{(1)}_{\nu,u}(\tau) + O(r^{\ep_1}),\]
where 
\[ R^{(1)}_{\nu,u}(\tau) =  - \nu_4 \cdot 12 |u| \tau^2 \frac{u_1^3 u_2^4}{|u|^7} +  \nu_4 \cdot 4(|u|^3 + 3|u|\tau^2 + 3|u|^2 \tau) \frac{u_1^3u_2^4}{|u|^7}
+ \nu_3 \cdot 3(|u|^2 + 2|u| \tau) \frac{u_1^2u_2^4}{|u|^6} \]
is a polynomial in $\tau$ with coefficients dependent only on $\nu,u$ (and independent of $\mu$).
We will show that $R_{\nu,u}^{(1)}(\tau)$ is large for almost all $u,\tau$.

The contribution to $R_{\nu,u}^{(1)}(\tau)$ that is constant with respect to $\tau$ may be written as $|u|^{-4}W_\nu^{(1)}(u)$ where we define
\[ W_\nu^{(1)}(u) = 4\nu_4 u_1^3 u_2^4 + 3\nu_3u_1^2 u_2^4.\]
(Note that this is independent of $\mu$.) We fix $\ep_2$ with $\ep_1<\ep_2<1$ and define an exceptional set 
\[ G^\nu = \{ u \in B_2(\R^2) : |W_\nu^{(1)}(u)| \leq r^{\ep_2} \}.\]
For $u \in B_2(\R^2) \setminus G^\nu$, we set 
\[ F_u^\nu = \{ \tau \in B_1(\R) : |R_{\nu,u}^{(1)}(\tau)| \leq C_0 r^{\ep_1} \}\]
for a sufficiently large constant $C_0$; for $u \in G^\nu$ we set $F_u^{\nu} = \emptyset$.

For all $(u,\tau) \in B_2(\R^2) \times B_1(\R)$ with $u \not\in G^\nu,$ $\tau \not\in F_u^\nu$, we have
\[ |R_{\nu,u}^{(1)}(\tau) | \geq C_0 r^{\ep_1},\]
so that for such $(u,\tau)$ we have
\begin{eqnarray*}
\|P_\nu(u+z) - P_\mu(z) \|_\sig 
&\geq & | C[\sig](u,\tau)| \\
&\geq & \left|  \frac{u_2^3}{|u|^3} C[\sig](u,\tau) \right|\\
&  \geq & |R_{\nu,u}^{(1)}(\tau) | - O(r^{\ep_1}) \\
&  \geq & C_0 r^{\ep_1} - O(r^{\ep_1})\\
& \geq & r^{\ep_1},
\end{eqnarray*}
if $C_0$ is sufficiently large. It then follows from the van der Corput estimate of Lemma \ref{lemma_Prop2.1} that 
\[ |K_{\sharp}^{\nu,\mu}(u,\tau)| \leq C r^{-\ep_1/4} \quad \text{if $u \not\in G^\nu$ and $\tau \not\in F_u^\nu$}.\]

On the other hand, the exceptional sets are small. Indeed $G^\nu$ is the set of $u\in B_2(\R^2)$ where $W_\nu^{(1)}(u)$ is small; since the polynomial $W_\nu^{(1)}(u)$ is visibly a sum of two monomials of distinct total degrees, we see that 
\[ \llbracket W_\nu^{(1)}(u) \rrbracket_u \geq |\nu_3| + |\nu_4| \geq r,\]
which implies by Lemma \ref{lemma_Prop2.2_0} that 
\[ |G^\nu| \leq C \left( \frac{r^{\ep_2}}{r} \right)^{1/4}  = Cr^{-(1-\ep_2)/4}.\]
Furthermore, if $u \in B_2 \setminus G^\nu$, then using the notation $\llbracket \cdot \rrbracket_\tau$ to denote the isotropic norm of the coefficients of $R_{\nu,u}^{(1)}(\tau)$ as a polynomial in $\tau$, we have
\[ \llbracket R_{\nu,u}^{(1)}(\tau) \rrbracket_\tau \geq |u|^{-4} |W_\nu^{(1)}(u)| \geq C  |W_\nu^{(1)}(u)|  \geq C r^{\ep_2}.\]
Thus if $u \in B_2 \setminus G^\nu$, it follows from Lemma \ref{lemma_Prop2.2_0} that 
\[ |F_u^\nu| \leq C \left( \frac{C_0 r^{\ep_1}}{r^{\ep_2}} \right)^{1/2} \leq C' r^{-(\ep_2 - \ep_1)/2}.\]
This suffices to prove Proposition \ref{prop_Ksharp} in this particular case.

\subsection{Case B1 example: a parabolic term is present and $p_2(y) \con 0$}
For our second example we consider the case where $
\Pscr = \{p_4(y)\}$, with $p_4(y) = |y|^4$.
In this case
\[ P_\lam(y)  = \lam_4 |y|^4,\]
a purely parabolic polynomial.
In this case we will see that it suffices to reduce our consideration to the coefficient of $\sigma^2$. We recall that by assumption $r \leq \| \nu \|, \| \mu \| \leq 2r$.
We use the identities $|u+z|^2 = (|u|+\tau)^2 + \sig^2$ and $|z|^2 = \tau^2 + \sig^2$ to compute explicitly that
\[ P_\nu(u+z) - P_\mu(z) = \sum_{l=0}^4 C[\sig^l](u,\tau) \sig^l \]
where
\begin{eqnarray*}
C[\sig^4](u,\tau) & = & \nu_4 - \mu_4 \\
C[\sig^3](u,\tau) & = & 0 \\
C[\sig^2](u,\tau) & = & \nu_4 2 (|u|+\tau)^2 - \mu_4 2 \tau^2 \\
C[\sig](u,\tau) & = & 0.
\end{eqnarray*}
Here clearly we cannot reduce our consideration to the coefficient of $\sig$, which is identically zero; we instead we use downward induction to eliminate the presence of $\nu_4 - \mu_4$ in the coefficient of $\sig^2$.

Fix $0 < \ep_1 < 1$. We may assume that $|C[\sig^4](u,\tau)| \leq r^{\ep_1}$, otherwise (\ref{Pcoeff}) and hence (\ref{big_l}) would already be known. Under this assumption,
\beq\label{P24}
|\nu_4 - \mu_4 | \leq r^{\ep_1}.
\eeq
We next consider the coefficient of $\sig^2$, which we re-write via (\ref{P24}) as 
\[ C[\sig^2](u,\tau) = (\nu_4 - \mu_4)2\tau^2 + \nu_4 2(|u|^2 + 2|u|\tau) =   R_{\nu,u}^{(2)}(\tau)  + O(r^{\ep_1}).\]
where 
\[ R_{\nu,u}^{(2)}(\tau) =   \nu_4 2(|u|^2 + 2|u|\tau) .\]
We note that the constant term in $R_{\nu,u}^{(2)}(\tau) $ with respect to $\tau$ is 
\[ W_\nu^{(2)}(u) = 2\nu_4 |u|^2.\]
We fix $\ep_2$ with $\ep_1< \ep_2 <1$ and set
\[ G^\nu  = \{ u \in B_2 (\R^2) : |W_\nu^{(2)}(u)| \leq r^{\ep_2} \},\]
and for $u \in B_2 \setminus G^\nu$ we set 
\[ F_u^\nu = \{ \tau \in B_1(\R): |R_{\nu,u}^{(2)}(\tau)| \leq C_0 r^{\ep_1} \}\]
for some large absolute constant $C_0$ to be determined later. If $u \in G^\nu$, we set $F_u^\nu = \emptyset$. 
For all $(u,\tau) \in B_2(\R^2) \times B_1(\R)$ with $u \not\in G^\nu,$ $\tau \not\in F_u^\nu$, we have
\[ |R_{\nu,u}^{(2)}(\tau) | \geq C_0 r^{\ep_1},\]
so that for such $(u,\tau)$
\begin{eqnarray*}
\|P_\nu(u+z) - P_\mu(z) \|_\sig 
&\geq & | C[\sig^2](u,\tau)| \\
&  \geq & |R_{\nu,u}^{(2)}(\tau) | - O(r^{\ep_1}) \\
&  \geq & C_0 r^{\ep_1} - O(r^{\ep_1})\\
& \geq & r^{\ep_1},
\end{eqnarray*}
if $C_0$ is sufficiently large. It then follows from the van der Corput estimate of Lemma \ref{lemma_Prop2.1} that 
\[ |K_{\sharp}^{\nu,\mu}(u,\tau)| \leq C r^{-\ep_1/4} \quad \text{if $u \not\in G^\nu$ and $\tau \not\in F_u^\nu$}.\]

It simply remains to verify that the exceptional sets are small. 
Since $G^\nu$ is the set of $u\in B_2$ where $W_\nu^{(2)}(u)$ is small and the polynomial $W_\nu^{(2)}(u)$ is visibly a homogeneous polynomial with coefficient $ \nu_4$, we see that 
\[ \llbracket W_\nu^{(2)}(u) \rrbracket_u \geq  |\nu_4| \geq r,\]
which implies by Lemma \ref{lemma_Prop2.2_0} that 
\[ |G^\nu| \leq C \left( \frac{r^{\ep_2}}{r} \right)^{1/2}  = Cr^{-(1-\ep_2)/2}.\]
Furthermore, if $u \in B_2 \setminus G^\nu$, then 
\[ \llbracket R_{\nu,u}^{(2)}(\tau) \rrbracket_\tau \geq |W_\nu^{(2)}(u)| \geq r^{\ep_2}.\]
Thus if $u \in B_2 \setminus G^\nu$, it follows from Lemma \ref{lemma_Prop2.2_0} that 
\[ |F_u^\nu| \leq C \left( \frac{C_0 r^{\ep_1}}{r^{\ep_2}} \right) \leq C' r^{-(\ep_2 - \ep_1)}.\]
This suffices to prove Proposition \ref{prop_Ksharp} in the case under consideration.

\subsection{Case B2 example: a parabolic term is present and $p_2(y) \not\con 0$}
For our third and final example we consider 
$\Pscr = \{p_4(y), p_2(y)\}$ where
$p_4(y) = |y|^4$, $p_2(y) = y_1^2$. Then
\[ P_\lam(y) = \lam_4 |y|^4 + \lam_2 y_1^2.\] This 
 requires a hybrid argument that considers the coefficients of both $\sig$ and $\sig^2$.
We again recall that $r \leq \| \nu \| , \| \mu \| \leq 2r$ and compute that 
\[ P_\nu(u+z) - P_\mu(z) = \sum_{l=0}^4 C[\sig^l](u,\tau) \sig^l \]
where
\begin{eqnarray*}
C[\sig^4](u,\tau) & = & \nu_4 - \mu_4 \\
C[\sig^3](u,\tau) & = & 0 \\
C[\sig^2](u,\tau) & = & 2 \nu_4 (|u| + \tau)^2 + \nu_2  \frac{u_2^2}{|u|^2}  - 2\mu_4 \tau^2 - \mu_2\frac{u_2^2}{|u|^2} \\
C[\sig](u,\tau) & = & 2 \nu_2 (|u| + \tau) \frac{u_1u_2}{|u|^2}- 2 \mu_2 \tau \frac{u_1u_2}{|u|^2}.
\end{eqnarray*}

Fix $0<\ep_1<1$. We may assume that $|C[\sig^4](u,\tau)| \leq r^{\ep_1}$ since otherwise (\ref{big_l}) would be known, and deduce that
\beq\label{P34}
|\nu_4 - \mu_4 | \leq r^{\ep_1}.
\eeq
We re-write the coefficient of $\sig^2$ as 
\[ C[\sig^2](u,\tau) =  (\nu_4 - \mu_4) 2\tau^2 + \nu_42(|u|^2 + 2|u| \tau) +( \nu_2 - \mu_2) \frac{u_2^2}{|u|^2} .
	   \]
We apply (\ref{P34})  to conclude that 
\beq\label{P22}
 C[\sig^2](u,\tau) = \nu_4 2(|u|^2 + 2|u| \tau) +( \nu_2 - \mu_2) \frac{u_2^2}{|u|^2} + O(r^{\ep_1}).
 \eeq
This still includes dependence on $\nu_2- \mu_2$, which we will eliminate by considering the coefficient of $\sig$. 

We may assume that $|C[\sig](u,\tau)| \leq r^{\ep_1}$, since otherwise (\ref{big_l}) would be known; this allows us to conclude that
\beq\label{P12}
(\nu_2 - \mu_2) 2 \tau \frac{u_1u_2}{|u|^2} = -  \nu_2 |u| \frac{u_1u_2}{|u|^2} + O(r^{\ep_1}).
\eeq
We multiply (\ref{P22}) by $2 \tau \frac{u_1}{|u|}$ and apply (\ref{P12}) to deduce that 
\[ C[\sig^2](u,\tau) \cdot 2 \tau \frac{u_1}{|u|} =  R_{\nu,u}^{(2)}(\tau) + O(r^{\ep_1}),\]
where we define 
\[ R_{\nu,u}^{(2)}(\tau) = \nu_4 \cdot 4 \tau  \frac{u_1}{|u|}(|u|^2 + 2|u| \tau)
	 - \nu_2 |u|(\frac{u_1u_2^2}{|u|^3}) .\]
	 
	 Our assumption is that $\|\nu\| = |\nu_4| + |\nu_2| \approx r$. 
We now must break into two further cases, depending on whether $|\nu_4| \geq r/2$ or $|\nu_2 | \geq r/2$. In the first case, we single out the coefficient of the linear term in $\tau$ in $R_{\nu,u}^{(2)}(\tau)$, which takes the form $|u|^{-1} W_\nu^{(2,1)}(u)$, where we define the polynomial 
\[ W_\nu^{(2,1)}(u) = 4\nu_4 u_1|u|^2 = 4\nu_4 (u_1^3 + u_1u_2^2).\]

We fix $\ep_2$ with $\ep_1< \ep_2 <1$ and set
\[ G^\nu  = \{ u \in B_2(\R^2) : |W_\nu^{(2,1)}(u)| \leq r^{\ep_2} \},\]
and for $u \in B_2 \setminus G^\nu$ we set 
\[ F_u^\nu = \{ \tau \in B_1 (\R): |R_{\nu,u}^{(2)}(\tau)| \leq C_0 r^{\ep_1} \}\]
for some large absolute constant $C_0$ to be determined. If $u \in G^\nu$, we set $F_u^\nu = \emptyset$. 
For all $(u,\tau) \in B_2 (\R^2) \times B_1(\R)$ with $u \not\in G^\nu,$ $\tau \not\in F_u^\nu$, 
\[ |R_{\nu,u}^{(2)}(\tau) | \geq C_0 r^{\ep_1},\]
so that for such $(u,\tau)$ we have
\begin{eqnarray*}
\|P_\nu(u+z) - P_\mu(z) \|_\sig 
&\geq & | C[\sig^2](u,\tau)| \\
&\gtrsim & \left| C[\sig^2](u,\tau) \cdot 2\tau \frac{u_1}{|u|} \right| \\
&  \geq & |R_{\nu,u}^{(2)}(\tau) | - O(r^{\ep_1}) \\
&  \geq & C_0 r^{\ep_1} - O(r^{\ep_1})\\
& \geq & r^{\ep_1},
\end{eqnarray*} 
if $C_0$ is sufficiently large.
 It then follows from the van der Corput estimate of Lemma \ref{lemma_Prop2.1} that 
\[ |K_{\sharp}^{\nu,\mu}(u,\tau)| \leq C r^{-\ep_1/4} \quad \text{if $u \not\in G^\nu$ and $\tau \not\in F_u^\nu$}.\]

It simply remains to verify that the exceptional sets are small. 
Since $G^\nu$ is the set of $u\in B_2$ where $W_\nu^{(2,1)}(u)$ is small, and 
\[ \llbracket W_\nu^{(2)}(u) \rrbracket_u \geq  |\nu_4| \geq r/2,\]
we see by Lemma \ref{lemma_Prop2.2_0} that 
\[ |G^\nu| \leq C \left( \frac{r^{\ep_2}}{r} \right)^{1/2}  = Cr^{-(1-\ep_2)/2}.\]
Furthermore, if $u \in B_2 \setminus G^\nu$, then 
\[ \llbracket R_{\nu,u}^{(2)}(\tau) \rrbracket_\tau \geq |u|^{-1} |W_\nu^{(2,1)}(u)| \geq C  |W_\nu^{(2,1)}(u)|  \geq C r^{\ep_2}.\]
Thus if $u \in B_2 \setminus G^\nu$, it follows from Lemma \ref{lemma_Prop2.2_0} that 
\[ |F_u^\nu| \leq C \left( \frac{C_0 r^{\ep_1}}{r^{\ep_2}} \right)^{1/2} \leq C' r^{-(\ep_2 - \ep_1)/2}.\]

In the remaining case when $|\nu_2| \geq r/2$, we single out the coefficient of the constant term in $\tau$ in 
$R_{\nu,u}^{(2)}(\tau)$, setting
\[ W_\nu^{(2,0)}(u) = - \nu_2 u_1u_2^2,\]
and proceed analogously.
This suffices to prove Proposition \ref{prop_Ksharp} in the case under consideration.

\section{Van der Corput estimates for kernels: Part II}\label{sec_high_dim}
\subsection{Partition of unity and change of variables}
We now state and prove a general version of Proposition \ref{prop_Ksharp} for $K_\sharp^{\nu,\mu}$ in all dimensions $n \geq 2$; for this we first require some notation.
In the case $n=2$, we were motivated to make the change of variables in (\ref{tau_sig}) because setting $\tau = u \cdot z /|u|$ captures the behavior of $z$ in the expression 
\beq\label{norm_diff}
 |u+z|^2 - |z|^2 = |u|^2 + 2 u \cdot z = |u|^2 + 2|u| \tau.
 \eeq
We continue to define $ \tau = u \cdot z/ |u|$ in the case of general dimension, but for $n \geq 3$ there is no longer a unique choice of $\sig$ orthogonal to $\tau,$ and we must be more careful.

We now set some notation. For any variable $u \in \R^n$, let $u^{(j)}$ denote the variable in $\R^{n-1}$ that omits the $j$-th coordinate of $u$. Similarly, for any multi-index $\al \in \Z_{\geq 0}^{n}$ we let $\al^{(j)} \in \Z_{\geq 0}^{n-1}$ denote the multi-index omitting $\al_j$.  Given coordinates $(u,z) \in \R^n \times \R^n \setminus \{(0,0)\}$, we may fix a coordinate $1 \leq l \leq n$ where $u_l \ne 0$, and make the change of variables $z \mapsto (\tau,\sig) \in \R \times \R^{n-1}$ defined by
\begin{eqnarray}
\tau &=& \frac{u\cdot z}{|u|} 	\label{tau_rel} \\
\sig & = & 
\frac{u^{(l)} \tau  - |u| z^{(l)} }{u_l}  . \label{sig_rel}
\end{eqnarray}
Correspondingly, $z$ is defined implicitly in terms of $(\tau,\sig)$ by the relations
\begin{eqnarray}
z^{(l)}& =& \frac{u^{(l)} \tau - u_l \sig}{|u|} \label{z_rel1}\\
z_l & = & \frac{u_l \tau + u^{(l)} \cdot \sig}{|u|} \label{z_rel2}.
\end{eqnarray}
The intuition behind these choices is as follows: $\tau$ again captures the behavior of $z$ in the relation (\ref{norm_diff}), while $\sig$ is defined so that (\ref{z_rel1}), which is the higher dimensional analogue of the equation defining $z_2$ in (\ref{z_dfns_z1z2}), holds. Indeed, the relation (\ref{z_rel1}) specifies $\sig$ uniquely, and one can then verify that (\ref{z_rel2}), which is analogous to the equation defining $z_1$ in (\ref{z_dfns_z1z2}), continues to hold.
By explicit computation, the Jacobian associated to this change of variables is $(|u_l|/|u| )^{n-2}$.\xtra{We note that, for example, in the case $k=1$, the relevant Jacobian matrix is:
\[ \frac{dz}{d\tau} = \left( \begin{array}{ccccc}
	\frac{u_1}{|u|} & \frac{u_2}{|u|} & \cdots &  \cdots &\frac{u_n}{|u|}\\
	\frac{u_2}{|u|} &- \frac{u_1}{|u|} & 0 &  \cdots & 0\\
	\frac{u_3}{|u|} & 0 & \cdots & - \frac{u_1}{|u|} & 0\\
	\vdots &&&& \vdots \\
	\frac{u_n}{|u|} &0 & \cdots &  \cdots & -\frac{u_1}{|u|}
	\end{array} 
	\right),
	\]
which has determinant $(-1)^{n-1} \left( -\frac{u_1}{|u|} \right)^{n-2}$.
	}

We choose a partition of unity
\beq\label{partition_1}
1 =  \sum_{l=1}^n W_l( s) 
\eeq
for $s \in \mathbb{S}^{n-1}\subseteq \R^n$ such that for each $1 \leq l \leq n$, $W_l \in C^\infty_c(\mathbb{S}^{n-1})$ and $W_l(s)$ is supported where
\beq\label{zeta_cond}
\left| \frac{s_l}{s} \right| \geq c_0 
\eeq
for some fixed $c_0 >0$. 


In general, given an integral of the form 
\[ \Lcal(u) =  \int_{\R^n} L(u,z) dz \]
with $u \in \R^n$ and $L(u,z)$ supported where $|u| \leq 2$, $|z| \leq 1$, 
we will decompose this as $\Lcal = \sum_{l=1}^n \Lcal_l$ where 
\[ \Lcal_l(u) = W_l \left( \frac{u}{|u|} \right) \int_{\R^n} L(u,z) dz.\]
For each $1 \leq l \leq n$, within the integral defining $\Lcal_l$ we will make the change of variables defined by (\ref{z_rel1}), (\ref{z_rel2}) according to the $l$-th coordinate, in order to obtain
\[ \Lcal_l(u) =  \left( \frac{|u_l|}{|u|} \right)^{n-2} W_l \left( \frac{u}{|u|} \right) \int_{\R} \int_{\R^{n-1}} \tilde{L}(u,\tau; \sig) d\sig  d \tau, \]
in which the function $\tilde{L}(u,\tau;\sig)$ is implicitly defined in terms of $L(u,z)$ by (\ref{z_rel1}) and (\ref{z_rel2}).

In the case of $\R^{2} \times \R$, $\tau$ and $\sig$ enjoyed the convenient relation that $|z|^2 = |\tau|^2 + |\sig|^2$, which allowed us to draw the conclusion that $|\tau|, |\sig| \leq 1$ as long as $|z| \leq 1$. This identity no longer need hold in higher dimensions; instead one can compute that according to the change of variables (\ref{z_rel1}), (\ref{z_rel2}),
\beq\label{sig_angle}
 |z|^2 = z_l^2 + |z^{(l)}|^2 = |\tau|^2 + \frac{ |u^{(l)}|^2 |\sig|^2 \cos^2 \theta_l + u_l^2 |\sig|^2}{|u|^2},
 \eeq
where $\theta_l$ is the angle between $u^{(l)}$ and $\sig$ in $\R^{n-1}$. In the case of $\R^2 \times \R$, this angle always had $\cos \theta = \pm 1$, so we trivially obtained the upper bound $|\sig|\leq 1$ for the range of integration in $\sig$, but in higher dimensions this may not be the case. Throwing away the non-negative cosine term and $|\tau|^2$, and using only $|z| \leq 1$, (\ref{sig_angle}) trivially yields the bound
\beq\label{sig_upper_bound}
 |\sig|  \leq \frac{|u|}{|u_l|}.
 \eeq
By construction, this is a bounded region, because of the key restriction provided by $W_l(u/|u|)$, which requires that $|u_l/u| \geq c_0>0$; this is our motivation for the partition of unity.

We will deploy this partition of unity and change of variables in two places: to treat the term $\mathrm{\bf I}$ in Section \ref{sec_term_I}, and to treat the kernel $^{(1)}\mathcal{K}^{\nu,\mu}$ in Section \ref{sec_I_2}.
In each case, we will require a corresponding van der Corput estimate for a family of oscillatory integrals $K^{\nu,\mu}_{\sharp,l}$ for $1 \leq l \leq n$, each corresponding to a different component in the partition of unity and the relevant change of variables. We state the necessary bounds, a generalization of Proposition \ref{prop_Ksharp}, below.

\begin{prop}[$K_\sharp^{\nu,\mu}$ van der Corput, general $n \geq 2$]\label{prop_Ksharp_n}
Fix a dimension $n \geq 2$ and a degree $d \geq 2$. Let $\Pscr = \{ p_2(y), p_3(y),\ldots, p_d(y)\}$ be a set of real-valued polynomials on $\mathbb{R}^n$, where each $p_j(y)$ is homogeneous of degree $j$, and $p_2(y) \ne C|y|^2$ for any nonzero constant $C$. Let 
$\Lambda = \Lambda (\Pscr) = \{ 2 \leq m \leq d \colon p_m(y) \not\con 0\}.$
 For $\nu = (\nu_2,\dots,\nu_d) \in \R^{d-1}$, let 
 \[ P_{\nu}(y) = \sum_{m=2}^d \nu_m p_m(y)\]
and
\[ \| \nu \| = \sum_{m \in \Lambda} |\nu_m|,\]
and define $P_\mu(y)$ and $\| \mu \|$ similarly for $\mu = (\mu_2, \dots, \mu_d) \in \R^{d-1}$.

Recall the partition of unity given by $W_l$ for $1 \leq l \leq n$. Given a $C^1$ function $\Psi(u,z)$ supported on $B_2 \times B_1 \subset \R^n \times \R^n$, 
define for each $1 \leq l \leq n$ the integral
\[ K^{\nu,\mu}_{\sharp,l}(u, \tau)
= \left(\frac{ |u_l|}{|u|} \right)^{n-2} W_l \left( \frac{u}{|u|} \right) \int_{\R^{n-1}} e^{iP_{\nu}(u+z)-iP_{\mu}(z)} \Psi(u,z)d\sigma, \]
where $z$ is defined implicitly in terms of $u, \tau,\sigma$ by 
\begin{eqnarray}
\tau &=& \frac{u\cdot z}{|u|} 	\label{tau_rel'} \\
\sig & = & 
\frac{u^{(l)} \tau  - |u| z^{(l)} }{u_l}  . \label{sig_rel'}
\end{eqnarray}
Suppose furthermore that 
\[\| \Psi \|_{C^1(\R)} := \sup_{(u,z) \in B_2(\R^n) \times B_1(\R^n)} \left( |\Psi(u,z)| + |\frac{\partial}{\partial \sig} \Psi(u,z)|\right) \leq 1.\]

Then there exists a small constant $\delta > 0$ (depending only on $d$) such that the following holds:
if $\nu$, $\mu$ satisfy
$$
r \leq \|\nu\|,\|\mu\| \leq 2r
$$
for some $r \geq 1$, then there exists a small set $G^\nu \subset B_2(\R^n)$, and for each $u \in B_2(\R^n)$ a small set $F_u^\nu \subset B_1(\R)$, such that
\beq\label{measure_conditions}
|G^\nu| \leq C r^{-\del}, \qquad |F^\nu_u| \leq C r^{-\del} \quad \text{for all $u \in B_2(\R^n)$},
\eeq
and such that for all $1 \leq l \leq n$,
\beq\label{K_sharp_bd_n}
|K^{\nu,\mu}_{\sharp,l}(u, \tau)| \leq C
	\left( r^{-\delta} \chi_{B_2}(u) \chi_{B_1}(\tau) + \chi_{G^{\nu}}(u) \chi_{B_1}(\tau) + \chi_{B_2}(u) \chi_{F_u^{\nu}}(\tau) \right).
\eeq
The choices of  the small sets $G^{\nu}$ and $F^{\nu}_u$ are independent of both $\mu$ and the amplitude $\Psi$.
\end{prop}
We note that due to the support of $W_l(u/|u|)$ we have the upper bound (\ref{sig_upper_bound}) for the support of the integral in $\sigma$, leading to the trivial bound 
\beq\label{K_triv_bd_n}
 |K^{\nu,\mu}_{\sharp,l}(u, \tau)| \leq C c_0^{-(n-1)} \chi_{B_2}(u) \chi_{B_1}(\tau),
 \eeq
for the finite nonzero constant $c_0$.
Thus the import of Proposition \ref{prop_Ksharp_n} is to extract decay in $r$ under the hypotheses of the proposition.

\subsection{Computing the phase}\label{sec_compute_phase_n}
We now prove Proposition \ref{prop_Ksharp_n} in full generality. We will construct for each index $1 \leq l \leq n$ a pair of exceptional sets $G^\nu$ and $F_u^\nu$; taking the union over $n$ of these sets will clearly give exceptional sets that work for all $l$ simultaneously and still satisfy the small measure conditions (\ref{measure_conditions}).
As it will be notationally convenient, we will focus on the case $l=n$, but the argument we present does not depend on this choice in any way beyond notation. In particular, from now on, $\tau$ and $\sig$ will refer to the definitions (\ref{tau_rel'}) and (\ref{sig_rel'}) with the choice $l=n$.

We recall the polynomials $P_{\nu}$, $P_{\mu}$ defined as in Proposition \ref{prop_Ksharp_n}, and compute the phase:
\begin{lemma}\label{lemma_P_phase_n}
The phase $P_{\nu}\left(u+ z\right) - P_\mu(z)$ of $K_{\sharp,n}^{\nu,\mu} (u,\tau)$ is a polynomial in $\sig \in \R^{n-1}$, which we will denote by 
\beq\label{phase_expansion_special_n}
 P_\nu(u+z)  - P_\mu(z) = \sum_{0 \leq |\be| \leq d} C[\sig^\be](u,\tau) \sig^\be,
 \eeq
where for each multi-index $\be$, if $|\be|=l$ then the coefficient $C[\sig^\be](u,\tau)$ is given by 
\beq\label{sig_coeff_n_res}
C[\sig^\be](u,\tau) = \sum_{m =l}^d \left(\nu_m \left(|u|+ \tau\right)^{m-l} - \mu_m \tau^{m-l} \right) B_{m,\be}\left(\frac{u}{|u|}\right) ,
\eeq
in which 
\beq
B_{m,\be} (w) =  \sum_{|\al| = m} c_\al A_{\alpha,\be} (w).
\eeq
Here  $A_{\al,\be}(w)$ is a polynomial in $w$ that is homogeneous of degree $|\al|$, and the coefficients $c_\al$ are fixed once and for all by the choice of the polynomials $\mathscr{P} = \{p_m : m=2,\ldots, d\}$. 
\end{lemma}

To prove Lemma \ref{lemma_P_phase_n}, we define a family of polynomials $A_{\al,\be}(w)$ acting on $w \in B_1( \R^n)$ and parametrized by multi-indices $\al \in \Z^{n}_{\geq 0}$ and $\be \in \Z^{n-1}_{\geq 0}$ (note the differing dimensions). For any multi-index $\al \in \Z^{n}_{\geq 0}$, we specify a polynomial $A_{\al,\be}$ for each $\be \in \Z^{n-1}_{\geq 0}$ with $|\be| \leq |\al|$ by the defining relation
\beq\label{defining_A_rel}
\sum_{|\be| \leq |\al| }A_{\al,\be}(w) \sig^\be = (w_n + w^{(n)} \cdot \sig)^{\al_n} (w^{(n)}  - w_n \sig)^{\al^{(n)}}.
\eeq
It is clear by inspection that with this definition, $A_{\al,\be}$ is a polynomial in $w$ homogeneous of degree $|\al|$. (Moreover, if $|\be| \geq |\al|$ then $A_{\al,\be}$ is the zero polynomial.) Note as well that these polynomials are defined purely in terms of combinatorial coefficients, and are independent of the fixed polynomials $p_m$ and of stopping-times.
For example, in the case of dimension $n=2$, we represent the one-dimensional index $\be \in \Z_{\geq 0}$ by $l$  and compute that for each $1 \leq l \leq |\al|$,
\beq\label{A_dfn_3}
A_{\alpha,l}(w) := \sum_{\substack{0 \leq j \leq \alpha_1 \\ 0 \leq k \leq \alpha_2 \\ j+k = l}} (-1)^j \binom{\alpha_1}{j} \binom{\alpha_2}{k} w_1^{\alpha_1-j+k} w_2^{\alpha_2-k+j}.
\eeq
This is visibly a homogeneous polynomial in $w\in \R^2$ of degree $|\alpha|$. 
\xtra{We would like to express $Q_{\vec{\nu}}\left(u+ z\right)$ as a ``reduced'' polynomial in $\sig$ by grouping all the coefficients of terms $\sig^\be$ for each multi-index $\be \in \Z_{\geq 0}^{n-1}$. To stimulate intuition, we remark that in the case of dimension $n=2$ the relations (\ref{z_dfns}) indicate the relations
\begin{align*}
u_1+ z_1 = \frac{u_1}{|u|} \left(|u|+\tau\right) + \frac{u_2}{|u|}  \sigma, \\
u_2+ z_2 = \frac{u_2}{|u|} \left(|u|+ \tau\right) - \frac{u_1}{|u|}  \sigma.
\end{align*}
Hence
\begin{eqnarray*}
Q_{\vec{\nu}}\left(u+ z\right)
&=& \sum_{\alpha} \vec{\nu}_{\alpha} \left(u_1+ z_1\right)^{\alpha_1} \left(u_2+ z_2\right)^{\alpha_2} \\
&=& \sum_{\alpha} \sum_{j=0}^{\alpha_1} \sum_{k=0}^{\alpha_2} \vec{\nu}_{\alpha} \binom{\alpha_1}{j} \binom{\alpha_2}{k} (-1)^k  \left(\frac{u_1}{|u|}\right)^{\alpha_1-j+k} \left(\frac{u_2}{|u|}\right)^{\alpha_2-k+j}\left(|u|+\tau\right)^{|\alpha|-j-k} \sigma^{j+k}.
\end{eqnarray*}
Similarly
\[Q_{\vec{\mu}}(z)
= \sum_{\alpha} \sum_{j=0}^{\alpha_1} \sum_{k=0}^{\alpha_2} \vec{\mu}_{\alpha}  \binom{\alpha_1}{j} \binom{\alpha_2}{k} (-1)^k \left(\frac{u_1}{|u|}\right)^{\alpha_1-j+k} \left(\frac{u_2}{|u|}\right)^{\alpha_2-k+j} \tau^{|\alpha|-j-k} \sigma^{j+k}.\]
Thus by inspection, the coefficient of $\sigma^l$ in the phase $Q_{\nu}\left(u+ z\right) - Q_{\mu}(z)$ is
\sum_{|\alpha| \geq l} A_{\alpha,l}\left(\frac{u}{|u|}\right) \left(\vec{\nu}_{\alpha} \left(|u|+ \tau\right)^{|\alpha|-l} - \vec{\mu}_{\alpha} \tau^{|\alpha|-l} \right)
\] 
where the function
\beq\label{A_dfn_3}
A_{\alpha,l}(w) := \sum_{\substack{0 \leq j \leq \alpha_1 \\ 0 \leq k \leq \alpha_2 \\ j+k = l}} (-1)^k \binom{\alpha_1}{j} \binom{\alpha_2}{k} w_1^{\alpha_1-j+k} w_2^{\alpha_2-k+j} 
\eeq
is a homogeneous polynomial in $w\in \R^2$ of degree $|\alpha|$. 
}
In arbitrary dimensions, it is too cumbersome to perform an explicit binomial expansion, so we use the implicit definition (\ref{defining_A_rel}) instead.

We note that one may also generalize the relation (\ref{defining_A_rel}) to
\beq\label{defining_A_rel_T}
\sum_{|\be| \leq |\al| }A_{\al,\be}(w) T^{|\al| - |\be|} \sig^\be = (w_n T+ w^{(n)} \cdot \sig)^{\al_n} (w^{(n)}T  - w_n \sig)^{\al^{(n)}},
\eeq
for a generic variable $T \in \R$.
To proceed with the proof of Lemma \ref{lemma_P_phase_n}, we use the relations (\ref{z_rel1}) and (\ref{z_rel2}) to compute that 
\begin{align*}
u_n+ z_n = \frac{u_n}{|u|} \left(|u|+ \tau\right) + \frac{u^{(n)}\cdot \sigma }{|u|}, \\
u^{(n)}+ z^{(n)} = \frac{u^{(n)}}{|u|} \left(|u|+ \tau\right) - \frac{u_n}{|u|}\sigma.
\end{align*}
We then write the expansion:
\[P_\nu \left(u+ z\right)
= \sum_{m=2}^d \nu_m \sum_{|\al| =m} c_\al \left(\frac{u_n}{|u|} \left(|u|+ \tau\right) + \frac{u^{(n)}}{|u|} \cdot \sigma\right)^{\alpha_n} \left( \frac{u^{(n)}}{|u|} \left(|u|+ \tau\right) - \frac{u_n}{|u|} \sigma\right)^{\alpha^{(n)}} .\]
Applying (\ref{defining_A_rel_T}) with $w = u/|u|$ and $T = (|u| + \tau)$ we now see that 
 \begin{eqnarray*}
P_\nu \left(u+ z\right) &=& \sum_{m=2}^d \nu_m \sum_{|\alpha|=m}c_\al \left( \sum_{|\ga| \leq |\al|} A_{\al,\ga} \left( \frac{u}{|u|} \right) (|u| + \tau)^{|\al| - |\ga|} \sig^\ga \right) \\
 & = & \sum_{0 \leq |\be| \leq d} \left( \sum_{m \geq |\be|} \nu_m \sum_{|\al| =m} c_\al A_{\al,\be} \left( \frac{u}{|u|} \right)  \left(|u|+ \tau\right)^{|\al| - |\be|} \right) \sig^\be.
 \end{eqnarray*}

We compute similarly that $P_\mu(z)$ can be written as 
\begin{eqnarray}
P_\mu(z) 
 &=& \sum_{m=2}^d \mu_m \sum_{|\al|=m} c_\al \left(\frac{u_n}{|u|} \tau + \frac{u^{(n)}}{|u|}  \cdot \sigma\right)^{\alpha_n} \left( \frac{u^{(n)}}{|u|} \tau - \frac{u_n}{|u|} \sigma\right)^{\alpha^{(n)}}  \nonumber \\
	 &=&  \sum_{0 \leq |\be| \leq d} \left( \sum_{m \geq |\be|} \mu_m \sum_{|\al| =m} c_\al A_{\al,\be} \left( \frac{u}{|u|} \right) \tau^{|\al| - |\be|} \right) \sig^\be\label{Q_A_dfn_n}
	 \end{eqnarray}
with the same functions $A_{\al,\be}$, and this completes the proof of Lemma \ref{lemma_P_phase_n}.

\subsection{Preliminary properties of $A_{\al,\be}$ and $B_{m,\be}$}
We summarize the key properties of the polynomials $A_{\al,\be}$, which we will prove in Section \ref{sec_A_B_lemmas}.
\begin{lemma}[Properties of $A_{\al,\be}$]\label{lemma_A_ab}
$\quad$
\begin{enumerate}
\item \label{item_A_1} Suppose that $|\be|=1$, so that for some $1 \leq j \leq n-1$ we may write $\be = e_j$ where $e_j = (0,\ldots, 1, \ldots, 0)$ is the $j$-th unit vector. Then for any $|\al| \geq 1$,
\[ A_{\al,e_j}(w) =\al_n w^{\al + e_j  - e_n} - \al_j w^{\al-e_j +e_n}. \]
\item \label{item_A_2} 
If $\be = 2e_j$ with $1 \leq j \leq n-1$, then for any $|\al| \geq 2$,
\[ A_{\al, 2e_j}(w) = {\al_n \choose 2} w^{\al + 2e_j - 2e_n} - \al_n \al_j w^\al + {\al_j \choose 2} w^{\al -2e_j + 2e_n}.\]
\xtra{Suppose that $|\be|=2$. If $\be = e_i + e_j$ with $i \neq j$, $1 \leq i,j \leq n-1$, then 
\[ A_{\al,e_i+e_j} (w) = {\al_n \choose 2} w^{\al+e_i+e_j - 2e_n} - \al_n \al_{i} w^{\al - e_i + e_j} - \al_n \al_j w^{\al +e_i - e_j} + \al_i \al_j w^{\al -e_i - e_j + 2e_n}.}
\item \label{item_A_be_al}
If $|\be|=|\al| $, then $A_{\al,\be}(\om)$ is a monomial, 
\[ A_{\al,\be}(w) = C_{\al,\be} (-1)^{|\al^{(n)}|} w_n^{|\al^{(n)}|}(w^{(n)})^{\be - \al^{(n)}},\]
for a non-negative combinatorial constant $C_{\al,\be}$. 
Moreover,  given any multi-index $\al \in \Z_{\geq 0}^n$ with $|\al| \geq 1$, there exists a multi-index $\be \in \Z_{\geq 0}^{n-1}$ with $|\be| =|\al|$ such that the coefficient $C_{\al,\be}$ is nonzero.
\end{enumerate} 
\end{lemma}

We also record the key properties of the polynomials $B_{m,\be}(w)$ which we require; these are proved in Sections \ref{sec_B_lemma_1} to \ref{sec_B_lemma_3}. While these appear to be simple properties of a combinatorial nature, they encode the advantages of the restricted class of polynomials we consider, and lie at the heart of the bound for $K_\sharp^{\nu,\mu}$.
As before, we let
\[
\Lambda = \{ 2 \leq m \leq d \colon p_m(y) \not\equiv 0 \}.
\]
We recall that for $m \in \Lambda$,  $p_m(y)$ is parabolic if  $p_m(y) = C|y|^m$ for some constant $C \ne 0$, and non-parabolic otherwise. 
\begin{lemma}[Properties of $B_{m,\be}$]\label{lemma_B_n}
$\quad$
\begin{enumerate}
\item \label{item_dependence_n} The polynomial $B_{m,\be}$ depends only on $p_m$ and $\be$ (and not on any other $p_j$, $j \ne m$). In particular, if $p_m \equiv 0$, then $B_{m,\be} \equiv 0$ for all $\be$. 
\item \label{item_bound_n} There is a constant $C_B$ depending only on 
the degree $d$, the dimension $n$, and the fixed coefficients $c_\al$ such that for all $1 \leq m \leq d$ and all $1 \leq |\be| \leq d$,
\[  |B_{m,\be}(w)| \leq C_B, \quad \text{for all $|w| \leq 1$}.\] 
\item $B_{m,\be}(w)$ is a homogeneous polynomial in $w$ of degree $m$. \label{item_hom_n}
\item \label{item_zero_n} If $m \in \Lambda$, then there exists some $|\be|=m$ such that $B_{m,\be}$ is not the zero polynomial.
\item If $m \in \Lambda$ and $p_m(y)$ is parabolic, then $B_{m,\be}(w)$ is the zero polynomial for all $\be$ with $|\be|$ odd. \label{item_parabolic_n}
\item If $m \in \Lambda$ and $p_m(y)$ is not parabolic, then there exists $\be$ with $|\be|=1$ such that $B_{m,\be}(w)$ is not the zero polynomial. \label{item_m1_n}
\item If $m \in \Lambda$ and $p_m(y)$ is parabolic, then there exists some $\be$ with $|\be|=2$ such that $B_{m,\be}(w)$ is not the zero polynomial. \label{item_m2_n}
\end{enumerate}
\end{lemma}
We now demonstrate how to derive the $K^{\nu,\mu}_{\sharp,n}$ bound of Proposition~\ref{prop_Ksharp_n} from these properties. As in the specific examples we considered in Section \ref{sec_vdC_n2},  the general strategy of the proof will be to understand when $|C[\sig^\be](u,\tau)|$ is small, since if $|C[\sig^\be](u,\tau)|$ is large for some $1 \leq |\be| \leq d$, then we may bound the kernel $K_{\sharp,n}^{\nu,\mu}$ by a van der Corput estimate.

For each $1 \leq |\be | \leq d$, 
We re-write the expression (\ref{sig_coeff_n_res}) defining $C[\sig^\be](u,\tau)$ for $|\be| =l$ as 
\begin{multline}\label{PC_dfn'_n}
C[\sig^\be](u,\tau) = \sum_{m=l}^d  (\nu_m - \mu_m) \tau^{m-l} B_{m,\be}\left(\frac{u}{|u|}\right)  \\
 + \; \sum_{m=l+1}^d \nu_m( (|u| + \tau)^{m-l} - \tau^{m-l})B_{m,\be} \left(\frac{u}{|u|}\right) .
\end{multline} 
We first use a downward induction process to eliminate the presence of $\mu = (\mu_2,\ldots, \mu_d)$ in the coefficient $C[\sig^\be](u,\tau)$ that we hope to show is large. 
The expression (\ref{PC_dfn'_n}) makes clear that for each $1 \leq |\be| \leq d$, $C[\sig^\be](u,\tau)$ depends only on $\mu_m$ with $m \geq |\be|$. Our strategy thus relies on the fact that if $|C[\sig^\be](u,\tau)|$ is small for all $|\be| \geq l_0$, then for all $l \geq l_0$ one can essentially rewrite $\mu_l$ in terms of $\nu = (\nu_2,\ldots, \nu_d)$, $u$ and $\tau$. Note also that $C[\sig^\be](u,\tau)$ is a polynomial of degree at most $d-|\be|$ in $\tau$. 

We now make precise the process of inductively eliminating the presence of $\mu$-coefficients.
We recall that we have by assumption $r \leq \| \nu \| , \| \mu \| \leq 2r$. We will consider three cases, motivated by the key examples we considered in Section \ref{sec_vdC_n2}.
Case A: a non-parabolic term dominates, namely there exists $m \in \Lambda$, with $p_m(y)$ not parabolic, such that $|\nu_m| \simeq r$; Case B: a parabolic term dominates, namely there exists $m \in \Lambda$ with $p_m(y)$ parabolic, such that $|\nu_m| \simeq r$; we then further subdivide this into case B1 when $p_2(y) \con 0$ and case B2 when $p_2(y) \not\con 0$.

\subsection{Case A: a non-parabolic term dominates}
In this case we aim to show that there exists a multi-index $\be = \be^*$ with $|\be^*|=1$ such that $C[\sig^{\be^*}](u,\tau)$ is large for ``most'' $u$ and $\tau$.

Fix $0< \ep_1< \ep_2 < 1$. 
If $|C[\sig^\be](u,\tau)| \geq r^{\ep_1}$ for some $2 \leq |\be| \leq d$, then $|K_{\sharp,n}^{\nu,\mu}(u,\tau)| \leq C r^{-\ep_1/d}$ by the van der Corput estimate of Lemma \ref{lemma_Prop2.1}, as desired. So we may assume that
\beq\label{C_assumed_small}
|C[\sig^\be](u,\tau)| \leq r^{\ep_1} \quad \text{for all $2 \leq |\be| \leq d$}.
\eeq
Applying the assumption (\ref{C_assumed_small}) to (\ref{PC_dfn'_n}) shows that for each $2 \leq l \leq d$, for every $|\be|=l$ we have
\begin{multline}\label{numu_B_n}
(\nu_l-\mu_l) B_{l,\be}\left(\frac{u}{|u|}\right) = - \sum_{m=l+1}^d (\nu_m-\mu_m) \tau^{m-l} B_{m,\be}\left(\frac{u}{|u|}\right) \\
	- \sum_{m=l+1}^d \nu_m ((|u|+\tau)^{m-l}-\tau^{m-l}) B_{m,\be}\left(\frac{u}{|u|}\right) + O(r^{\ep_1}).
\end{multline}
We now use Statement \ref{item_zero_n} of Lemma \ref{lemma_B_n} to conclude that for each $2 \leq l \leq d$ with $l \in \Lambda$, there exists some $\be$ with $|\be|=l$ for which $B_{l,\be}$ is not the zero polynomial. For each $l \in \Lambda$ we will pick such a distinguished index $\be$ and denote it by $\be(l)$. We thus obtain for each $l \in \Lambda$ a relation (\ref{numu_B_n}) specialized to the distinguished index $\be(l)$ that will allow us to express $\nu_l - \mu_l$  in terms of $\nu_m-\mu_m$ for $m >l$, and some harmless terms involving only $u, \tau$ and $\nu = (\nu_2,\ldots, \nu_d)$. 

By a downward induction argument on $l$, we conclude that for all $2 \leq l \leq d$ we may write
\beq\label{S_nu_u_1_n}
(\nu_l - \mu_l) \prod_{\substack{j \geq l \\ j \in \Lambda}} B_{j,\be(j)}\left(\frac{u}{|u|}\right) = S_{\nu,u}^{(l)}(\tau) + O(r^{\ep_1}),
\eeq
 where $S_{\nu,u}^{(l)}(\tau)$ is a polynomial in $\tau$ whose coefficients depend only on $\nu$ and $u$ (but not $\mu$). 

We now feed this information into an analysis of $C[\sig^\be](u,\tau)$ for indices $|\be|=1$. We fix any index $\be$ with $|\be|=1$ and see that in (\ref{PC_dfn'_n}), $C[\sig^\be](u,\tau)$ may be expressed as
\begin{multline}\label{C_sig_case1_n}
C[\sig^\be](u,\tau) = \sum_{m=2}^d  (\nu_m - \mu_m) \tau^{m-1} B_{m,\be}\left(\frac{u}{|u|}\right)  \\
+\; 
\sum_{m=2}^d \nu_m( (|u| + \tau)^{m-1} - \tau^{m-1})B_{m,\be} \left(\frac{u}{|u|}\right).
\end{multline}
Note in particular that the first sum begins with $m = 2$,  since $p_1(y) \equiv 0$ and so we have $B_{1,\be} \equiv 0$ for all $\be$; thus all terms in the first sum are linear or higher order with respect to $\tau$.  We will use (\ref{S_nu_u_1_n}) to eliminate the presence of $\mu$ in the first sum in (\ref{C_sig_case1_n}), and then we will single out the constant terms with respect to $\tau$ in the resulting polynomial (which will come only from the second sum on the right hand side of (\ref{C_sig_case1_n})). To proceed with this plan, we multiply (\ref{C_sig_case1_n}) through by the  polynomial
\beq\label{prod_poly}
 \prod_{\substack{j \geq 2 \\ j \in \Lambda}} B_{j,\be(j)}\left(\frac{u}{|u|}\right),
 \eeq
 and substitute (\ref{S_nu_u_1_n}) wherever possible. (Recall that by construction, (\ref{prod_poly}) is not identically zero.) One then concludes that
\beq\label{C_R_exprA}
C[\sig^\be](u,\tau) \prod_{\substack{j \geq 2 \\ j \in \Lambda}} B_{j,\be(j)}\left(\frac{u}{|u|}\right) 
= R_{\nu,u}^{(\be)}(\tau) + O(r^{\ep_1}), 
\eeq
where $R_{\nu,u}^{(\be)}(\tau)$ is a polynomial in $\tau$ of degree $\leq d-1$ whose coefficients depend only on $\nu$ and $u$ (but not $\mu$). (The superscript $\be$ reflects that this comes from the coefficient $C[\sig^\be](u,\tau)$.) In fact, we will write $R_{\nu,u}^{(\be)}$ as a constant term in $\tau$, plus higher powers of $\tau$. The constant term in $\tau$ arise from the contribution from the second term of (\ref{C_sig_case1_n}) only, so
\[
R_{\nu,u}^{(\be)}(\tau) = 
\sum_{m=2}^d \nu_m |u|^{m-1} B_{m,\be} \left(\frac{u}{|u|}\right)  \prod_{\substack{j \geq 2 \\ j \in \Lambda}} B_{j,\be(j)}\left(\frac{u}{|u|}\right) + R^{(\be,1)}_{\nu,u}(\tau),
\] 
where  the first term is constant with respect to $\tau$, and $R^{(\be,1)}_{\nu,u}(\tau)$ is a polynomial in $\tau$ with no constant term and with coefficients that depend only on $u,\nu$ (and not $\mu$).
Let $W_{\nu}^{(\be)}(u)$ be the polynomial in $u$ defined by
\beq\label{W_nu_def_case1_n}
W_{\nu}^{(\be)}(u) := 
\sum_{m=2}^d \nu_m B_{m,\be}(u) \prod_{\substack{j \geq 2 \\ j \in \Lambda}} B_{j,\be(j)}(u).
\eeq
We note that since $B_{m,\be}$ is homogeneous of degree $m$, the degree of $W_\nu^{(\be)}(u)$ with respect to $u$ is at most 
\beq\label{s0_dfn1_n}
 s_0 = d+ \sum_{\bstack{j \geq 2}{j \in \Lambda}} j.
 \eeq
\xtra{Furthermore, we recall that for each $j \in \Lambda$ we chose an index $\be(j)$ such that $B_{j,\be(j)}(u)$ is not identically zero, by Statement \ref{item_zero_n} of Lemma~\ref{lemma_B_n}.} 
Moreover, we see that
\[
R_{\nu,u}^{(\be)}(\tau) = |u|^{-s_1} W_{\nu}^{(\be)}(u) +R^{(\be,1)}_{\nu,u}(\tau),
\]
where 
\[ s_1 = 1 + \sum_{\bstack{j\geq 2}{ j \in \Lambda}} j.\]

Recall that in the current case, we assumed that there exists $m \in \Lambda$, say $m_0$, with $p_{m_0}(y)$ not parabolic, such that $|\nu_{m_0}| \simeq r$. For that $m_0$, we know that there exists an index $\be$ with $|\be|=1$ such that $B_{m_0,\be}(u)$ is not identically zero by Statement \ref{item_m1_n} of Lemma~\ref{lemma_B_n}; we will denote this distinguished $\be$ by $\be^*$; this choice of a distinguished index $\be^*$ depends only on the original choice of polynomial $p_{m_0}$.
 We now define our exceptional sets,  with respect to the fixed $\be^*$ with $|\be^*|=1$ determined above. Recall that $0<\ep_1<\ep_2<1$ and set
\[
G^{\nu} := \{ u \in B_2 (\R^n)\colon |W_{\nu}^{(\be^*)}(u)| \leq r^{\ep_2} \},
\]
and for $u \in B_2 \setminus G^{\nu}$, let 
\[
F^{\nu}_u := \{ \tau \in B_1(\R) \colon |R_{\nu,u}^{(\be^*)}(\tau)| \leq C_0 r^{\ep_1} \}
\]
for some large absolute constant $C_0$ to be determined. Also define $F^{\nu}_u := \emptyset$ if $u \in G^{\nu}$. Then for all $(u,\tau) \in B_2 (\R^n)\times B_1 (\R)$ with $u \notin G^{\nu}$, $\tau \notin F^{\nu}_u$, we have
\[
|R_{\nu,u}^{(\be^*)}(\tau)| \geq C_0 r^{\ep_1},
\]
so that for such $(u,\tau)$, we conclude from (\ref{C_R_exprA}) that
\begin{eqnarray*}
|C[\sig^{\be^*}](u,\tau)| 
&\geq& C_B^{-|\Lambda|}\left|C[\sig^{\be^*}](u,\tau) \prod_{\substack{j \geq 2 \\ j \in \Lambda}} B_{j,\be(j)}\left(\frac{u}{|u|}\right)\right|  \\
&\geq &C_B^{-|\Lambda|} \left( |R_{\nu,u}^{(\be^*)}(\tau)| - O(r^{\ep_1})  \right)\\
&\geq &C_B^{-|\Lambda|} \left( C_0 r^{\ep_1} - O(r^{\ep_1}) \right)\\
&\geq &C_B^{-|\Lambda|} r^{\ep_1} 
\end{eqnarray*}
if $C_0$ is sufficiently large. It follows from the van der Corput estimate of Lemma \ref{lemma_Prop2.1} that
\[
|K_{\sharp,n}^{\nu,\mu}(u,\tau)| \leq C r^{-\ep_1/d} \quad \text{if $u \notin G^{\nu}$ and $\tau \notin F_u^{\nu}$},
\]
for some fixed constant $C$ dependent only on the initial choice of the polynomials $p_m$.

On the other hand, $K_{\sharp,n}^{\nu,\mu}(u,\tau)$ is supported on $B_2(\R^n) \times B_1(\R)$, and is bounded by some uniform constant $C$ (see (\ref{K_triv_bd_n})). Hence to complete the proof of Proposition~\ref{prop_Ksharp_n} in this case, we only need to show that $G^{\nu}$ and $F_u^{\nu}$ are sets of small measures. Now $G^{\nu}$ is the set of $u \in B_2$ where $W_{\nu}^{(\be^*)}(u)$ is small, and in (\ref{W_nu_def_case1_n}), $W_{\nu}^{(\be^*)}(u)$ is represented as a sum of homogeneous polynomials of different total degrees. In particular, the coefficient $\nu_{m_0}$ for which $|\nu_{m_0}| \simeq r$ appears in $W_{\nu}^{(\be^*)}(u)$ as the coefficient of a homogeneous polynomial that is by construction not identically zero, and which has different total degree from all other terms in $W_{\nu}^{(\be^*)}(u)$. Thus we see that
\[
\llbracket W_{\nu}^{(\be^*)} \rrbracket_u \geq C |\nu_{m_0}| \geq C r,
\]
which implies by Lemma \ref{lemma_Prop2.2_0} that
\[
|G^{\nu}| \leq C \left( \frac{r^{\ep_2}}{r} \right)^{1/s_0} = C r^{-(1-\ep_2)/s_0},
\]
since the degree of $W_{\nu}^{(\be^*)}(u)$ is at most  $s_0$. Furthermore, if $u \in B_2 \setminus G^{\nu}$, then 
\[
\llbracket R_{\nu,u}^{(\be^*)} \rrbracket_{\tau} 
\geq |u|^{-s_1} |W_{\nu}^{(\be^*)}(u)|
\geq C |W_{\nu}^{(\be^*)}(u)|
\geq C r^{\ep_2}.
\] 
Thus if $u \in B_2 \setminus G^{\nu}$, then by Lemma \ref{lemma_Prop2.2_0},
\[
|F_u^{\nu}| \leq C \left( \frac{C_0 r^{\ep_1}}{r^{\ep_2}} \right)^{1/(d-1)} = C r^{-(\ep_2-\ep_1)/(d-1)}.
\] 
The same inequality is clearly true if $u \in G^{\nu}$, since then $F_u^\nu = \emptyset$. Thus we have $|F_u^{\nu}|$ being small for all $u \in B_2$. This concludes the proof of Proposition~\ref{prop_Ksharp_n} in the case under consideration.

\subsection{Case B1: a parabolic term dominates and $p_2(y) \equiv 0$}
In Case B there exists $m \in \Lambda$ with $p_m(y)$ parabolic such that $|\nu_m| \simeq r$.
We recall the division into two subcases: Case B1, in which $p_2(y) \equiv 0$, and Case B2, in which $p_2(y) \not\equiv 0$.

We first consider Case B1. Similarly to Case A, we fix $0<\ep_1< \ep_2 < 1$ and start by assuming  without loss of generality that $|C[\sig^\be](u,\tau)| \leq r^{\ep_1}$ for all $3 \leq |\be| \leq d$. For each $3 \leq l \leq d$ we choose (via Statement \ref{item_zero_n} of Lemma \ref{lemma_B_n}) a distinguished index $\be(l)$ with $|\be(l)|=l$ such that $B_{l,\be(l)}$ is not the zero polynomial. For each pair $(l,\be(l))$ we use the relation (\ref{numu_B_n}) to provide an expression for $(\nu_l - \mu_l)B_{l,\be(l)}$ which we then feed into a downward induction argument, with the result that for all $3 \leq l \leq d$,
\beq\label{S_nu_u_2a_n}
(\nu_l - \mu_l) \prod_{\substack{j \geq l \\ j \in \Lambda}} B_{j,\be(j)}\left(\frac{u}{|u|}\right) = S_{\nu,u}^{(l)}(\tau) + O(r^{\ep_1}),
\eeq
where $S_{\nu,u}^{(l)}(\tau)$ is a polynomial in $\tau$ whose coefficients depend only on $\nu$ and $u$ (but not $\mu$). We now feed this information into an analysis of $C[\sig^\be](u,\tau)$ for some $|\be|=2$ to be chosen precisely later. 

For now, fix any $\be$ with $|\be|=2$. Since $p_2(y) \equiv 0$, we have $B_{2,\be} \equiv 0$ for all $\be$, so by (\ref{PC_dfn'_n}), $C[\sig^\be](u,\tau)$ reduces to
\begin{multline}\label{C_sig2_case2a_n}
C[\sig^\be](u,\tau) = \sum_{m=3}^d  (\nu_m - \mu_m) \tau^{m-2} B_{m,\be}\left(\frac{u}{|u|}\right)  \\
+\; 
\sum_{m=3}^d \nu_m( (|u| + \tau)^{m-2} - \tau^{m-2})B_{m,\be} \left(\frac{u}{|u|}\right).
\end{multline}
In particular, since the first sum starts with $m \geq 3$, all terms from the first sum are linear or higher order with respect to $\tau$. The idea is to use (\ref{S_nu_u_2a_n}) to eliminate the role of $\mu$ in the first sum, and then to consider the constant term with respect to $\tau$ coming from the second sum. To proceed with this, we multiply (\ref{C_sig2_case2a_n}) through by $\prod_{\substack{j \geq 3 \\ j \in \Lambda}} B_{j,\be(j)}\left(\frac{u}{|u|}\right)$, and substitute (\ref{S_nu_u_2a_n}) wherever possible. One then concludes that
\beq\label{C_exprB1}
C[\sig^\be](u,\tau) \prod_{\substack{j \geq 3 \\ j \in \Lambda}} B_{j,\be(j)}\left(\frac{u}{|u|}\right) 
= R_{\nu,u}^{(\be)}(\tau) + O(r^{\ep_1}), 
\eeq
where $R_{\nu,u}^{(\be)}(\tau)$ is a polynomial in $\tau$ whose coefficients depend only on $\nu$ and $u$ (but not $\mu$). (The superscript $\be$ again reflects that this comes from the coefficient $C[\sig^\be](u,\tau)$.) In fact, 
\[
R_{\nu,u}^{(\be)}(\tau) = 
\sum_{m=3}^d \nu_m |u|^{m-2} B_{m,\be} \left(\frac{u}{|u|}\right) \prod_{\substack{j \geq 3 \\ j \in \Lambda}} B_{j,\be(j)}\left(\frac{u}{|u|}\right) + R_{\nu,u}^{(\be,1)}(\tau),
\] 
where the first term is constant with respect to $\tau$ and $R_{\nu,u}^{(\be,1)}(\tau)$ is a polynomial in $\tau$ with no constant term and with coefficients dependent only on $\nu,u$ (and not $\mu$).
Let $W_{\nu}^{(\be)}(u)$ be the polynomial in $u$ defined by
\beq\label{W_nu_def_case2a_n}
W_{\nu}^{(\be)}(u) := 
\sum_{m=3}^d \nu_m B_{m,\be}(u)  \prod_{\substack{j \geq 3 \\ j \in \Lambda}} B_{j,\be(j)}(u),
\eeq
which we note has degree at most 
\[s_2 = d + \sum_{\bstack{j \geq 3}{j \in \Lambda}} j.\]
We also recall that each $B_{j,\be(j)}(u)$ is not identically zero whenever $j \in \Lambda$, because we have chosen $\be(j)$ using Statement \ref{item_zero_n} of Lemma~\ref{lemma_B_n}.
Then we can rewrite
\[
R_{\nu,u}^{(\be)}(\tau) = |u|^{-s_3} W_{\nu}^{(\be)}(u) + R_{\nu,u}^{(\be,1)}(\tau),
\]
where 
\[ s_3 = 2 + \sum_{\bstack{j \geq3}{ j \in \Lambda}} j.
\]

Recall that in the current case, we assumed that there exists $m \in \Lambda$, say $m_0$, with $p_{m_0}(y)$ parabolic such that $|\nu_{m_0}| \simeq r$. Since $p_2(y) \equiv 0$, we have $2 \notin \Lambda$, so the $m_0$ above cannot be 2; furthermore, by Statement \ref{item_m2_n} of Lemma~\ref{lemma_B_n}, for this $m_0$, there exists a $\be$ with $|\be|=2$ such that $B_{m,\beta}(u)$ is not identically zero. We will denote this distinguished index by $\be^*$; the choice of $\be^*$ depends only on $p_{m_0}$.
We now define our exceptional sets as follows: recall that $0<\ep_1<\ep_2<1$ and set
\[
G^{\nu} := \{ u \in B_2(\R^n) \colon |W_{\nu}^{(\be^*)}(u)| \leq r^{\ep_2} \},
\]
and for $u \in B_2 \setminus G^{\nu}$, let 
\[
F^{\nu}_u := \{ \tau \in B_1(\R) \colon |R_{\nu,u}^{(\be^*)}(\tau)| \leq C_0 r^{\ep_1} \}
\]
for some large absolute constant $C_0$ to be determined. Also define $F^{\nu}_u := \emptyset$ if $u \in G^{\nu}$. Then for all $(u,\tau) \in B_2(\R^n) \times B_1(\R)$ with $u \notin G^{\nu}$, $\tau \notin F^{\nu}_u$, we have
\[
|R_{\nu,u}^{(\be^*)}(\tau)| \geq C_0 r^{\ep_1},
\]
so for such $(u,\tau)$, we conclude from (\ref{C_exprB1}) that
\begin{eqnarray*}
|C[\sig^{\be^*}](u,\tau)| 
&\geq & C_B^{-|\Lambda|}\left|C[\sig^{\be^*}](u,\tau) \prod_{\substack{j \geq 3 \\ j \in \Lambda}} B_{j,\be(j)}\left(\frac{u}{|u|}\right)\right|  \\
&\geq & C_B^{-|\Lambda|} \left( |R_{\nu,u}^{(\be^*)}(\tau)| - O(r^{\ep_1}) \right)\\
&\geq& C_B^{-|\Lambda|}\left( C_0 r^{\ep_1} - O(r^{\ep_1}) \right)\\
&\geq& C_B^{-|\Lambda|}r^{\ep_1} 
\end{eqnarray*}
if $C_0$ is sufficiently large. It follows from the van der Corput estimate of Lemma \ref{lemma_Prop2.1} that
\[
|K_{\sharp,n}^{\nu,\mu}(u,\tau)| \leq C r^{-\ep_1/d} \quad \text{if $u \notin G^{\nu}$ and $\tau \notin F_u^{\nu}$}.
\]
On the other hand, we will now show that $G^{\nu}$ and $F_u^{\nu}$ are sets of small measures. First, $G^{\nu}$ is the set of $u \in B_2$ where $W_{\nu}^{(\be^*)}(u)$ is small, and in (\ref{W_nu_def_case2a_n}), $W_{\nu}^{(\be^*)}(u)$ is represented as a sum of  homogeneous polynomials of different total degrees; in particular $\nu_{m_0}$ appears as the coefficient of a homogeneous polynomial in $W_{\nu}^{(\be^*)}(u)$ that is by construction not identically zero, and has different total degree than all other terms in $W_{\nu}^{(\be^*)}(u)$.  Thus in the current case,
\[
\llbracket W_{\nu}^{(\be^*)} \rrbracket_u \geq C |\nu_{m_0}| \geq C r,
\]
which implies by Lemma \ref{lemma_Prop2.2_0} that
\[
|G^{\nu}| \leq C \left( \frac{r^{\ep_2}}{r} \right)^{1/s_2} = C r^{-(1-\ep_2)/s_2},
\]
where  the degree of $W^{(\be^*)}$ is at most $s_2$. Furthermore, if $u \in B_2 \setminus G^{\nu}$, then 
\[
\llbracket R_{\nu,u}^{(\be^*)} \rrbracket_{\tau} 
\geq |u|^{-s_3} |W_{\nu}^{(\be^*)}(u)|
\geq C |W_{\nu}^{(\be^*)}(u)|
\geq C r^{\ep_2}.
\] 
Thus if $u \in B_2 \setminus G^{\nu}$, then by Lemma \ref{lemma_Prop2.2_0},
\[
|F_u^{\nu}| \leq C \left( \frac{C_0 r^{\ep_1}}{r^{\ep_2}} \right)^{1/d} = C r^{-(\ep_2-\ep_1)/d}.
\] 
The same inequality is clearly true if $u \in G^{\nu}$. Thus we have $|F_u^{\nu}|$ being small for all $u \in B_2$. This concludes the proof of Proposition~\ref{prop_Ksharp_n} in the case under consideration.

 \subsection{Case B2: a parabolic term dominates and $p_2(y) \not\equiv 0$}

In this case, we will need to consider the coefficients of two terms $\sig^{\be_1^*}$ and $\sig^{\be_2^*}$ for some $|\be_1^*|=1$ and $|\be_2^*|=2$ to be chosen precisely later. We fix $0< \ep_1< \ep_2 <1$ and 
suppose first that $|C[\sig^\be](u,\tau)| \geq r^{\ep_1}$ for some $1 \leq |\be| \leq d$ with $|\be| \ne 2$. Then $|K_{\sharp,n}^{\nu,\mu}(u,\tau)| \leq C r^{-\ep_1/d}$ as desired. Thus we may assume that $|C[\sig^\be](u,\tau)| \leq r^{\ep_1}$ for all $\be$ with $|\be|=1$ and all $\be$ with $3 \leq |\be| \leq d$. From the latter set of conditions, we use the relations (\ref{numu_B_n}) for each $3 \leq l \leq d$ and an appropriate choice $\be(l)$ such that $B_{l,\be(l)}$ is not the zero polynomial, to provide expressions for $(\nu_l - \mu_l) B_{l,\be(l)}$ in terms of $\nu_m - \mu_m$ for $m >l$. We then feed these expressions into a downward induction in $l$ in order to show that for all $3 \leq l \leq d$,
\beq\label{S_nu_u_2b_1_n}
(\nu_l - \mu_l) \prod_{\substack{j \geq l \\ j \in \Lambda}} B_{j,\be(j)}\left(\frac{u}{|u|}\right) = S_{\nu,u}^{(l)}(\tau) + O(r^{\ep_1}),
\eeq
where $S_{\nu,u}^{(l)}(\tau)$ is a polynomial in $\tau$ whose coefficients depend only on $\nu$ and $u$ (but not $\mu$). 

Next for each $|\be| = 1$, we apply the assumption that $|C[\sig^\be](u,\tau)|\leq r^{\ep_1}$, to solve for $(\nu_2 - \mu_2)\tau B_{2,\be}$ in (\ref{PC_dfn'_n}). We conclude that for each $\be$ with $|\be|=1$ we have
\begin{multline}\label{Be21_n}
(\nu_2 - \mu_2) \tau B_{2,\be}\left( \frac{u}{|u|} \right) = -\sum_{m=3}^d (\nu_m-\mu_m) \tau^{m-1} B_{m,\be} \left( \frac{u}{|u|} \right) \\
 - \sum_{m=2}^d \nu_m ((|u|+\tau)^{m-1}-\tau^{m-1}) B_{m,\be}\left( \frac{u}{|u|} \right) + O(r^{\ep_1}).
\end{multline}
Now we recall that by Statement \ref{item_m1_n} of Lemma \ref{lemma_B_n}, if $m \in \Lambda$ and $p_m(y)$ is not parabolic, then there exists an index $\be$ with $|\be|=1$ such that $B_{m,\be}$ is not the zero polynomial. By the assumption of our main theorem, we always have $p_2$ not parabolic, and since in the case under consideration $2 \in \Lambda$,  there exists by Statement \ref{item_m1_n} of Lemma \ref{lemma_B_n} an index $|\be|=1$, which we will denote by $\be_1^*$, for which $B_{2,\be_1^*}$ is not the zero polynomial. We will use the relation (\ref{Be21_n}) specialized to the choice $\be = \be_1^*$. 

We multiply through (\ref{Be21_n}), specialized to $\be = \be_1^*$, by the polynomial 
\[\prod_{\substack{j \geq 3 \\ j \in \Lambda}} B_{j,\be(j)}\left(\frac{u}{|u|}\right),\]
 use (\ref{S_nu_u_2b_1_n}) to eliminate the presence of $\mu$ in the first sum on the right hand side, and then group terms according to powers of $\tau$. We note that the second term on the right hand side of (\ref{Be21_n}) contributes constant and  linear terms with respect to $\tau$, while the first sum contributes only terms that are at least order 2 with respect to $\tau$. The result is that 
\begin{eqnarray}
&& \hspace{-2cm} (\nu_2 - \mu_2) \tau B_{2,\be_1^*}\left( \frac{u}{|u|} \right)\prod_{\substack{j \geq 3 \\ j \in \Lambda}} B_{j,\be(j)}\left(\frac{u}{|u|}\right) \nonumber \\
&=&
 - \sum_{m = 2}^d \nu_m |u|^{m-1} B_{m,\be_1^*}\left( \frac{u}{|u|} \right) \prod_{\substack{j \geq 3 \\ j \in \Lambda}} B_{j,\be(j)}\left(\frac{u}{|u|}\right) \nonumber \\
&&-\; \tau \sum_{m = 3}^d \nu_m (m-1) |u|^{m-2} B_{m,\be_1^*}\left( \frac{u}{|u|} \right)  \prod_{\substack{j \geq 3 \\ j \in \Lambda}} B_{j,\be(j)}\left(\frac{u}{|u|}\right) \nonumber \\
&&+ \; \tau^2 S_{\nu,u}^{(2)}(\tau) + O(r^{\ep_1}),  \label{S_nu_u_2b_2_n}
\end{eqnarray}
where $S_{\nu,u}^{(2)}(\tau)$ is a polynomial in $\tau$ whose coefficients depend only on $\nu$ and $u$ (but not $\mu$).

We now feed this information into an analysis of $C[\sig^\be](u,\tau)$ for some $\be$ with $|\be|=2$ to be chosen later. For now we fix any $\be$ with $|\be|=2$ and recall that by definition, for this $\beta$
\begin{multline}\label{C_sig2_case2b_n}
C[\sig^\be](u,\tau) = \sum_{m=2}^d  (\nu_m - \mu_m) \tau^{m-2} B_{m,\be}\left(\frac{u}{|u|}\right)  \\
 + \;
\sum_{m=3}^d \nu_m( (|u| + \tau)^{m-2} - \tau^{m-2})B_{m,\be} \left(\frac{u}{|u|}\right).
\end{multline}
(Note that unlike Case B1, the first sum here begins with $m = 2$.) We want to use (\ref{S_nu_u_2b_1_n}) and (\ref{S_nu_u_2b_2_n}) to make the first sum on the right hand side independent of $\mu$: hence we multiply (\ref{C_sig2_case2b_n}) through by 
\[
\tau B_{2,\be_1^*}\left(\frac{u}{|u|}\right) \prod_{\substack{j \geq 3 \\ j \in \Lambda}} B_{j,\be(j)}\left(\frac{u}{|u|}\right),\]
 and substitute  (\ref{S_nu_u_2b_2_n}) to treat the term including $(\nu_2 - \mu_2)$ and (\ref{S_nu_u_2b_1_n}) to treat terms including $(\nu_m - \mu_m)$ for all $m \geq 3$. One then concludes that
\beq\label{C_exprB2}
C[\sig^\be](u,\tau) \tau B_{2,\be_1^*}\left(\frac{u}{|u|}\right)  \prod_{\substack{j \geq 3 \\ j \in \Lambda}} B_{j,\be(j)}\left(\frac{u}{|u|}\right) 
= R_{\nu,u}^{(\be)}(\tau) + O(r^{\ep_1}), 
\eeq
where $R_{\nu,u}^{(\be)}(\tau)$ is a polynomial in $\tau$ whose coefficients depends only on $\nu$ and $u$ (but not $\mu$). (The superscript $\be$ reflects that this comes from $C[\sig^\be](u,\tau)$.) 

In fact we may compute the coefficient of $\tau$ in $R_{\nu,u}^{(\be)}(\tau)$ explicitly. We need only note that the term in $R_{\nu,u}^{(\be)}(\tau)$ that is linear in $\tau$ comes from the terms in $C[\sig^\be](u,\tau)$ that are constant with respect to $\tau$. Then we use the fact that the first term on the right hand side of (\ref{C_sig2_case2b_n}) contributes a constant with respect to $\tau$ with the $m=2$ summand; the second term on the right hand side of (\ref{C_sig2_case2b_n}) contributes a constant term in $\tau$ for each $3 \leq m \leq d$.
This shows that the coefficient of $\tau$ in $R_{\nu,u}^{(\be)}(\tau)$ is given by
\begin{multline}\label{B2_coeff_tau}
-\sum_{m = 3}^d \nu_m (m-1) |u|^{m-2} B_{m,\be_1^*}\left( \frac{u}{|u|} \right) B_{2,\be}\left(\frac{u}{|u|}\right)  \prod_{\substack{j \geq 3 \\ j \in \Lambda}} B_{j,\be(j)}\left(\frac{u}{|u|}\right) \\
+ \sum_{m=3}^d \nu_m |u|^{m-2} B_{m,\be} \left(\frac{u}{|u|}\right) B_{2,\be_1^*}  \left(\frac{u}{|u|}\right)  \prod_{\substack{j \geq 3 \\ j \in \Lambda}} B_{j,\be(j)}\left(\frac{u}{|u|}\right).
\end{multline}
Using the homogeneity property of $B_{m,\be}$, we see that (\ref{B2_coeff_tau}) simplifies to
\[
|u|^{-s_4} \left( \prod_{\substack{j \geq 3 \\ j \in \Lambda}} B_{j,\be(j)}(u)\right) \sum_{m = 3}^d \nu_m [-(m-1) B_{m,\be_1^*}(u) B_{2,\be}(u) + B_{m,\be}(u) B_{2,\be_1^*}(u)] 
\] 
where 
\[ s_4 = 4 + \sum_{\bstack{j \geq 3}{ j \in \Lambda}} j.\]
Let $W_{\nu}^{(\be)}(u)$ be the polynomial in $u$ defined by
\beq\label{W_nu_def_case2b_n}
W_{\nu}^{(\be)}(u) := 
\left( \prod_{\substack{j \geq 3 \\ j \in \Lambda}} B_{j,\be(j)}(u)\right)  \sum_{m = 3}^d \nu_m [-(m-1) B_{m,\be_1^*}(u) B_{2,\be}(u) + B_{m,\be}(u) B_{2,\be_1^*}(u)],
\eeq
which has total degree at most 
\[ s_5 = d+2 + \sum_{\bstack{j \geq 3}{j \in \Lambda}} j.\]
Then we can conclude that
\[
\llbracket R_{\nu,u}^{(\be)} \rrbracket_{\tau} \geq |u|^{-s_4} |W_{\nu}^{(\be)}(u)| \geq  C |W_{\nu}^{(\be)}(u)| .
\]
Recall that in the current case, we assumed that there exists an $m \in \Lambda$, say $m_0$, with $p_{m_0}(y)$ parabolic, such that $|\nu_{m_0}| \simeq r$.  Since we assumed $p_2(y) \ne C|y|^2$ for any $C \ne 0$ in Proposition~\ref{prop_Ksharp_n}, we know that $m_0 \geq 4$. Thus we can
single out the contribution to $W_{\nu}^{(\be)}(u)$ from $m_0$, which by (\ref{W_nu_def_case2b_n}) is 
 \beq\label{W_m0_contribution}
 \nu_{m_0} \left( \prod_{\substack{j \geq 3 \\ j \in \Lambda}} B_{j,\be(j)}(u)\right)  [-(m_0-1) B_{m_0,\be_1^*}(u) B_{2,\be}(u) + B_{m_0,\be}(u) B_{2,\be_1^*}(u)].
 \eeq
 It is clear that the total degree of this contribution is distinct from that of all other terms in $W_{\nu}^{(\be)}(u)$.
Our aim now is to show that we can pick a particular index $\be$ with $|\be|=2$ such that (\ref{W_m0_contribution}) is a nonzero polynomial with respect to $u$, so that in particular we can conclude that $W_{\nu}^{(\be)}(u)$ contains a coefficient of size $|\nu_{m_0}| \simeq r$.

Since $m_0 \in \Lambda$ and $p_{m_0}$ is parabolic we know by Statement \ref{item_parabolic_n} of Lemma \ref{lemma_B_n} that $B_{m_0,\be_1^*} \equiv 0$ since $|\be_1^*|=1$, which is odd. 
Thus the contribution (\ref{W_m0_contribution}) to $W_{\nu}^{(\be)}(u)$ from $m_0$ is in fact precisely
 \beq\label{m0_n}
 \nu_{m_0} \left( \prod_{\substack{j \geq 3 \\ j \in \Lambda}} B_{j,\be(j)}(u)\right) B_{m_0,\be}(u) B_{2,\be_1^*}(u).
 \eeq
We recall that by construction,  $B_{2,\be_1^*} \not \equiv 0$ and $B_{j,\be(j)}(u)\not\equiv 0$ for all $3 \leq j \leq d$, $j \in \Lambda$.
Next, we note by Statement \ref{item_m2_n} of Lemma \ref{lemma_B_n} that there exists $\be$ with $|\be|=2$ such that $B_{m_0,\be}$ is not the zero polynomial; we will call this distinguished index $\be_2^*$.

We now define our exceptional sets as follows: recall that $0<\ep_1<\ep_2<1$, and for the choice $\be = \be_2^*$ with $|\be_2^*|=2$ distinguished above, set
\[
G^{\nu} := \{ u \in B_2(\R^n) \colon |W_{\nu}^{(\be_2^*)}(u)| \leq r^{\ep_2} \},
\]
and for $u \in B_2 \setminus G^{\nu}$, let 
\[
F^{\nu}_u := \{ \tau \in B_1(\R) \colon |R_{\nu,u}^{(\be_2^*)}(\tau)| \leq C_0 r^{\ep_1} \}
\]
for some large absolute constant $C_0$ to be determined. Also define $F^{\nu}_u := \emptyset$ if $u \in G^{\nu}$. Then for all $(u,\tau) \in B_2(\R^n) \times B_1(\R)$ with $u \notin G^{\nu}$, $\tau \notin F^{\nu}_u$, we have
\[
|R_{\nu,u}^{(\be_2^*)}(\tau)| \geq C_0 r^{\ep_1},
\]
so for such $(u,\tau)$ we conclude by (\ref{C_exprB2}) that
\begin{eqnarray*}
|C[\sig^{\be_2^*}](u,\tau)| 
&\geq & C_B^{-|\Lambda|} \left|C[\sig^{\be_2^*}](u,\tau) \tau B_{2,\be_1^*}\left(\frac{u}{|u|}\right) \prod_{\substack{j \geq 3 \\ j \in \Lambda}} B_{j,\be(j)}\left(\frac{u}{|u|}\right)\right| \\
& \geq& C_B^{-|\Lambda|}  \left( |R_{\nu,u}^{(\be_2^*)}(\tau)| - O(r^{\ep_1}) \right) \\
&\geq &C_B^{-|\Lambda|}  \left(C_0 r^{\ep_1} - O(r^{\ep_1}) \right) \\
& \geq & C_B^{-|\Lambda|} r^{\ep_1} 
\end{eqnarray*}
if $C_0$ is sufficiently large. It follows from the van der Corput estimate of Lemma \ref{lemma_Prop2.1} that
\[
|K_{\sharp,n}^{\nu,\mu}(u,\tau)| \leq C r^{-\ep_1/d} \quad \text{if $u \notin G^{\nu}$ and $\tau \notin F_u^{\nu}$}.
\]

On the other hand, we will now show that $G^{\nu}$ and $F_u^{\nu}$ are sets of small measures. First, $G^{\nu}$ is the set of $u \in B_2$ where $W_{\nu}^{(\be_2^*)}(u)$ is small, and in (\ref{W_nu_def_case2b_n}), $W_{\nu}^{(\be_2^*)}(u)$ is represented as a sum of homogeneous polynomials of different total degrees. 
We have already noted above that the contribution to $W_{\nu}^{(\be_2^*)}$ from $m_0$ is the term (\ref{m0_n}) specialized to $\be = \be_2^*$, and that this has total degree different from all other terms in $W_{\nu}^{(\be_2^*)}$. We have also noted that by construction, the term (\ref{m0_n}) with $\be = \be_2^*$ is not the zero polynomial; thus in particular it contributes a term in the polynomial $W_{\nu}^{(\be_2^*)}$ that has coefficient $|\nu_{m_0}| \simeq r$.
We therefore see that
\[
\llbracket W_{\nu}^{(\be_2^*)} \rrbracket_u \geq C |\nu_{m_0}| \geq C r,
\]
which implies by Lemma \ref{lemma_Prop2.2_0} that
\[
|G^{\nu}| \leq C \left( \frac{r^{\ep_2}}{r} \right)^{1/s_5} = C r^{-(1-\ep_2)/s_5}.
\]
Furthermore, if $u \in B_2 \setminus G^{\nu}$, then 
\[
\llbracket R_{\nu,u}^{(\be_2^*)} \rrbracket_{\tau} 
\geq C |W_{\nu}^{(\be_2^*)}(u)|
\geq C r^{\ep_2}.
\] 
Thus if $u \in B_2 \setminus G^{\nu}$, then by Lemma \ref{lemma_Prop2.2_0},
\[
|F_u^{\nu}| \leq C \left( \frac{C_0 r^{\ep_1}}{r^{\ep_2}} \right)^{1/d} = C r^{-(\ep_2-\ep_1)/d}.
\] 
The same inequality is clearly true if $u \in G^{\nu}$. Thus we have $|F_u^{\nu}|$ being small for all $u \in B_2$. This concludes the proof of Proposition~\ref{prop_Ksharp_n} in this final case.

\subsection{Proof of Lemma \ref{lemma_A_ab} for $A_{\al,\be}$}\label{sec_A_B_lemmas}
We now turn to the proof of the key properties of $A_{\al,\be}$ given in Lemma \ref{lemma_A_ab}. Recall the expansion (\ref{defining_A_rel}) that defines $A_{\al,\be}(w)$:
   \beq\label{A_expansion}
   \sum_{ |\be| \leq |\al|} A_{\al,\be} (w)  \sig^{\be} = ( w_n + w^{(n)}\cdot \sig)^{\al_n} (w^{(n)}- w_n \sig)^{\al^{(n)}}.
   \eeq
   We recall that $w \in \R^n$ while $\sig = (\sig_1, \ldots, \sig_{n-1}) \in \R^{n-1}$, so that the multi-index $\al$ is $n$-dimensional while the multi-index $\be$ is $(n-1)$-dimensional.
   
   To prove Statement \ref{item_A_1}, we use the fact that $|\be|=1$ to compute $A_{\al,\be}$ explicitly. Suppose that $\be = e_j$, where $e_j = (0,\ldots, 1, \ldots, 0)$ is the $j$-th unit vector, with $1 \leq j \leq n-1$. Then to compute $A_{\al,\be}(w)$, we must pick out the coefficient of $\sig^\be = \sig_j$ on the right hand side of the expansion (\ref{A_expansion}), namely\xtra{There are two ways a term in the expansion of the product can be precisely of the form $\sig_j$: first, if it consists of choosing one $\sig_j$-term in the first factor and choosing the remaining $\al_n-1$ terms in the factor to be $w_n$ (there are $\al_n$ ways to do this), and choosing all $w^{(n)}$-terms in all the remaining factors. Second, a $\sig_j$ monomial can be composed by choosing $w^{(n)}$-terms in all factors but the $\al_j$ factor, in which we choose one $\sig_j$-term and the remaining $\al_j-1$ terms are chosen to be $w_j$ (there are $\al_j$ ways to do this).} 
\[ \al_n w_j w_n^{\al_n-1} w_1^{\al_1} \cdots w_{n-1}^{\al_{n-1}} 
	+  \al_j w_n^{\al_n} w_1^{\al_1} \cdots (-w_n)w_{j}^{\al_j-1} \cdots w_{n-1}^{\al_{n-1}}.
\]
This coefficient of $\sig_j$ is the polynomial $A_{\al,e_j}(w)$ we seek; we write it more efficiently as
\[ A_{\al,e_j}(w) =\al_n w^{\al + e_j  - e_n} - \al_j w^{\al-e_j +e_n}. \]
(We pause here to note that if we were working with $K^{\nu,\mu}_{\sharp,l}(u,\tau)$ for any $1 \leq l \leq n$, we would consider all $1 \leq j \leq n$ with $j \neq l$ and the analogous expression would clearly be 
\[A_{\al,e_j}(w) =\al_l w^{\al + e_j  - e_l} - \al_j w^{\al-e_j +e_l}. \]
Thus although we appear to be privileging the $n$-th component, nothing in the proof depends more than notationally on this choice.)
Statement \ref{item_A_2} follows in a similar fashion from examining the  expansion (\ref{A_expansion}).

To prove Statement \ref{item_A_be_al}, fix $\al \in \Z^n_{\geq 0}$ and fix any $\be \in \Z^{n-1}_{\geq 0}$ with $|\be| = |\al|$. To single out the term $A_{\al,\be}(w)$ from the left hand side of (\ref{A_expansion}), we must simply pick out the coefficient of $\sig^{\be}$ on the right hand side. Expand the right hand side as 
 \beq\label{A_expansion1}  ( w_n + w_1\sig_1 + \cdots + w_{n-1} \sig_{n-1})^{\al_n} (w_1 - w_n \sig_1)^{\al_1} \cdots (w_{n-1}  - w_n \sig_{n-1})^{\al_{n-1}}.
 \eeq
To prove that $A_{\al,\be}$ is a monomial, we must verify that no more than one term in this expansion is of the form $\sig^\be$. Clearly, since $|\al|=|\be|$, in order for the total degree of such a monomial to reach $|\be|$, we must always choose $\al_j$ copies of the $\sig_j$-factor in each of the last $n-1$ terms, so that the product of the last $n-1$ terms contributes 
\beq\label{first_monomial}
(-1)^{|\al^{(n)}|}w_n^{|\al^{(n)}|}\sig^{\al^{(n)}}.
\eeq
The remaining part of the $\sig^\be$ monomial comes from the first factor, which we write as 
\begin{multline}\label{w_gamma_sum}
  ( w_n + w_1\sig_1 + \cdots + w_{n-1} \sig_{n-1})^{\al_n} \\
	= \sum_{\ga_1 + \cdots + \ga_n=\al_n} \binom{\al_n}{\ga_1,\ldots, \ga_n} (w_1\sig_1)^{\ga_1} \cdots (w_{n-1}\sig_{n-1})^{\ga_{n-1}}(w_n)^{\ga_n}.
\end{multline}
Here we are using the usual multinomial coefficient, defined for $k_1 + \cdots + k_n = k$ by
\[ \binom{k}{k_1,\ldots, k_n} = \frac{k!}{k_1! \cdots k_n!}.\]

A multi-index $\ga \in \Z^n_{\geq 0}$ in the sum (\ref{w_gamma_sum}) will contribute to the $\sig^\be$ monomial (when combined with the second factor (\ref{first_monomial})) if and only if it solves the equations 
\begin{eqnarray}
|\ga| & = & \al_n \label{ga_1_n} \\
\ga^{(n)} + \al^{(n)} &=& \be. \label{ga_2_n}
\end{eqnarray}
There is at most one solution to this system. 
Indeed, if $\al_j > \be_j$ for any $1 \leq j \leq n-1$ then there is no solution $\ga$, since (\ref{ga_2_n}) cannot be satisfied by any $\ga^{(n)} \in \Z_{\geq 0}^{n-1}$, and in this case $ A_{\al,\be} \equiv 0$ (that is, the constant $C_{\al,\be}$ is zero). 
But if 
\beq\label{albe}
\al_j \leq \be_j \quad \text{for all $1 \leq j \leq n-1$},
\eeq
then we may solve (\ref{ga_2_n}) uniquely for $\ga^{(n)}$. The condition (\ref{albe}) along with $|\be| = |\al|$ also guarantees that 
\[ |\ga^{(n)}| =  |\be - \al^{(n)}| = |\be| - |\al^{(n)}| = \al_n ,\] 
which implies in (\ref{ga_1_n}) that $\ga_n=0$. We have shown that if (\ref{albe}) holds, there is precisely the unique solution $\ga = (\be - \al^{(n)},0)$ for the system (\ref{ga_1_n})-(\ref{ga_2_n}) and consequently $A_{\al,\be}$ is the monomial 
\[ A_{\al,\be} (w) = C_{\al,\be}(-1)^{|\al^{(n)}|}w_n^{|\al^{(n)}|}(w^{(n)})^{\be - \al^{(n)}} 
 \]
with combinatorial coefficient
\[ C_{\al,\be} =  \binom{\al_n}{(\be_1 - \al_1),\ldots, (\be_{n-1} - \al_{n-1}),0} \neq 0.\]
In particular, $A_{\al,\be}$ is a monomial with nonzero coefficient for any $\be$ such that $|\be| = |\al|$ and $\be_j \geq \al_j$ for all $1 \leq j \leq n-1$, and otherwise is the zero polynomial. (Certainly there is at least one such $\be$ for every $\al$.)

\subsection{Proof of Lemma \ref{lemma_B_n} for $B_{m,\be}$, Statements \ref{item_dependence_n} - \ref{item_parabolic_n}}\label{sec_B_lemma_1}
Statements \ref{item_dependence_n} and \ref{item_bound_n} are clear from the definition of $B_{m,\be}$. 
Statement \ref{item_hom_n} is a consequence of the fact, already observed, that for each $\al,\be$, the polynomial $A_{\al,\be}(w)$ is homogeneous of degree $|\al|$. To see Statement \ref{item_zero_n}, we fix an $m \in \Lambda$ and recall that 
\[ B_{m,\be}(w) = \sum_{|\al|=m} c_\al A_{\al,\be}(w).\]
By Lemma \ref{lemma_A_ab}, for each $|\al|=|\be|$, $A_{\al,\be}(w)$ is a monomial, and in particular we can write
\[ B_{m,\be}(w) = \sum_{|\al|=m} c_\al C_{\al,\be} w_n^{|\al^{(n)}|} (w^{(n)})^{\be - \al^{(n)}},\]
for non-negative combinatorial constants $C_{\al,\be}$. For $\be$ fixed, each monomial 
\[w_n^{|\al^{(n)}|} (w^{(n)})^{\be - \al^{(n)}}\]
 with $|\al|=|\be|$ is distinct. Thus it suffices to show that there is some $\be$ with $|\be|=m$ and some $|\al|=m$ such that $c_\al C_{\al,\be} \neq 0$. We first note that since $m \in \Lambda$, there is some $\al$ with $|\al|=m$ such that $c_\al \neq 0$; we will call this $\al(m)$. It now suffices to find $\be$ with $|\be|=m$ such that $C_{\al(m),\be} \neq 0$; such a choice of $\be$ is in fact guaranteed by Lemma \ref{lemma_A_ab}. This proves Statement \ref{item_zero_n}.

We now turn to Statement \ref{item_parabolic_n}.  Suppose $p_m(y)$ is parabolic for some index $2 \leq m \leq d$, say $p_m(y) = C|y|^m$ for some nonzero constant $C$; in particular $m$ must be even, so we write $m=2k$. Recall that (\ref{phase_expansion_special_n}) holds for all choices of $\nu$ and $\mu$; we will apply that equation to expand $C|z|^m$ alone by setting $\nu_j = 0$ for all $j$, $\mu_m = 1$, and $\mu_j = 0$ for all $j \ne m$. With these choices in (\ref{phase_expansion_special_n}), we see that 
\beq\label{dya_exp_n}
-C|z|^m = -\sum_{0 \leq |\be| \leq m} \tau^{m-|\be|} B_{m,\be} \left(\frac{u}{|u|}\right) \sig^\be,
\eeq
in which $z$ is defined implicitly by (\ref{z_rel1}) and (\ref{z_rel2}) (with index $l=n$). On the other hand, if $z$ is defined by (\ref{z_rel1}) and (\ref{z_rel2}), then 
\[ |z|^{2k} = (z_n^2 + |z^{(n)}|^2)^k = \frac{1}{|u|^{2k}} \left\{ \tau^2 |u|^2 + (u^{(n)} \cdot \sig)^2 + u_n^2 |\sig|^2 \right\}^k.
\]
If we denote  the right hand side by $T(\sigma)$, say, it is clear by inspection that $T(\sig)$ is an even polynomial in $\sigma$, that is to say $T(\sigma) = T(-\sigma)$. 
Now fix any multi-index $\be$; the coefficient of $\sigma^\be$ in $T(\sigma)$ is of course $\left. \partial^\be T(\sigma) \right|_{\sigma=0}$. But since $T$ is an even polynomial, we see that 
\[ \left. \partial^\be T(\sigma) \right|_{\sigma=0} =\left. \partial^\be (T(-\sigma) )\right|_{\sigma=0} = (-1)^{|\be|} \left. (\partial^\be T)(-\sigma) \right|_{\sigma=0} = (-1)^{|\be|} \left. (\partial^\be T)(\sigma) \right|_{\sigma=0},\]
which shows that the coefficient of $\sigma^\be$ in $T(\sigma)$ must be zero whenever $|\be|$ is odd.
In the expansion (\ref{dya_exp_n}) this shows that $\tau^{m-|\be|}B_{m,\be}(u/|u|)$ is identically zero in $\tau, u$ whenever $|\be|$ is odd, which implies that $B_{m,\be}(w)$ is the zero polynomial whenever $|\be|$ is odd.

\xtra{We can also argue more explicitly as follows:
Raising the expansion for $|z|^2$ to the $k$-th power, we have
\beq\label{z_exp_n}
|z|^m = |z|^{2k} =  \frac{1}{|u|^{2k}} \sum_{k_1 + k_2+k_3 =k} \binom{k}{k_1,k_2,k_3} (\tau |u|)^{2k_1} (u^{(n)} \cdot \sig)^{2k_2} (u_n|\sig|)^{2k_3},
\eeq
where the sum ranges over all partitions $k_1 + k_2 + k_3 = k$ with $k_i \geq 0$ for $i=1,2,3$ and we use the notation standard in multi-nomial expansions for the combinatorial coefficient, namely for $k_1 + \cdots +k_r=k$ we set
\[ \binom{k}{k_1, \ldots, k_r} = \frac{k !}{k_1! \cdots k_r!}.\]
The contribution to $|z|^m$ from a fixed partition $k_1 + k_2 + k_3$ is a monomial (ignoring constants and some factors of $u$ for the moment) of the form 
\[ (u^{(n)} \cdot \sig)^{2k_2} |\sig|^{2k_3} =  \left\{ \sum_{j_1 + \cdots +j_{n-1} = 2k_2} \binom{2k_2}{j_1,\ldots, j_{n-1}} (u_1\sig_1)^{j_1} \cdots (u_{n-1}\sig_{n-1})^{j_{n-1}}  \right\} (\sig_1^2 + \cdots + \sig_{n_1}^2)^{k_3}.
\]
Fix any monomial $\sig^\be$ arising in this expression; we can write it as $\sig^\be = \sig^{\be_1} \sig^{\be_2}$ where $\sig^{\be_1}$ is the contribution of the first factor and $\sig^{\be_2}$ is the contribution of the second factor; then clearly  $|\be_1| = j_1 + \cdots + j_{n-1} = 2k_2$ which is even, and $|\be_2|$ is visibly even, hence $|\be|$ is even. We may conclude that any monomial $\sig^\be$ appearing in the expression (\ref{z_exp_n}) for $|z|^{m}$ must have $|\be|$ even. Thus 
it follows from (\ref{dya_exp_n}) that for any $\be$ with $|\be|$ odd, we have $\tau^{m-|\be|} B_{m,\be}(u/|u|) \equiv 0$ for all $\tau$ and $u$, which shows that $B_{m,\be}(w)$ is the zero polynomial. 
}

\subsection{Proof of Lemma \ref{lemma_B_n} for $B_{m,\be}$, Statement \ref{item_m1_n}}\label{sec_B_lemma_2}
We will prove that if $m \in \Lambda$ and $p_m(y)$ is not parabolic, then there exists $\be$ with $|\be|=1$ such that $B_{m,\be}$ is not the zero polynomial.
In fact, it is easier to prove the equivalent statement that if $m \in \Lambda$ and $B_{m,\be} \equiv 0$ for all $|\be|=1$, then $p_m$ is parabolic. We recall that being parabolic puts a constraint on the relationships between the coefficients $c_\al$ (fixed once and for all) in the definition $p_m(y) = \sum_{|\al|=m} c_\al y^\al$. 
Precisely, if $p_m(y)= C|y|^m$ for $m=2k$, then 
\[ C|y|^{2k} = C(y_1^2 + \cdots + y_n^2)^k = C \sum_{k_1 + k_2 + \cdots +k_n = k} \binom{k}{k_1,k_2, \ldots, k_n} y_1^{2k_1} y_2^{2k_2} \cdots y_n^{2k_n}.\]
Thus we see for $m=2k \in \Lambda$, $p_{2k}(y)$ being parabolic is characterized the property that there exists a nonzero constant $C$ such that for any partition $k_1 + \cdots+ k_n = k$ and corresponding multi-index $\al = (2k_1,\ldots, 2k_n)$, the coefficient $c_\al$ must satisfy
\[ c_\al = C \binom{k}{k_1,k_2, \ldots, k_n} ,\]
and for all $\al$ with $|\al|=2k$ that have an odd entry, $c_\al=0$.

We now fix any $m \in \Lambda$ and assume that $B_{m,\be} \equiv 0$ for all $|\be|=1$. We will show that $p_m$ must be parabolic. Fix any $\be$ with $|\be|=1$, which we will denote by $\be = e_j$ for some $1 \leq j \leq n-1$.
We apply Statement \ref{item_A_1} of Lemma \ref{lemma_A_ab} to compute that
\[ B_{m,e_j}(w) = \sum_{|\al|=m} c_\al A_{\al,e_j} (w) = \sum_{|\al|=m} c_\al ( \al_n w^{\al + e_j  - e_n} - \al_j w^{\al-e_j +e_n}) .\]
We now re-write this by grouping coefficients for each fixed monomial $w^\rho$, as 
\beq\label{Bmej}
 B_{m,e_j}(w) = \sum_{|\rho|=m} [c_{\rho+ e_n -e_j} (\rho_n+1) - c_{\rho - e_n + e_j} (\rho_j+1)] w^\rho.
\eeq
By notational convention, if $\rho_j=0$ then $c_{\rho + e_n - e_j}=0$ and if $\rho_n=0$ then $c_{\rho - e_n +e_j}=0$, so the corresponding terms do not actually appear in the expression above.

Under our assumption that $B_{m,e_j} \equiv 0$ for all $1 \leq j \leq n-1$, we see that for each $j$ all the coefficients in (\ref{Bmej}) are identically zero, so that for all $|\rho|=m$ and all $1 \leq j \leq n-1$ we have 
\beq\label{Bmej_relation}
c_{\rho + e_n -e_j} (\rho_n+1) - c_{\rho - e_n + e_j} (\rho_j+1) =0.
\eeq
Again, if either $\rho_j$ or $\rho_n$ is zero, then the corresponding term does not appear in the expression above. First we note that since $(\rho_n+1)$ and $(\rho_j+1)$ will never vanish for multi-indices $\rho \in \Z_{\geq 0}^n$, (\ref{Bmej_relation}) shows that the coefficients $c_{\rho+e_n-e_j}$ and $c_{\rho-e_n+e_j}$ are either both zero or both nonzero. Second, we note that (\ref{Bmej_relation}) provides relations between coefficients with indices sharing a certain parity property.
For any $m \geq 1$, let $\Ga_{\odd}(m)$ denote  the set of multi-indices $\ga \in \Z_{\geq 0}^n$ with $|\ga|=m$ such that at least one coordinate of $\ga$ is odd. We will prove:
\begin{lemma}\label{lemma_cCD_odd}
Let $m \in \Lambda$  and suppose that $B_{m,e_j} \con 0$ for all $1 \leq j \leq n-1$.
Then $c_{\ga} =0$ for all $ \ga \in \Ga_{\odd}(m)$.
\end{lemma}

Assume this for the moment. If $m$ is odd, then every $\ga$ with $|\ga|=m$ contains an odd entry, so as a consequence of Lemma \ref{lemma_cCD_odd}, if $m \in \Lambda$ yet $B_{m,e_j}\con 0$ for all $1 \leq j \leq n-1$ then we must have $m$ even. Thus assuming $m=2k$ is even,
we let $\Ga_\even(m)$  denote the set of multi-indices $\ga \in \Z_{\geq 0}^n$ with $|\ga|=m$ such that all entries in $\ga$ are even; in this case we will write $\ga = (2\ga_1, \ldots, 2 \ga_n)$. Then we will deduce from the identity (\ref{Bmej_relation}) the following lemma:
\begin{lemma}\label{lemma_cCD}
Let $m \in \Lambda$ with $m=2k$ for some $k \geq 1$ and suppose that $B_{m,e_j} \con 0$ for all $1 \leq j \leq n-1$.
Let $\ga^*$ denote the element  $(0,0,\ldots,2k) \in \Ga_\even(2k)$ and define the constant
\[ C :=c_{\ga^*}.\]
Then for  all $\ga = (2\ga_1, \ldots, 2\ga_n) \in \Ga_\even(2k)$, the coefficient $c_\ga$ of $p_{2k}(y)$ must satisfy
\beq\label{cDn}
 c_\ga= C  \binom{k}{\ga_1, \ldots, \ga_n}.
 \eeq
\end{lemma}
With Lemmas \ref{lemma_cCD_odd} and \ref{lemma_cCD} in hand,  we can deduce Statement  \ref{item_m1_n}
of Lemma \ref{lemma_B_n} immediately.
Indeed, we have already seen that if $m \in \Lambda$ but $B_{m,\be} \con 0$ for all $|\be|=1$ then $m$ must be even, and now for $m=2k$ we may conclude
 \[ p_m(y) = \sum_{\al \in \Ga_\even(m)} c_\al y^\al = C\sum_{\bstack{|\al|=2k}{\al = (2\al_1,\ldots, 2\al_n)}} \binom{ k}{\al_1,\ldots, \al_n} y^{2\al_1} \cdots y^{2\al_n} = C|y|^{2k}.\]
 From this we also deduce that $C$ must be nonzero, since otherwise we would have $p_m \equiv 0$, which would contradict $m \in \Lambda$.

 We now prove Lemmas \ref{lemma_cCD_odd} and \ref{lemma_cCD}.
 To prove Lemma \ref{lemma_cCD_odd} in a notationally economical fashion, we define the equivalence relation
 $\ga \sim \ga'$ to indicate that $c_\ga =0$ if and only if $c_{\ga'}=0$. 
 Suppose $\ga \in \Ga_\odd(m)$ is given; we will prove first that $\ga \sim \ga'$ for some $\ga'$ with an entry equal to 1, and then that for any $\ga'$ with an entry equal to 1 we necessarily have $c_{\ga'}=0$.
 
For the first statement,  we assume that $\ga$ is given with $\ga_j$ odd; there are two cases, depending on whether $1 \leq j \leq n-1$ or $j=n$. We first assume that $1 \leq j \leq n-1$, in which case we will observe that $\ga \sim \ga - 2\theta e_j + 2\theta e_n$ for any non-negative integer $\theta$ such that $2\theta \leq \ga_j$.
This follows from a simple induction argument; the base case $\theta=0$ is clear. Assuming the induction hypothesis that $\ga \sim \ga - 2(\theta-1)e_j + 2(\theta-1)e_n$ for some $\theta \geq 1$, we apply (\ref{Bmej_relation}) with the choice
$\rho  = \ga -(2\theta-1) e_j + (2\theta-1) e_n$ to see that 
\begin{multline*}
c_{\rho + e_n -e_j} (\rho_n+1) - c_{\rho - e_n + e_j} (\rho_j+1) \\
=c_{\ga +2\theta e_n - 2\theta e_j}(\ga_n +2\theta) - c_{\ga +2(\theta-1)e_n - 2(\theta-1)e_j}(\ga_j-2\theta +2)  =0.
\end{multline*}
Since $(\ga_j - 2\theta + 2) \neq 0$ and $(\ga_n + 2\theta )\neq 0$, it follows  that $ \ga - 2(\theta-1)e_j + 2(\theta-1)e_n \sim  \ga - 2\theta e_j + 2\theta e_n $, which completes the induction.

With this fact in hand, given $\ga \in \Ga_\odd(m)$ with $\ga_j = 2k+1$ for some $k \geq 0$, we may conclude that $\ga \sim \ga' :=\ga - 2k e_j + 2ke_n$, where $\ga'$ has $j$-th coordinate equal to 1 by construction. 
If we instead had $\ga_j$ odd for $j=n$, we would similarly use induction to show that  for any $1 \leq i \leq j-1$, $\ga \sim \ga - 2\theta e_n + 2\theta e_i$, for any non-negative integer $\theta$ such that $2\theta \leq \ga_n$; then we would apply this to show $\ga \sim \ga'$ for some $\ga'$ with $\ga_n'=1$.

Next we show that any $\ga$ with an entry equal to 1 has $c_\ga=0$.  Assume that $\ga$ is such that $\ga_j=1$ for some $1 \leq j \leq n$. There are again two cases to consider: $1 \leq j \leq n-1$ and $j=n$. 
In the case that $\ga_j=1$ for some $1 \leq j \leq n-1$, we apply (\ref{Bmej_relation}) with the choice $\rho = \ga+e_n-e_j$, to conclude that 
 \[ c_\ga = c_{\rho-e_n+e_j}=0.\]
  (Here we are using the fact that with our particular choice of $\rho$, the $j$-th entry of $\rho+e_n-e_j$ is negative, so that by convention the coefficient $c_{\rho+e_n-e_j}=0$.)
  Similarly, in the case that $j=n$, we choose any $1 \leq i \leq n-1$ we like and apply (\ref{Bmej_relation}) with the choice $\rho = \ga -e_n+e_i$, to conclude that 
 \[ c_\ga = c_{\rho+e_n-e_i} =0.
 \]
This completes the proof of Lemma \ref{lemma_cCD_odd}.

In Lemma \ref{lemma_cCD}, we only consider $m=2k$ even and indices $\ga \in \Gamma_{\even}(2k)$, which we denote by $\ga = (2\ga_1,\ldots, 2\ga_n)$. For $\ga, \ga' \in \Gamma_\even(2k)$, we define the equivalence relation $\ga \sim \ga'$ to represent that $c_\ga$ satisfies (\ref{cDn}) if and only if $c_{\ga'}$ satisfies (\ref{cDn}).

We note first of all that (\ref{cDn}) holds true for $\ga^*$,
since 
\[ c_{\ga^*} = C \binom{k}{0, \cdots 0, k} = C \frac{k!}{0! \cdots 0! k!} = C.\]
Thus it suffices to prove $\ga^* \sim \ga$ for any $\ga \in \Ga_\even(2k)$. Thus we fix any $\ga  = (2\ga_1,\ldots, 2\ga_n) \in \Ga_\even(2k)$ and note that
\[ \ga  = \ga^* + 2 \ga_1 e_1 + \cdots + 2\ga_{n-1}e_{n-1} - (2k - 2\ga_n)e_n
	= \ga^* + \sum_{j=1}^{n-1}2 \ga_j e_j - (\sum_{j=1}^{n-1} 2\ga_j)e_n .\]
We will prove:
\begin{lemma}\label{lemma_imply}
Given $\ga \in \Ga_\even(2k)$, suppose $\ga' = \ga + 2\theta e_j - 2 \theta e_n$ for some $1 \leq j \leq n-1$ and some non-negative integer $\theta$ with $2\theta \leq 2\ga_n$. Then $\ga \sim \ga'$.
\end{lemma}
 With this lemma in hand, we see that 
\begin{multline*}
 \ga^* \sim \ga^* + 2\ga_1 e_1 - 2\ga_1 e_n \sim (\ga^* + 2\ga_1 e_1 - 2\ga_1 e_n) + 2\ga_2e_2 - 2\ga_2 e_n  \sim
\\
\cdots  \sim  \ga^* + \sum_{j=1}^{n-1} 2\ga_je_j - (\sum_{j=1}^{n-1} 2 \ga_j) e_n.
\end{multline*}
Thus $\ga^*  \sim  \ga$, so that Lemma \ref{lemma_cCD} will be proved as soon as we have proved Lemma \ref{lemma_imply}.

We now prove Lemma \ref{lemma_imply} by induction. The statement clearly holds for the base case $\theta=0$. 
We make the inductive hypothesis that $\ga \sim  \ga+2(\theta-1) e_j - 2(\theta -1) e_n$ for some integer $\theta \geq 1$, and show that this implies 
$\ga \sim  \ga+2 \theta e_j - 2 \theta e_n$. Here we may naturally assume we are in the case $2\theta \leq 2\ga_n$. Without loss of generality we assume we are dealing with the case $j=1$, for notational simplicity.
Under the inductive hypothesis, $c_{\ga}$ satisfies (\ref{cDn}) if and only if $c_{\ga+2(\theta-1) e_1 - 2(\theta -1) e_n}$ also satisfies (\ref{cDn}). We need only show that  $c_{\ga+2(\theta-1) e_1 - 2(\theta -1) e_n}$ satisfies (\ref{cDn}) if and only if $c_{\ga + 2\theta e_1 - 2\theta e_n}$ satisfies (\ref{cDn}).

By (\ref{Bmej_relation}) applied with $\rho$ chosen to be $\ga +(2\theta-1)e_1 - (2 \theta-1)e_n$, we see that 
\begin{multline*}
 c_{\rho + e_n -e_1} (\rho_n+1) - c_{\rho - e_n + e_1} (\rho_1+1)  \\
=c_{\ga-2(\theta-1)e_n + 2(\theta-1)e_1}(2\ga_n - 2(\theta-1)) - c_{\ga -2\theta e_n + 2\theta e_1}(2\ga_1 + 2\theta)
=0.
\end{multline*}
Since $\ga \in \Lambda$ by assumption, we have $c_\ga \neq 0$, and so by the induction hypothesis, $c_{\ga-2(\theta-1)e_n + 2(\theta-1)e_1} \neq 0$.
Thus we may write 
\beq\label{relation_cc}
 \frac{c_{\ga + 2\theta e_1 - 2\theta e_n}}{c_{\ga+2(\theta -1) e_1 - 2(\theta  -1)e_n} } =  \frac{(\ga_n - (\theta-1))}{(\ga_1 + \theta)} .
 \eeq
It is now clear from (\ref{relation_cc}) that 
\[ c_{\ga + 2\theta e_1 - 2 \theta e_n} = C  \frac{k!}{(\ga_1+\theta)!\ga_2! \cdots \ga_{n-1}! (\ga_n - \theta)!} \]
if and only if 
\[ c_{\ga+2(\theta -1) e_1 - 2(\theta -1) e_n}  = C \frac{k!}{(\ga_1+(\theta-1))!\ga_2! \cdots \ga_{n-1}! (\ga_n -( \theta -1))!}.\]
We may conclude that $\ga + 2(\theta-1) e_1 - 2(\theta-1)e_n \sim \ga + 2\theta e_1 - 2\theta e_n$, which proves Lemma \ref{lemma_imply}.
This completes the proof of Statement \ref{item_m1_n} of Lemma \ref{lemma_B_n}.

\subsection{Proof of Lemma \ref{lemma_B_n} for $B_{m,\be}$, Statement \ref{item_m2_n}}\label{sec_B_lemma_3}
We will prove that if $m \in \Lambda$ and $p_m(y)$ is parabolic, then there exists some $\be$ with $|\be|=2$ such that $B_{m,\be}$ is not the zero polynomial. 
In fact, it suffices to consider $\be=2e_j$ for any fixed $1 \leq j \leq n-1$; to fix ideas we choose $j=1$. 
We recall the definition 
\[ B_{m,2e_1}(w) = \sum_{|\al|=m} c_\al A_{\al,2e_1}(w).\]
By Statement \ref{item_A_2} of Lemma \ref{lemma_A_ab}, for each $\al$ with $|\al|=m$, 
\[ A_{\al, 2e_1}(w) = {\al_n \choose 2} w^{\al + 2e_1 - 2e_n} - \al_n \al_1 w^\al + {\al_1\choose 2} w^{\al -2e_1 + 2e_n},\]
with the understanding that any term with a multi-index  with a negative entry does not actually appear.
Thus after regrouping terms,
\[B_{m,2e_1} (w) = \sum_{|\ga|=m} \left\{ c_{\ga-2e_1+2e_n} { \ga_n +2 \choose 2} -c_\ga \ga_n \ga_1 + c_{\ga + 2e_1-2e_n} {\ga _1 +2 \choose 2} \right\} w^\ga,\]
still with the understanding that any term involving $c_\al$ with a multi-index $\al$ with a negative entry does not appear.
In particular, we see that that coefficient of $w_1^m$, that is $w^{\ga}$ with $\ga = (m,0,\ldots,0)$, is precisely 
\beq\label{specific_c}
 c_{(m-2)e_1 + 2e_n}.
 \eeq
By assumption, $m \in \Lambda$ and $p_m$ is parabolic, so that $m$ is even and $c_{\al} \neq 0$ for all indices $\al$ such that $|\al|=2k$ and all entries in $\al$ are even. Thus we see that (\ref{specific_c}) is a nonzero constant, and thus in particular $B_{m,2e_1}$ is not the zero polynomial. This suffices for Statement \ref{item_m2_n}, and completes the proof of Lemma \ref{lemma_B_n}, and hence of Proposition \ref{prop_Ksharp_n}.

\section{Final treatment of Propositions \ref{prop_Case_1}, \ref{prop_Case_2a}, \ref{prop_Case_2b} in general dimension}\label{sec_gen_dim}
\subsection{Completing the proof of Proposition \ref{prop_Case_1}  (general $n \geq 2)$}\label{sec_high_dim_termI}
With Proposition \ref{prop_Ksharp_n} in hand,  we return briefly to the treatment of the term $\mathrm{\bf I}$ in Section \ref{sec_term_I} for all dimensions $n \geq 2$. 
We recall that this term is represented as in (\ref{I_dfn_n2_2}) by 
\begin{multline*}
\mathrm{\bf I}= \iint_{\R^{n+1}} e^{iP_{\nu}(u+z)-iP_{\mu}(z)} \eta(u+z) \eta(z) \zeta(\xi + |z|^2)
\\
\cdot \; \left( \uDel_{j-k}(\theta - |u+z|^2 + |z|^2) - \uDel_{j-k}(\theta - |u+z|^2 + |z|^2 + \xi) \right)  
dz d\xi .
\end{multline*}
We would again like to isolate an oscillatory integral within this term that is independent of the $\uDel_{j-k}$ factors. Regardless of the dimension, we are still motivated to define a new variable $\tau = (u \cdot z)/|u|$ in order to capture the behavior of $z$ in $|u+z|^2 - |z|^2$. We fix the partition of unity $\sum_l W_l(s)$ given in (\ref{partition_1})
and write 
\[ \mathrm{\bf I} = \sum_{l=1}^n \mathrm{\bf I}_l \]
where for each $1 \leq l \leq n$
\begin{multline*}
 \mathrm{\bf I}_l = W_l\left( \frac{u}{|u|} \right)  \iint_{\R^{n+1}} e^{iP_{\nu}(u+z)-iP_{\mu}(z)} \eta(u+z) \eta(z) \zeta(\xi + |z|^2)
\\
\cdot \left( \uDel_{j-k}(\theta - |u+z|^2 + |z|^2) - \uDel_{j-k}(\theta - |u+z|^2 + |z|^2 + \xi) \right)  
dz d\xi .
\end{multline*}
For each index $l$ we make the change of variables $z \mapsto (\tau,\sigma)$ relevant to the $l$-th coordinate as given in (\ref{tau_rel}) and (\ref{sig_rel}); under this transformation the term $\mathrm{\bf I}_l$ becomes
\begin{multline*}
 \mathrm{\bf I}_l = \left( \frac{|u_l|}{|u|} \right)^{n-2} W_l\left( \frac{u}{|u|} \right)  \int_{\R^{2}} K_{\sharp,l}^{\nu,\mu}( u,\tau;\xi) \left( \uDel_{j-k}(\theta - |u|^2 - 2|u| \tau) \right. \\
  - \; \left. \uDel_{j-k}(\theta - |u|^2 -2|u| \tau + \xi) \right)  
d\tau d\xi ,
\end{multline*}
where 
\[ K_{\sharp,l}^{\nu,\mu}(u,\tau;\xi) =  \int_{\R^{n-1}} e^{iP_{\nu}(u+z)-iP_{\mu}(z)} \eta(u+z) \eta(z) \zeta(\xi + |z|^2) d \sigma.\]
Here $z^{(l)} \in \R^{n-1}, z_l \in \R$ are implicitly defined in terms of $u \in \R^n$ and $\sig \in \R^{n-1},\tau \in \R$ by (\ref{z_rel1}) and (\ref{z_rel2}). 
As previously observed in (\ref{sig_upper_bound}), the range of integration for $\sig$ is in the compact set $|\sig| \leq c_0$, where $c_0$ is the absolute constant specified in (\ref{zeta_cond}). Moreover, $K_{\sharp,l}^{\nu,\mu}(u,\tau;\xi)$ has support where $u \in B_2( \R^n)$,  $|\xi| \leq 2$, $|\tau| \leq 1$. 

Applying the mean-value theorem to $\uDel_{j-k}$, as recorded in Lemma \ref{lemma_uDel_diff}, we have
$$
|\mathrm{\bf I}_l| \leq C 2^{-2(j-k)} \int_{\R^2} |K^{\nu,\mu}_{\sharp,l}(u, \tau; \xi)| \;
|\xi| \chi_{B_2}(\xi) \psi_{j-k} (\theta - |u|^2 - 2|u|\tau) \chi_{B_1}(\tau) d\tau d\xi ,
$$
where $\psi_{j-k}(t)$, as defined in Lemma \ref{lemma_uDel_diff}, is an $L^1$ dilation of $(1+t^2)^{-1}$, and is thus uniformly in $L^1$, independent of $j-k$.  
We now apply Proposition \ref{prop_Ksharp_n} to bound $K_{\sharp,l}^{\nu,\mu}$ and conclude that there exists $\delta > 0$ and a small set $G^\nu \subset B_2(\R^n)$ with 
$|G^{\nu}| \leq C r^{-\delta} $, and for each $u \in B_2(\R^n)$ a small set $F^\nu_u \subset B_1(\R)$ with $|F^{\nu}_u| \leq C r^{-\delta},$
such that
\beq\label{K_sharp_bd_app_n}
|K^{\nu,\mu}_{\sharp,l}(u, \tau; \xi)| \leq C \left[ r^{-\delta} \chi_{B_2}(u) \chi_{B_1}(\tau) + \chi_{G^{\nu}}(u) \chi_{B_1}(\tau) + \chi_{B_2}(u) \chi_{F_u^{\nu}}(\tau) \right].
\eeq
Moreover, these estimates are uniform in $\xi$ and the index $l$, as the small sets do not depend on $\xi$ or $l$, and neither do the bounds.
Hence we may sum over $1\leq l \leq n$ to obtain (for some universal constant $C$)
\begin{multline*}
 |\mathrm{\bf I}| 
\leq C2^{-2(j-k)} \int_{\R^2} \left( r^{-\delta} \chi_{B_2}(u) \chi_{B_1}(\tau) + \chi_{G^{\nu}}(u) \chi_{B_1}(\tau) + \chi_{B_2}(u) \chi_{F_u^{\nu}}(\tau) \right) \\
\cdot \;  \chi_{B_2}(\xi)
\psi_{j-k} (\theta - |u|^2 - 2|u|\tau) d\tau d\xi .
\end{multline*}
This is the exact analogue of (\ref{I_bound}) and we may proceed via the argument used to treat (\ref{I_bound}) in Section \ref{sec_term_I}. We conclude that the contribution of the term $\mathrm{\bf I}$ to the kernel of $TT^* f(x,t)$ leads to an operator with $L^2$ norm bounded above by $C2^{-2(j-k)}r^{-\del/2} ,$ as in (\ref{ITT}).

\subsection{Completing the proof of (\ref{case2aI}) and (\ref{case2bI}) for $I_a^\lam$ (general $n \geq 2$)}\label{sec_high_dim_2a2b}
We next briefly return to the treatment of the kernel $^{(1)}\mathcal{K}^{\nu,\mu}$ in Section \ref{sec_I_2}, now treating the case of general dimension $n \geq 2$. We need to prove  (\ref{case2aI}) and (\ref{case2bI}). 

Assume $j < k$.
We first prove (\ref{case2aI}).  Recall from (\ref{gen_kernel_11}) that the kernel relevant to $T_1T_1^*$ is  
\[ ^{(1)}\Kcal^{\nu,\mu} (u,\theta) = \int_{\R^{n}} e^{iP_\nu(u+z) - iP_\mu(z)} \eta(u+z) \eta(z) \uDel_{j-k}(\theta - |u+z|^2 + |z|^2) dz.\]
Again using the partition of unity (\ref{partition_1}), we write 
\beq\label{partition_K11}
 ^{(1)}\Kcal^{\nu,\mu}  (u,\theta) = \sum_{l=1}^n {}^{(1)}\mathcal{K}_l^{\nu,\mu} (u,\theta)
 \eeq
with 
\[ ^{(1)}\Kcal^{\nu,\mu}_l (u,\theta) = W_l\left( \frac{u_l}{|u|} \right) \int_{\R^{n}} e^{iP_\nu(u+z) - iP_\mu(z)} \eta(u+z) \eta(z) \uDel_{j-k}(\theta - |u+z|^2 + |z|^2) dz.\]
For each index $l$ we make the change of variables $z \mapsto (\tau,\sig)$ defined with respect to the $l$-th coordinate in (\ref{tau_rel}) and (\ref{sig_rel}), so that 
\beq\label{partition_K11_sharp}
 ^{(1)}\Kcal^{\nu,\mu}_l (u,\theta) =  \left( \frac{|u_l|}{|u|} \right)^{n-2} W_l\left( \frac{u}{|u|} \right)  \int_{\R} K_{\sharp,l}^{\nu,\mu} (u,\tau)\uDel_{j-k}( \theta - |u|^2 - 2|u| \tau) d\tau,
\eeq
where 
\[ K_{\sharp,l}^{\nu,\mu}(u,\tau) = \int_{\R^{n-1}} e^{i P_\nu(u+z) - iP_\mu(z) } \eta (u+z) \eta(z) d\sig.\]
We apply the nontrivial bound of Proposition \ref{prop_Ksharp_n} to $K_{\sharp,l}^{\nu,\mu}$ and the trivial bound to $\uDel_{j-k}$, to conclude that the analogue of (\ref{K_11_0}) holds, uniformly in $l$. From here the analysis of $T_1T_1^*$ 
proceeds as in (\ref{T1T1_0}), and we may conclude that this portion of the operator has norm bounded by $Cr^{-\del/2}$.

We next prove (\ref{case2bI}). 
We again partition $^{(1)}\Kcal^{\nu,\mu}_l (u,\theta)$ as in (\ref{partition_K11}) with components given by (\ref{partition_K11_sharp}) after the appropriate change of variables.
We then use the identity (\ref{delta_expression_0}) for $\underline{\Delta}_{j-k}$ as before, so that we may write each component as
\[  ^{(1)}\Kcal_l^{\nu,\mu} (u,\theta) = \frac{2^{2(j-k)}}{2|u|}  \left( \frac{|u_l|}{|u|} \right)^{n-2} W_l\left( \frac{u}{|u|} \right)  \int_{\mathbb{R}}  \partial_\tau K_{\sharp,l}^{\nu,\mu} (u,\tau)\widetilde{\uDel}_{j-k}(\theta -  2|u| \tau - |u|^2) d\tau,\]
with the Schwartz function $\widetilde{\uDel}_{j-k}$ constructed in Lemma \ref{lemma_Del_anti}.
We now note that since $K_{\sharp,l}^{\nu,\mu}(u,\tau)$ is supported where $|u|\leq 2$, $|\tau| \leq 1$ and is a smooth function of $\tau$, 
\[| \partial_\tau K_{\sharp,l}^{\nu,\mu}(u,\tau)| \leq c_0^{-(n-1)}r\chi_{B_2}(u) \chi_{B_1}(\tau), \]
where $c_0$ is the absolute constant (\ref{zeta_cond}) coming from the restriction $|\sigma|\leq c_0^{-1}$, and the factor of $r$ comes from bringing down coefficients of size $\| \nu \|, \| \mu \| \approx r$ when differentiating the phase $P_\nu(u+z) - P_\mu(z)$ with respect to $\tau$. 
Hence
\beq\label{K11_l_n}
 |^{(1)}\Kcal_l^{\nu,\mu} (u,\theta) | \leq \frac{c_0^{-(n-1)} r 2^{2(j-k)}}{2|u|}\chi_{B_2}(u) \int \chi_{B_1}(\tau) \widetilde{\uDel}_{j-k}(\theta - 2|u|\tau-|u|^2) d\tau,
 \eeq
uniformly in $l$. By the uniformity in $l$ we see that we may sum the bounds provided by (\ref{K11_l_n}) over $1 \leq l \leq n$ so that a bound of the order (\ref{K11_l_n}) holds for the full kernel $^{(1)}\Kcal^{\nu,\mu}$. Thus from here on the analysis of $T_1T_1^*$ may proceed as from (\ref{T1T1_f_2}) to (\ref{T1T1_final_2}). 
This completes the proof.

\section{Appendix: Proof of Theorem \ref{thm_trunc_Radon}}\label{sec_appendix}

Let $\ga(t) = (t,|t|^2) \subset \R^{n+1}$ denote a parametrization of the paraboloid and define 
 $\mathcal{H}_{\lambda}$ acting on Schwartz functions $f$ on $\R^{n+1}$ by
\beq\label{H_lower_lam}
\mathcal{H}_{\lambda} f(x) = \int_{|t| > \lambda} f(x-\gamma(t)) K(t) dt ,
\eeq
where $K$ is a Calder\'{o}n-Zygmund kernel.
Our goal is to prove Theorem  \ref{thm_trunc_Radon}, which states that the  maximal truncated singular Radon transform defined by
$$\mathcal{H}^* f (x) = \sup_{\lambda > 0} |\mathcal{H}_{\lambda} f(x)|$$
is bounded on $L^p$ for all $1<p<\infty$. The following simple proof is based on the approach of \cite{DuoRub86}.

We recall the non-truncated singular integral operator $\Hcal$ along the paraboloid defined in (\ref{Hcal_dfn}) and the maximal Radon transform $\Mcal_{\Rad}$ along the paraboloid, defined in (\ref{MRad_dfn}), both of which are known to be bounded on $L^p$ for $1<p<\infty$.
We first make the simple observation that if $\lam, \lam' >0$ with $\lam \leq \lam' \leq 2\lam,$ then
\[ |\Hcal_{\lam'} f | \lesssim |\Hcal_{\lam} f| + \Mcal_\Rad f,\]
so that it suffices to prove that the operator 
\beq\label{H2k}
 f \mapsto \sup_{k \in \Z} | \Hcal_{2^k} f | 
 \eeq
is bounded on $L^p$ for $1<p < \infty$.

In preparation for the proof of this, observe that the Fourier multiplier $m(\xi)$ of $\Hcal$ is given by 
\[ m(\xi) = \int_{\R^n} K(t) e^{-2\pi i \ga(t) \cdot \xi} dt.\]
Similarly the Fourier multiplier $m_\lam(\xi)$ of $\Hcal_\lam$ is given by
\beq\label{m_trunc} 
m_\lam(\xi) = \int_{|t| >\lam} K(t) e^{-2\pi i \ga(t) \cdot \xi} dt.
\eeq
It will be convenient to define a non-isotropic dilation of $\xi  = (\xi',\xi_{n+1}) \in \mathbb{R}^{n+1}$  by 
\[ \lambda \circ \xi = (\lambda \xi', \lambda^2 \xi_{n+1}).\]
Accordingly, we denote the non-isotropic norm of $\xi$ by  $\|\xi\| = |\xi'| + |\xi_{n+1}|^{1/2}$, and the non-isotropic ball of radius $r$ by 
\[\Ball_r = \{y: \| y \| \leq r\}.\]
We may also define a maximal function averaging over non-isotropic balls by
\[ \Mcal_\Ball f(x) =  \sup_{a>0} \int f(x-y) \frac{1}{a^{n+2}} \chi_{\Ball_1}(a^{-1} \circ y) dy;\]
this is known to be a bounded operator on $L^p$ for $1 < p \leq \infty$ (see Chapter I of \cite{SteinHA}).

Now let $\eta \in C_c^\infty (\R^{n+1})$ be a smooth bump function that is identically one for $\| \xi \|\leq 1/2$ and vanishes for $\| \xi \| \geq 1$.
Define $\Rcal_\lam$ to be the operator with Fourier multiplier $(1-\eta (\lam \circ \xi)) m_\lam(\xi)$; this may be compared pointwise to $\Hcal_\lam$ as follows (we defer the proof for the moment):  
\begin{lemma}\label{HRMR}
\[ |\Hcal_\lam f| \lesssim |\Rcal_\lam f| + \Mcal_\Ball (\Hcal f) + \Mcal_\Ball f.\]
\end{lemma}
As the last two terms are known to be bounded on $L^p$ for $1<p<\infty$, the proof  of (\ref{H2k})  reduces to showing that 
\beq\label{R_theorem}
 f \mapsto \sup_{k \in \Z} | \Rcal_{2^k} f |  \quad \text{is bounded on $L^p$ for $1<p < \infty$.} 
\eeq

For each $l \in \Z$ we set 
\[ K^{(l)}(u) = \frac{1}{2^{ln}} K(\frac{u}{2^l});\]
as is well-known, each such kernel is also a Calder\'{o}n-Zygmund kernel satisfying the same bounds as $K$, with constants uniform in $l$. We now define 
\beq\label{tilde_m_dfn}
 \widetilde{m}^{(l)} (\xi) = \int_{1 \leq |t| \leq 2} K^{(l)}(t) e^{-2\pi i \ga(t) \cdot \xi} dt.
 \eeq
Then for each $k \in \Z$ we may decompose the multiplier $m_{2^k}$ defined in (\ref{m_trunc}) as
\[ m_{2^k}(\xi) = \sum_{j=0}^\infty  \int_{2^{j+k} \leq |t| \leq 2^{j+k+1}} K(t) e^{-2\pi i \ga (t) \cdot \xi} dt = \sum_{j=0}^\infty \widetilde{m}^{(j+k)}(2^{j+k} \circ \xi).\]
Correspondingly,  $\Rcal_{2^k}$ has the Fourier multiplier
\[ \sum_{j=0}^\infty (1-\eta(2^k \circ \xi)) \widetilde{m}^{(j+k)}(2^{j+k} \circ \xi),\]
so that if we let $T_{j,k}$ denote the operator with multiplier $(1-\eta(2^k \circ \xi)) \widetilde{m}^{(j+k)}(2^{j+k} \circ \xi)$, we have
\[ \Rcal_{2^k} = \sum_{j=0}^\infty T_{j,k}.\]

The proof of (\ref{R_theorem}) will now follow quickly once we have the following two claims:
\begin{lemma}\label{lemma_Tjk1}
For each integer $j \geq 0$ and all $1<p<\infty$,
\[ \| \sup_{k \in \Z} |T_{j,k} | \|_{L^p} \leq C_p \|f \|_{L^p}.\]
\end{lemma}

\begin{lemma}\label{lemma_Tjk2}
For each integer $j \geq 0$,
\[\| \sup_{k \in \Z} |T_{j,k}| \|_{L^2} \leq  \| ( \sum_{k \in \Z} |T_{j,k} f|^2 )^{1/2} \|_{L^2} \leq C2^{-j/2} \|f \|_{L^2}. \]
\end{lemma}

Assuming these facts for the moment, we note that 
\beq\label{RTjk}
 \sup_{k \in \Z} |\Rcal_{2^k} f| \leq \sum_{j=0}^\infty \sup_{k \in \Z} |T_{j,k} f|.
 \eeq
But interpolating the results of Lemmas \ref{lemma_Tjk1} and \ref{lemma_Tjk2}, we see that for all $1<p<\infty$ there exists some $\ep = \ep(p)>0$ such that 
\[ \| \sup_{k \in \Z} |T_{j,k}| \|_{L^p}  \leq C 2^{-j\ep} \| f\|_{L^p}.\]
Taking norms and adding these estimates in (\ref{RTjk}), we see that $ \|\sup_{k \in \Z} |\Rcal_{2^k} f | \|_{L^p}$
is finite  for any $1<p<\infty$, as desired. This completes the proof of Theorem \ref{thm_trunc_Radon}, aside from the proofs of the three small lemmas, which we treat in the next section.

\subsection{Proof of the lemmas}

To prove Lemma \ref{HRMR}, we note that certainly 
\beq\label{three-terms}
m_\lam(\xi) = (1-\eta(\lam \circ \xi)) m_\lam(\xi) + \eta (\lam \circ \xi) m(\xi) + \eta (\lam \circ \xi) (m_\lam(\xi) - m(\xi)).
\eeq
The first term on the right hand side corresponds to $\Rcal_\lam$ by definition. The second term is the Fourier multiplier of the operator 
\[  f\mapsto \int \Hcal f(x-y) \lam^{-(n+2)}\check{\eta}(\lam^{-1} \circ y) dy,\]
which may be majorized by $\Mcal_{\Ball} (\Hcal f)$,  with a constant independent of $\lam$.
To handle the third term, we first note that we can write
\[  m_\lam(\xi) - m(\xi)= \int_{|t| \leq \lam} K(t) e^{-2\pi i \ga(t) \cdot \xi} dt = \int_{|t| \leq 1} K^{(\lam^{-1})}(t) e^{-2\pi i \ga(t) \cdot (\lam \circ \xi)} dt,\]
which we denote by $m_1^{(\lam^{-1})}(\lam \circ \xi).$ (Here we are slightly abusing notation by letting $K^{(\rho)}$ denote $\rho^{-n} K(\rho^{-1} \cdot)$ for any parameter $\rho>0$, not necessarily dyadic.)
Thus the third term in (\ref{three-terms}) is the dilation of the function 
\beq\label{eta_m}
\eta(\xi) m_1^{(\lam^{-1})}(\xi),
\eeq
where  $m_1^{(\lam^{-1})}(\xi)$ is the Fourier transform of the compactly supported distribution defined (up to sign) for $(t,s) \in \R^{n+1}$ by
 \[ \chi_{|t| \leq 1} (t) K^{(\lam^{-1})}(t) \del_{s = |t|^2}.\]
 The Fourier transform of this compactly supported distribution is a $C^\infty$ function, and due to the uniformity of $K^{(\lam^{-1})}$ in $\lam$, the bound for the Fourier transform and all its derivative holds uniformly in $\lam$. This verifies that (\ref{eta_m}) is a compactly supported $C^\infty$ function (uniformly in $\lam$), and hence we see that its inverse Fourier transform may be bounded on $\R^{n+1}$ by $C_N (1+\|x\|)^{-N}$ for any $N$ (uniformly in $\lam$). Thus upon dilating the last term in (\ref{three-terms}) by $\lam$,  multiplying with $\hat{f}(\xi)$,  and taking the inverse Fourier transform, we see that this term may be bounded by $\Mcal_{\text{Ball}} f$, as claimed.

To prove Lemma \ref{lemma_Tjk1}, we let $S_j$ denote the operator with multiplier $\widetilde{m}^{(j)}(2^j \circ \xi)$, that is,
\[ S_j f(x) = \int_{2^j \leq |t| \leq 2^{j+1}} f(x-\ga(t)) K(t)dt.\]
In particular, we may see immediately from the upper bound (\ref{K_prop}) for $K$ that 
\[ |S_jf(x) |\lesssim \frac{1}{2^{nj}} \int_{2^j \leq |t| \leq 2^{j+1}} |f(x-\ga(t))| dt\lesssim \Mcal_{\Rad} f,\]
uniformly in $j$.
In addition, by adapting the argument given above for Lemma \ref{HRMR}, we see that uniformly in $j,k$,
\[ |T_{j,k} f| \leq |S_{j+k} f| + \Mcal_\Ball (S_{j+k} f) \lesssim \Mcal_{\Rad} f + \Mcal_\Ball (\Mcal_{\Rad} f).\]
Thus we immediately obtain that for each $j \geq 0$ and any $1<p<\infty$,
\[ \| \sup_{k \in \Z} |T_{j,k} f| \||_{L^p}  \leq C_p \| f\|_{L^p}.\]

To prove Lemma \ref{lemma_Tjk2}, it suffices by Plancherel's theorem to bound 
\[ \sum_{k \in \Z} | (1 -\eta(2^k \circ \xi)) \widetilde{m}^{(j+k)}(2^{j+k} \circ \xi) |^2 \]
in $L^\infty$. 
For this we will use the following van der Corput estimate:
\beq\label{Vdc_m}
\widetilde{m}^{(j)}(\xi) \leq \frac{C}{\| \xi \|^{1/2}}. 
\eeq
With this in hand, and applying the vanishing of $1-\eta$ near the origin, we conclude that 
\begin{eqnarray*}
 \sum_{k \in \Z} | (1 -\eta(2^k \circ \xi)) \widetilde{m}^{(j+k)}(2^{j+k} \circ \xi) |^2 
	&\leq &\sum_{2^k \geq \| \xi \|^{-1}} | \widetilde{m}^{(j+k)}(2^{j+k} \circ \xi) |^2 \\
		&\leq&  \sum_{2^k \geq \| \xi \|^{-1}} \frac{C}{2^{j+k} \| \xi \|} = C2^{-j},
		\end{eqnarray*}
		which is sufficient for Lemma \ref{lemma_Tjk2}.
Finally, to prove (\ref{Vdc_m}) we recall that 
\[ \widetilde{m}^{(j)}(\xi) = \int_{1 \leq |t| \leq 2} K^{(j)}(t) e^{-2\pi i (t \cdot \xi' + |t|^2 \xi_{n+1})} dt.\]
If $|\xi'| \geq 8 n |\xi_{n+1}|,$ then the first derivative test gives 
\[ |\widetilde{m}^{(j)}(\xi)| \leq  |\xi'|^{-1},\]
uniformly in $j$.
For indeed, without loss of generality we may assume that 
\[|\xi'_1| \geq \frac{1}{n} |\xi'| \geq 8 |\xi_{n+1}|.\]
 Then the first partial with respect to $t_1$ of the phase is $-2\pi i(\xi'_1 + 2 t_1 \xi_{n+1}),$ which is bounded below in absolute value by $2\pi \cdot \frac{1}{2}|\xi'_1|$ if $|\xi'_1| \geq 8 |\xi_{n+1}|$ and $1 \leq |t| \leq 2$. Writing 
$$e^{-2\pi i(t \cdot \xi' + |t|^2 \xi_{n+1})}=\frac{\partial_{t_1} (e^{-2\pi i(t \cdot \xi' + |t|^2 \xi_{n+1})})}{-2\pi i(\xi'_1 + 2 t_1 \xi_{n+1})}$$
 and integrating by parts once, our claim follows; note that the boundary terms also obey this bound.
On the other hand, if $|\xi'| \leq 8 n |\xi_{n+1}|$ then the second derivative test applied to the integral in $t_1$ gives 
\[  |\widetilde{m}^{(j)}(\xi)|\leq  |\xi_{n+1}|^{-1/2};\]
see for example the Corollary of Proposition 2 of Chapter VIII in \cite{SteinHA}.
Combining these bounds proves (\ref{Vdc_m}).

\subsection*{Acknowledgements}
The authors would like to thank Elias M. Stein for suggesting this problem and for a number of helpful conversations---as well as for many years of generous mentoring. The authors also thank Petru Mironescu and Jean van Schaftingen for a helpful discussion on stopping-times. Pierce is partially supported by NSF DMS-1402121, and during this collaboration has been partially supported by NSF DMS-0902658, DMS-0635607, the Simonyi Fund at the Institute for 
Advanced Study and a Marie Curie Incoming International Fellowship at Oxford University, with additional support in Oxford via an EPSRC Developing Leaders  grant. Yung was partially supported by NSF grant DMS-1201474, a Titchmarsh Fellowship at Oxford University, a Junior Research Fellowship from St. Hilda's College, and a direct grant for research from the Chinese University of Hong Kong (4053120). The authors thank the Institute for Advanced Study, Princeton University, the University of Oxford, and the Hausdorff Center for Mathematics for providing pleasant working environments during a number of collaborative visits in 2010 -- 2015.

\bibliographystyle{amsplain}
\bibliography{AnalysisBibliography}
\label{endofproposal}
\end{document}

%% file: format.tex
\newtheorem{letterthm}{Theorem}

\newtheorem{thm}{Theorem}[section]
\newtheorem{prop}[thm]{Proposition}
\newtheorem{lemma}[thm]{Lemma}
\newtheorem{cor}[thm]{Corollary}

\newcommand{\con}{\equiv}

\newcommand{\bstack}[2]{#1 \atop #2}

\newcommand{\maps}{\rightarrow}

\newcommand{\cH}{\mathcal{H}}

\newcommand{\al}{\alpha}
\newcommand{\be}{\beta}

\newcommand{\ga}{\gamma}
\newcommand{\del}{\delta}
\newcommand{\ep}{\varepsilon}

\newcommand{\om}{\omega}
\newcommand{\Om}{\Omega}

\newcommand{\sig}{\sigma}
\newcommand{\lam}{\lambda}

\newcommand{\Ga}{\Gamma}

\newcommand{\Hcal}{\mathcal{H}}

\newcommand{\Lcal}{\mathcal{L}}
\newcommand{\Mcal}{\mathcal{M}}

\newcommand{\Scal}{\mathcal{S}}

\newcommand{\N}{\mathbb{N}}

\newcommand{\R}{\mathbb{R}}

\newcommand{\Z}{\mathbb{Z}}

\newcommand{\beq}{\begin{equation}}
\newcommand{\eeq}{\end{equation}}